\documentclass[3p,times]{elsarticle}
\usepackage[T1]{fontenc}
\usepackage[utf8]{inputenc}
\usepackage[british]{babel}
\usepackage[autostyle]{csquotes}
\usepackage[]{geometry} 
\usepackage{microtype}
\usepackage{amsmath}
\usepackage{amsfonts}
\usepackage{amssymb}
\usepackage[normalem]{ulem} 
\usepackage{mathtools} 
\usepackage[shortlabels]{enumitem}
\usepackage{tikz}
\usetikzlibrary{cd,arrows,graphs,petri,snakes,calc}
\usepackage[colorlinks=true,urlcolor=blue,linkcolor=blue,citecolor=red]{hyperref}
\usepackage{rotating}
\usepackage{graphicx}
\usepackage{ifthen}
\usepackage{proof}
\usepackage{subdepth} 
\usepackage{xparse} 

\usepackage{JPAAsettings}

\makeatletter
\def\ps@pprintTitle{%
	\let\@oddhead\@empty
	\let\@evenhead\@empty
	\def\@oddfoot{\centerline{\thepage}}%
	\let\@evenfoot\@oddfoot}
\makeatother

\begin{document}

	\begin{frontmatter}
	\title{Composing Dinatural Transformations: Towards a Calculus of Substitution\tnoteref{t1}}
	\tnotetext[t1]{This is the \emph{accepted manuscript} version of a published journal article accessible at \url{https://doi.org/10.1016/j.jpaa.2021.106689}.}	
	\author{Guy McCusker}
	\ead{G.A.McCusker@bath.ac.uk}
	
	\author{Alessio Santamaria\corref{cor1}\fnref{fn1}}
	\ead{alessio.santamaria@di.unipi.it}

	\cortext[cor1]{Corresponding author}
    
    \fntext[fn1]{Present address: Dipartimento di Informatica, Università degli Studi di Pisa, Largo B.\ Pontecorvo 3, 56127 Pisa, Italy}
	
	\address{Department of Computer Science, University of Bath, BA2 7AY Bath, United Kingdom}
	
	\date{}
	
	\begin{abstract}

Dinatural transformations, which generalise the ubiquitous natural transformations to the case where the domain and codomain functors are of mixed variance, fail to compose in general; this has been known since they were discovered by Dubuc and Street in 1970. Many ad hoc solutions to this remarkable shortcoming have been found, but a general theory of compositionality was missing until Petri\'c, in 2003, introduced the concept of g-dinatural transformations, that is, dinatural transformations together with an appropriate graph: he showed how acyclicity of the composite graph of two arbitrary dinatural transformations is a sufficient and essentially necessary condition for the composite transformation to be in turn dinatural. Here we propose an alternative, semantic rather than syntactic, proof of Petri\'c's theorem, which the authors independently rediscovered with no knowledge of its prior existence; we then use it to define a generalised functor category, whose objects are functors of mixed variance in many variables, and whose morphisms are transformations that happen to be dinatural only in some of their variables.

We also define a notion of horizontal composition for dinatural transformations, extending the well-known version for natural transformations, and prove it is associative and unitary. Horizontal composition embodies substitution of functors into transformations and vice-versa, and is intuitively reflected from the string-diagram point of view by substitution of graphs into graphs.

This work represents the first, fundamental steps towards a substitution calculus for dinatural transformations as sought originally by Kelly, with the intention then to apply it to describe coherence problems abstractly. There are still fundamental difficulties that are yet to be overcome in order to achieve such a calculus, and these will be the subject of future work; however, our contribution places us well in track on the path traced by Kelly towards a calculus of substitution for dinatural transformations.
	\end{abstract}

	
	\end{frontmatter}

\section{Introduction}

The problem of coherence for a certain theory (like monoidal, monoidal closed\dots) consists in understanding which diagrams necessarily commute as a consequence of the axioms. One of the most famous results is Mac Lane's theorem on coherence for monoidal categories~\cite{mac_lane_natural_1963}: every diagram built up only using associators and unitors, which are the data that come with the definition of monoidal category, commutes. One of the consequences of this fact is that every monoidal category is monoidally equivalent to a \emph{strict} monoidal category, where associators and unitors are, in fact, identities. What this tells us is that those operations that one would like to regard as not important--such as the associators and unitors etc.--really are not important. Solving the coherence problem for a theory, therefore, is fundamental to the complete understanding of the theory itself.

In this article we aim to set down the foundations for the answer to an open question left by Kelly in his task to study the coherence problem abstractly, started with~\cite{kelly_abstract_1972,kelly_many-variable_1972}.  
Kelly argued that coherence problems are concerned with categories carrying an extra structure: a collection of functors and natural transformations subject to various equational axioms. For example, in a monoidal category $\A$ we have $\otimes \colon \A^2\! \to \A$, $I \colon \A^0 \to \A$; if $\A$ is also closed then we would have a functor of mixed variance $(-) \implies (-) \colon \Op\A\! \times \A \to \A$. The natural transformations that are part of the data, like associativity in the monoidal case:
\[
\alpha_{A,B,C} \colon (A \otimes B) \otimes C \to A \otimes (B \otimes C),
\]
connect not the basic functors directly, but rather functors obtained from them by \emph{iterated substitution}. By ``substitution'' we mean the process where, given functors
\[
K \colon \A \times \Op\B \times \C \to \D, \quad F \colon \E \times \mathbb G \to \A, \quad G \colon \mathbb H \times \Op\L \to \B, \quad H \colon \Op\M \to \C
\]
we obtain the new functor
\[
K(F,\Op G, H) \colon \E \times \mathbb G \times \Op{\mathbb H} \times \L \times \Op\M \to \D\label{substitution functors example}
\]
sending $(A,B,C,D,E)$ to $K(F(A,B),\Op G (C,D), H(E))$. Hence substitution generalises composition of functors, to which it reduces if we only consider one-variable functors. In the same way, the equational axioms for the structure, like the pentagonal axiom for monoidal categories:
\[
\begin{tikzcd}[column sep={3.5em,between origins},row sep=2em]
& & & (A \otimes B) \otimes (C \otimes D) \ar[drr,"\alpha_{A,B,C\otimes D}"] \\
\bigl( (A \otimes B) \otimes C \bigr) \otimes D \ar[urrr,"\alpha_{A\otimes B, C, D}"] \ar[dr,"\alpha_{A,B,C} \otimes D"'] & & & & & A \otimes \bigl( B \otimes (C \otimes D) \bigr) \\
& \bigl( A \otimes (B \otimes C) \bigr) \otimes D \ar[rrr,"\alpha_{A,B \otimes C, D}"'] & & & A \otimes \bigl( (B \otimes C) \otimes D \bigr) \ar[ur,"A \otimes \alpha_{B,C,D}"']
\end{tikzcd}
\]
involve natural transformations obtained from the basic ones by ``substituting functors into them and them into functors'', like $\alpha_{A \otimes B, C, D}$ and $\alpha_{A,B,C} \otimes D$ above. 

By substitution of functors into transformations and transformations into functors we mean therefore a generalised \emph{whiskering} operation or, more broadly, a generalised \emph{horizontal composition} of transformations. For these reasons Kelly argued in~\cite{kelly_many-variable_1972} that an abstract theory of coherence requires ``a tidy calculus of substitution'' for functors of many variables and appropriately general kinds of natural transformations, generalising the usual Godement calculus~\cite[Appendice]{godement_topologie_1958} for ordinary functors in one variable and ordinary natural transformations. (The ``five rules of the functorial calculus'' set down by Godement are in fact equivalent to saying that sequential composition of functors and vertical and horizontal composition of natural transformations are associative, unitary and satisfy the usual interchange law; see~\cite[Introduction]{santamaria_towards_2019} for more details.)

One could ask why bother introducing the notion of substitution, given that it is not primitive, as the functor  $K(F,\Op G, H)$ above can be easily seen to be the usual composite $K \circ (F \times \Op G \times H)$. Kelly's argument is that there is \emph{no need} to consider functors whose codomain is a product of categories, like $F \times \Op G \times H$, or the twisting functor $T(A,B) = (B,A)$, or the diagonal functor $\Delta \colon \A \to \A \times \A$ given by $\Delta(A) = (A,A)$, if we consider substitution as an operation on its own. However, take a Cartesian closed category $\A$, and consider the diagonal transformation $\delta_A \colon A \to A \times A$, the symmetry $\gamma_{A,B} \colon A \times B \to B \times A$ and the evaluation transformation $\eval A B \colon A \times (A \implies B) \to B$. It is true that we can see $\delta$ and $\gamma$ as transformations $\id \A \to \Delta$ and $\times \to \times \circ T$, but there is no way to involve $\Delta$ into the codomain of $\eval{}{}$, given that the variable $A$ appears covariantly and contravariantly at once. Kelly suggested adapting the notion of \emph{graph} for \emph{extranatural} transformations that he had introduced with Eilenberg~\cite{eilenberg_generalization_1966} to handle the case of natural transformations; that is, he proposed to consider natural transformations $\phi \colon F \to G$ between functors of many variables together with a graph $\Gamma(\phi)$ that tells us which arguments of $F$ and $G$ are to be equated when we write down the general component of $\phi$. The information carried by the graph is what allows us to get by without explicit mention of functors like $T$ and $\Delta$ and, moreover, it paves the way to the substitution calculus he sought.

With the notion of ``graph of a natural transformation'', Kelly constructed a full Godement calculus for covariant functors only. His starting point was the observation that the usual Godement calculus essentially asserts that $\Cat$ is a 2-category, but this is saying less than saying that $\Cat$ is actually \emph{Cartesian closed}, $- \times \B$ having a right adjoint $[\B,-]$ where $[\B,\C]$ is the functor category. Since every Cartesian closed category is enriched over itself, we have that $\Cat$ is a $\Cat$-category, which is just another way to say 2-category. Now, vertical composition of natural transformations is embodied in $[\B,\C]$, but sequential composition of functors and horizontal composition of natural transformations are embodied in the functor
\[
M \colon [\B,\C] \times [\A,\B] \to [\A,\C]
\]
given by the closed structure (using the adjunction and the evaluation map twice). What Kelly does, therefore, is to create a generalised functor category $\FC \B \C $ over a category of graphs $\Per$ and to show that the functor $\FC - -$ is the internal-hom of $\catover\Per$, which is then monoidal closed (in fact, far from being Cartesian or even symmetric), the left adjoint of $\FC \B -$ being denoted as $\ring - \B$. The analogue of the $M$ above, now of the form $\ring {\FC \B \C} {\FC \A \B} \to \FC \A \C,$ is what provides the desired substitution calculus.

When trying to deal with the mixed-variance case, however, Kelly ran into problems. He considered the every-variable-twice extranatural transformations of~\cite{eilenberg_generalization_1966} and, although he got ``tantalizingly close'', to use his words, to a sensible calculus, he could not find a way to define a category of graphs that can handle cycles in a proper way. This is the reason for the ``I'' in the title \emph{Many-Variable Functorial Calculus, I} of~\cite{kelly_many-variable_1972}: he hoped to solve these issues in a future paper, which sadly has never seen the light of day.

What we do in this article is, in fact, consider transformations between mixed-variance functors whose type is even more general than Eilenberg and Kelly's, corresponding to $\text{\uuline{G}}^*$ in~\cite{kelly_many-variable_1972}, recognising that they are a straightforward generalisation of \emph{dinatural transformations}~\cite{dubuc_dinatural_1970} in many variables. This poses an immediate, major obstacle: dinatural transformations notoriously fail to compose, as already observed by Dubuc and Street when they introduced them in 1970. There are certain conditions, known already to their discoverers, under which two dinatural transformations $\phi$ and $\psi$ compose: if either of them is natural, or if a certain square happens to be a pullback or a pushout, then the composite $\psi\circ\phi$ turns out to be dinatural. However, these are far from being satisfactory solutions for the compositionality problem, for either they are too restrictive (as in the first case), or they speak of properties enjoyed not by $\phi$ and $\psi$ themselves, but rather by other structures, namely one of the functors involved.  Many studies have been conducted about them~\cite{bainbridge_functorial_1990,blute_linear_1993,freyd_dinaturality_1992,girard_normal_1992,lataillade_dinatural_2009,mulry_categorical_1990,pare_dinatural_1998,pistone_dinaturality_2017,plotkin_logic_1993,simpson_characterisation_1993,wadler_theorems_1989}, and many attempts have been made to find a proper calculus for dinatural transformations, but until recently only \emph{ad hoc} solutions have been found and, ultimately, they have remained poorly understood.

In 2003, Petri\'c~\cite{petric_g-dinaturality_2003} studied coherence results for bicartesian closed categories, and found himself in need, much like Kelly in his more general case, of understanding the compositionality properties of \emph{g-dinatural transformations}, which are slightly more general dinatural transformations than those of Dubuc and Street~\cite{dubuc_dinatural_1970} in what their domain and codomain functors are allowed to have different variance and, moreover, they always come with a graph (whence the ``g'' in ``g-dinatural'') which reflects their signature. Petri\'c successfully managed to find a sufficient and essentially necessary condition for two consecutive g-dinatural transformations $\phi$ and $\psi$ to compose: if the composite graph, obtained by appropriately ``glueing'' together the graphs of $\phi$ and $\psi$, is acyclic, then $\psi\circ\phi$ is again g-dinatural. This result, which effectively solves the compositionality problem of dinatural transformations, 
surprisingly does not appear to be well known: 
fifteen years after Petri\'c's paper, the authors of the present article, completely oblivious to Petri\'c's contribution, independently re-discovered the same theorem, which was one of the results of~\cite{mccusker_compositionality_2018} and of the second author's PhD thesis~\cite{santamaria_towards_2019}\footnotemark. We, too, associated to each dinatural transformation a graph, inspired by Kelly's work of~\cite{kelly_many-variable_1972}, such graph being slightly different from Petri\'c's; we also proved that acyclicity of the composite graph of $\phi$ and $\psi$ is ``essentially enough'' for $\psi\circ\phi$ to be dinatural. The proof of our and Petri\'c's theorem are, deep down, following the same argument, but the main difference is in the approach we took to formalise it: Petri\'c's went purely syntactic, using re-writing rules to show how the arbitrary morphism of the universal quantification of the dinaturality property for $\psi\circ\phi$ can ``travel through the composite graph'' when the graph is acyclic, whereas we showed this by interpreting the composite graph as a \emph{Petri Net}~\cite{petri_kommunikation_1962} and re-casting the dinaturality property of $\psi\circ\phi$ into a \emph{reachability} problem. We then proceeded to solve it by exploiting the general theory of Petri Nets: in other words, we took a more semantic approach.

\footnotetext{We also presented our result as novel in various occasions, including in a plenary talk at the Category Theory conference in Edinburgh in 2019, yet nobody redirected us to Petri\'c's paper, which we found by chance only in September 2019.}

Because of this appreciable difference of Petri\'c's and our proof of the compositionality result for dinatural transformations, we believe it is worth presenting in this paper our theorem despite the non-novelty of its statement; moreover, we give here a more direct proof for it than the one in~\cite{mccusker_compositionality_2018}: this is done in Section~\ref{section vertical compositionality}. In Section~\ref{chapter horizontal}, we define a working notion of horizontal composition, that we believe will play the role of substitution of dinaturals into dinaturals, precisely as horizontal composition of natural transformation does, as shown by Kelly in~\cite{kelly_many-variable_1972}. Next, we form a generalised functor category $\FC \B \C$ for these transformations (Definition~\ref{def: generalised functor category}). Finally, we prove that $\FC \B -$ has indeed a left adjoint $\ring - \B$, which gives us the definition of a category of formal substitutions $\ring \A \B$ generalising Kelly's one. Although the road paved by Kelly towards a substitution calculus for dinatural transformations still stretches a long way, our work sets the first steps in the right direction for a full understanding of the compositionality properties of dinaturals, which hopefully will be achieved soon.

\paragraph{Notations} 
$\N$ is the set of natural numbers, including 0, and we shall ambiguously write $n$ for both the natural number $n$ and the set $\{1,\dots,n\}$. We denote by $\I$ the category with one object and one morphism. Let $\alpha\in \List{\{+,-\}}$, $\length\alpha=n$, with $\length{-}$ denoting the length function (and also the cardinality of an ordinary finite set). We refer to the $i$-th element of $\alpha$ as $\alpha_i$. Given a category $\C$, if $n\ge 1$, then we define $\C^\alpha=\C^{\alpha_1} \times \dots \times \C^{\alpha_n}$, with $\C^+=\C$ and $\C^-=\Op\C$, otherwise $\C^\alpha=\I$. 

Composition of morphisms $f \colon A \to B$ and $g \colon B \to C$ will be denoted by $g\circ f$, $gf$ or also $f;g$. The identity morphism of an object $A$ will be denoted by $\id A$, $1_A$ (possibly without subscripts, if there is no risk of confusion), or $A$ itself. 
Given $A$, $B$ and $C$ objects of a category $\C$ with coproducts, and given $f \colon A \to C$ and $g \colon B \to C$, we denote by $[f,g] \colon A + B \to C$ the unique map granted by the universal property of $+$.

We use boldface capital letters $\bfA,\bfB\dots$ for tuples of objects, whose length will be specified in context. Say $\bfA=(A_1,\dots,A_n) \in \C^n$: we can see $\bfA$ as a function from the set $n$ to the objects of $\C$. If $\sigma \colon k \to n$ is a function of sets, the composite $\bfA \sigma$ is the tuple $(A_{\sigma 1}, \dots, A_{\sigma k})$. For $\bfB \in \C^n$ and $i \in \{1,\dots,n\}$, we denote by $\subst B X i$ the tuple obtained from $\bfB$ by replacing its $i$-th entry with $X$, and by $\subst B \cdot i$ the tuple obtained from $\bfB$ by removing its $i$-th entry altogether. In particular, the tuple $\subst A X i \sigma$ is equal to $(Y_1,\dots,Y_k)$ where
\[
Y_j = \begin{cases}
X & \sigma j=i \\
A_{\sigma j} & \sigma j \ne i
\end{cases}.
\]
Let $\alpha \in \List\{+,-\}$, $\bfA = (A_1,\dots,A_n)$, $\sigma \colon \length\alpha \to n$, $i \in \{1,\dots,n\}$. We shall write $\substMV A X Y i \sigma$ for the tuple $(Z_1,\dots,Z_{\length\alpha})$ where
\[
Z_j = \begin{cases}
X & \sigma j = i, \alpha_j = - \\
Y & \sigma j = i, \alpha_j = + \\
A_{\sigma j} & \sigma j \ne i
\end{cases}\label{not:A[X,Y/i]sigma}
\]
We shall also write $\subst B {\bfA} i$ for the tuple obtained from $\bf B$ by substituting $\bfA$ into its $i$-th entry. For example, if $\bfA = (A_1,\dots,A_n)$ and $\bfB = (B_1,\dots,B_m)$, we have
\[
\subst B {\bfA} i = (B_1,\dots, B_{i-1},A_1,\dots A_n, B_{i+1}, \dots B_m).
\]

If $F \colon \B^{\alpha} \to \C$ is a functor, we define $\funminplus F {A_i} {B_i} i {\length\alpha}$ to be the following object (if $A_i$, $B_i$ are objects) or morphism (if they are morphisms) of $\C$:
\[
\funminplus F {A_i} {B_i} i {\length\alpha}= F(X_1,\dots,X_{\length\alpha}) \text{ where } X_i =
\begin{cases}
A_i & \alpha_i = - \\
B_i & \alpha_i = +
\end{cases}
\]
If $A_i = A$ and $B_i = B$ for all $i \in \length\alpha$, then we will simply write $\funminplusconst F A B$ for the above.

We denote by $\Not\alpha$ the list obtained from $\alpha$ by swapping the signs. Also, we call $\Op F \colon \B^{\Not\alpha} \to \Op\C$ the \emph{opposite functor}, which is the obvious functor that acts like $F$ between opposite categories.

\section{Vertical compositionality of dinatural transformations}\label{section vertical compositionality}

We begin by introducing the notion of \emph{transformation} between two functors of arbitrary variance and arity, which is simply a family of morphisms that does not have to satisfy any naturality condition. (This simple idea is, unsurprisingly, not new: it appears, for example, in~\cite{power_premonoidal_1997}.) A transformation comes equipped with a cospan in $\finset$ that tells us which variables of the functors involved are to be equated to each other in order to write down the general component of the family of morphisms. 

\begin{definition}\label{def:transformation}
	Let $\alpha$, $\beta \in \List\{+,-\}$, $F \colon \B^\alpha \to \C$, $G \colon \B^\beta \to \C$ be functors. A \emph{transformation} $\phi \colon F \to G$ \emph{of type} 
	$
	\begin{tikzcd}[cramped,sep=small]
	\length\alpha \ar[r,"\sigma"] & n & \length\beta \ar[l,"\tau"']
	\end{tikzcd}
	$ 
	(with $n$ a positive integer) is a family of morphisms in $\C$
	\[
	\bigl( \phi_{\bfA} \colon F(\bfA\sigma) \to G(\bfA\tau) \bigr)_{\bfA \in \B^n}
	\]
	(i.e., according to our notations, a family $\phi_{A_1,\dots,A_n} \colon F(A_{\sigma 1}, \dots, A_{\sigma\length\alpha}) \to G(A_{\tau1},\dots,A_{\tau\length\beta})$). Notice that $\sigma$ and $\tau$ need not be injective or surjective, so we may have repeated or unused variables.
	
	Given another transformation $\phi' \colon F' \to G'$ of type
	$
	\begin{tikzcd}[cramped,sep=small]
	\length\alpha \ar[r,"\sigma'"] & n & \length\beta \ar[l,"\tau'"']
	\end{tikzcd},
	$
	we say that
	\[
	\phi \sim {\phi'} \text{ if and only if there exists } \pi \colon n \to n \text{ permutation such that }
	\begin{cases}
	\sigma' = \pi\sigma \\
	\tau' = \pi\tau \\
	\phi'_\bfA = \phi_{\bfA\pi}
	\end{cases}. 
	\]
	$\sim$ so defined is an equivalence relation and we denote by $\class\phi$ the equivalence class of $\phi$. 
\end{definition}

\begin{remark}
	Two transformations are equivalent precisely when they differ only by a permutation of the indices in the cospan describing their type: they are ``essentially the same''. For this reason, from now on we shall drop an explicit reference to the equivalence class $\class\phi$ and just reason with the representative $\phi$, except when defining new operations on transformations, like the vertical composition below.
\end{remark}

\begin{definition}\label{def:vertical composition}
	Let $\phi \colon F \to G$ be a transformation as in Definition~\ref{def:transformation}, let $H \colon \B^\gamma \to \C$ be a functor and $\psi \colon G \to H$ a transformation of type
	$
	\begin{tikzcd}[cramped,sep=small]
	\length\beta \ar[r,"\eta"] & m & \ar[l,"\theta"'] \length\gamma
	\end{tikzcd}
	$. The \emph{vertical composition} $\class\psi \circ \class\phi$ is defined as the equivalence class of the transformation $\psi\circ\phi$ of type
	$
	\begin{tikzcd}[cramped,sep=small]
	\length\alpha \ar[r,"\zeta\sigma"] & l & \ar[l,"\xi\theta"'] \length\gamma
	\end{tikzcd}
	$,
	where $\zeta$ and $\xi$ are given by a choice of a pushout
	\begin{equation}\label{eqn:pushout composite type}
		\begin{tikzcd}
			& & \length\gamma \ar[d,"\theta"] \\
			& \length\beta \ar[d,"\tau"'] \ar[r,"\eta"] \ar[dr,phantom,very near end,"\ulcorner"] & m \ar[d,"\xi",dotted] \\
			\length\alpha \ar[r,"\sigma"] & n \ar[r,"\zeta",dotted] & l
		\end{tikzcd}
	\end{equation}
	and the general component $(\psi\circ\phi)_{\bfA}$, for $\bfA \in \B^l$, is the composite:
	\[
	\begin{tikzcd}
	F(\bfA\zeta\sigma) \ar[r,"\phi_{\bfA\zeta}"] & G(\bfA\zeta\tau)=G(\bfA\xi\eta) \ar[r,"\psi_{\bfA\xi}"] & H(\bfA\xi\theta)
	\end{tikzcd}.
	\]
	(Notice that by definition $\phi_{\bfA\zeta} = \phi_{(A_{\zeta1},\dots,A_{\zeta n})}$ requires that the $i$-th variable of $F$ be the $\sigma i$-th element of the list $(A_{\zeta1},\dots,A_{\zeta n})=\bfA\zeta$, which is $A_{\zeta\sigma i}$, hence the domain of $\phi_{\bfA\zeta}$ is indeed $F(\bfA\zeta\sigma)$.)
\end{definition}

Before giving some examples, we introduce the definition of dinaturality of a transformation in one of its variables, as a straightforward generalisation of the classical notion of dinatural transformation in one variable. Recall from p.~\pageref{not:A[X,Y/i]sigma} the meaning of the notation $\substMV \bfA X Y i \sigma$ for $\bfA\in\B^n$, $\sigma \colon \length\alpha \to n$ and $i \in \{1,\dots,n\}$.

\begin{definition}\label{def:dinaturality in i-th variable}
	Let $\phi = (\phi_{A_1,\dots,A_n}) \colon F \to G$ be a transformation as in Definition~\ref{def:transformation}. For $i \in \{1,\dots,n\}$, we say that $\phi$ is \emph{dinatural in $A_i$} (or, more precisely, \emph{dinatural in its $i$-th variable}) if and only if for all $A_1,\dots,A_{i-1}, A_{i+1},\dots,A_n$ objects of $\B$ and for all $f \colon A \to B$ in $\B$ the following hexagon commutes:
	\[
	\begin{tikzcd}
	& F(\subst \bfA A i \sigma) \ar[r,"\phi_{\subst \bfA A i}"] & G(\subst \bfA A i \tau) \ar[dr,"G(\substMV \bfA A f i \tau)"] \\
	F(\substMV \bfA B A i \sigma) \ar[ur,"F(\substMV \bfA f A i \sigma)"] \ar[dr,"F(\substMV \bfA B f \sigma)"'] & & & G(\substMV \bfA A B i \tau) \\
	& F(\subst \bfA B i \sigma) \ar[r,"\phi_{\subst \bfA B i}"'] & G(\subst \bfA B i \tau) \ar[ur,"G(\substMV \bfA f B i \tau)"']
	\end{tikzcd}
	\]
	where $\bfA$ is the $n$-tuple $(A_1,\dots,A_n)$ of the objects above with an additional (unused in this definition) object $A_i$ of $\B$.
\end{definition}

Definition~\ref{def:dinaturality in i-th variable} 
reduces to the well-known notion of dinatural transformation when $\alpha=\beta=[-,+]$ and $n=1$. Our generalisation allows multiple variables at once and the possibility for $F$ and $G$ of having an arbitrary number of copies of $\B$ and $\Op\B$ in their domain, for each variable $i \in \{1,\dots,n\}$. 

\begin{example}\label{ex:delta}
	Let $\C$ be a cartesian category. The diagonal transformation $\delta=(\delta_A \colon A \to A \times A)_{A \in \C}$, classically a natural transformation from $\id\C$ to the diagonal functor, can be equivalently seen in our notations as a transformation $\delta \colon \id\C \to \times$ of type
	$
	\begin{tikzcd}[cramped,sep=small]
	1 \ar[r] & 1 & \ar[l] 2
	\end{tikzcd}.
	$ Of course $\delta$ is dinatural (in fact, natural) in its only variable.
\end{example}

\begin{example}\label{ex:eval}
	Let $\C$ be a cartesian closed category and consider the functor
	\[
	\begin{tikzcd}[row sep=0em]
	\C \times \Op\C \times \C  \ar[r,"T"] & \C \\
	(X,Y,Z) \ar[r,|->] & X \times (Y \Rightarrow Z)
	\end{tikzcd}
	\]
	The evaluation 
	$
	\eval{}{} = \left(\eval A B \colon A \times (A \implies B) \to B\right)_{A,B \in \C} \colon T \to \id\C
	$
	is a transformation of type
	\[
	\begin{tikzcd}[row sep=0em]
	3 \ar[r] & 2 & 1 \ar[l] \\
	1 \ar[r,|->] & 1 & 1 \ar[dl,|->,out=180,in=30] \\[-3pt]
	2 \ar[ur,|->,out=0,in=210]& 2 & \\[-3pt]
	3 \ar[ur,|->,out=0,in=210]
	\end{tikzcd}
	\]
	which is dinatural in both its variables.
\end{example}

\begin{example}\label{ex:Church numeral}
	Let $\C$ be any category, and call $\hom\C \colon \Op \C \times \C \to \Set$ the hom-functor of $\C$. The $n$-th numeral~\cite{dubuc_dinatural_1970}, for $n \in \N$, is the transformation $n \colon \hom\C \to \hom\C$ of type
	$
	\begin{tikzcd}[cramped,sep=small]
	2 \ar[r] & 1 & \ar[l] 2
	\end{tikzcd}
	$
	whose general component $n_A \colon \C(A,A) \to \C(A,A)$ is given, for $A \in \C$ and $g \colon A \to A$, by
	\[
	n_A (g) = g^n,
	\]
	with $0_A (g) = \id A$. 	Then $n$ is dinatural because for all $f \colon A \to B$ the following hexagon commutes:
	\[
	\begin{tikzcd}
	& \C(B,B) \ar[r,"n_B"] & \C(B,B) \ar[dr,"-\circ f"] \\
	\C(B,A) \ar[ur,"f\circ -"] \ar[dr,"-\circ f"'] & & & \C(A,B) \\
	& \C(A,A) \ar[r,"n_A"'] & \C(A,A) \ar[ur,"f \circ -"']
	\end{tikzcd}	
	\]
	It is indeed true that for $h \colon B \to A$, $(f \circ h)^n \circ f = f \circ (h \circ f)^n$: for $n=0$ it follows from the identity axiom; for $n \ge 1$ it is a consequence of associativity of composition.
\end{example}

\paragraph{The graph of a transformation} Given a transformation $\phi$, we now define a graph that reflects its signature, which we shall use to prove our version of Petri\'c's theorem on compositionality of dinatural transformations~\cite{petric_g-dinaturality_2003}. This graph is, as a matter of fact, a \emph{string diagram} for the transformation. String diagrams were introduced by Eilenberg and Kelly in~\cite{eilenberg_generalization_1966} (indeed our graphs are inspired by theirs) and have had a great success in the study of coherence problems (\cite{kelly_coherence_1980,mac_lane_natural_1963}) and monoidal categories in general (\cite{joyal_geometry_1991,joyal_traced_1996}, a nice survey can be found in~\cite{selinger_survey_2010}).

\begin{definition}\label{def:standard graph}
	Let $F \colon \B^\alpha \to \C$ and $G \colon \B^\beta \to \C$ be functors, and let $\phi \colon F \to G$ be a transformation of type
	$
	\begin{tikzcd}[cramped,sep=small]
	\length\alpha \ar[r,"\sigma"] & n & \ar[l,"\tau"'] \length\beta
	\end{tikzcd}
	$. We define its \emph{standard graph} $\graph\phi = (P,T,\inp{(-)},\out{(-)})$ as a directed, bipartite graph as follows:
	\begin{itemize}
		\item $P=\length{\length\alpha + \length\beta}$ and $T=n$ are distinct finite sets of vertices;
		\item $\inp{(-)},\out{(-)} \colon T \to \parts P$ are the input and output functions for elements in $T$: there is an arc from $p \in P$ to $t \in T$ if and only if $p \in \inp t$, and there is an arc from $t$ to $p$ if and only if $p \in \out t$. Indicating with $\injP {\length\alpha} \colon \length\alpha \to P$ and $\injP {\length\beta} \colon \length\beta \to P$ the injections defined as follows:
		\[
		\injP{\length\alpha} (x) = x, \quad \injP{\length\beta} (x) = \length\alpha + x,
		\]
		we have:
		\begin{align*}
		\inp{t}  &=
		\{ \injP {\length\alpha} (p) \mid \sigma (p) = t,\, \alpha_p = +  \} \, \cup \, \{ \injP {\length\beta} (p) \mid \tau (p) = t,\, \beta_p = -  \} \\
		\out{t} &= 
		\{ \injP {\length\alpha} (p) \mid \sigma(p) = t,\, \alpha_p = -  \} \, \cup
		\, \{ \injP {\length\beta} (p) \mid \tau (p) = t,\, \beta_p = +  \}
		\end{align*}
	\end{itemize}
	In other words, elements of $P$ correspond to the arguments of $F$ and $G$, while those of $T$ to the variables of $\phi$. For $t \in T$, its inputs are the covariant arguments of $F$ and the contravariant arguments of $G$ which are mapped by $\sigma$ and $\tau$ to $t$; similarly for its outputs (swapping `covariant' and `contravariant').
\end{definition}

Graphically, we draw elements of $P$ as white or grey boxes (if corresponding to a covariant or contravariant argument of a functor, respectively), and elements of $T$ as black squares. The boxes for the domain functor are drawn at the top, while those for the codomain at the bottom; the black squares in the middle. The graphs of the transformations given in examples \ref{ex:delta}-\ref{ex:Church numeral} are the following:
\begin{itemize}
	\item $\delta=(\delta_A \colon A \to A \times A)_{A \in \C}$ (example \ref{ex:delta}):
	\[
	\begin{tikzpicture}
	\matrix[row sep=1em,column sep=0.5em]{
		& \node (1) [category] {}; \\
		& \node (A) [component] {}; \\
		\node (2) [category] {}; & & \node (3) [category] {}; \\
	};
	\graph[use existing nodes]{
		1 -> A -> {2,3}; 
	};
	\end{tikzpicture}
	\]
	\item $\eval{}{} = \left(\eval A B \colon A \times (A \implies B) \to B\right)_{A,B \in \C}$ (example \ref{ex:eval}):
	\[
	\begin{tikzpicture}
	\matrix[row sep=1em, column sep=1em]{
		\node (1) [category] {}; & & \node (2) [opCategory] {}; & & \node (3) [category] {}; \\
		& \node (A) [component] {}; & & & \node (B) [component] {}; \\
		& & & & \node (4) [category] {}; \\
	};
	\graph[use existing nodes]{
		1 -> A -> 2; 3 -> B -> 4;
	};
	\end{tikzpicture}
	\] 
	\item $n=(n_A \colon \C(A,A) \to \C(A,A))_{A \in \C}$ (example \ref{ex:Church numeral}):
	\[
	\begin{tikzpicture}
	\matrix[row sep=1em, column sep=1em]{
		\node (1) [opCategory] {}; & & \node (2) [category] {};\\
		& \node (A) [component] {};\\
		\node (3) [opCategory] {}; & & \node (4) [category] {};\\
	};
	\graph[use existing nodes]{
		2 -> A -> 1;
		3 -> A -> 4;
	};
	\end{tikzpicture}
	\]
\end{itemize}

\begin{remark}
	Each connected component of $\graph\phi$ corresponds to one variable of $\phi$: the arguments of the domain and codomain of $\phi$ corresponding to (white, grey) boxes belonging to the same connected component are all computed on the same object, when we write down the general component of $\phi$.
	
\end{remark}

\label{discussion:informal-reading-morphisms-in-a-box}This graphical counterpart of a transformation $\phi \colon F \to G$ permits us to represent, in an informal fashion, the dinaturality properties of $\phi$. By writing inside a box a morphism $f$ and reading a graph from top to bottom as ``compute $F$ in the morphisms as they are written in its corresponding boxes, compose that with an appropriate component of $\phi$, and compose that with $G$ computed in the morphisms as they are written in its boxes (treating an empty box as an identity)'', we can express the commutativity of a dinaturality diagram as an informal equation of graphs. (We shall make this precise in Proposition~\ref{prop:fired labelled marking is equal to original one}.) For instance, the dinaturality of examples~\ref{ex:delta}-\ref{ex:Church numeral} can be depicted as follows, where the upper leg of the diagrams are the left-hand sides of the equations:
\begin{itemize}
	\item $\delta=(\delta_A \colon A \to A \times A)_{A \in \C}$ (example \ref{ex:delta}):
	\[
	\begin{tikzcd}
	A \ar[r,"f"] \ar[d,"\delta_A"'] & B \ar[d,"\delta_B"] \\
	A \times A \ar[r,"f \times f"] & B \times B
	\end{tikzcd}
	\qquad
	\begin{tikzpicture}
	\matrix[row sep=1em,column sep=0.5em]{
		& \node (1) [category] {$f$}; \\
		& \node (A) [component] {}; \\
		\node (2) [category] {}; & & \node (3) [category] {}; \\
	};
	\graph[use existing nodes]{
		1 -> A -> {2,3}; 
	};
	\end{tikzpicture}
	\quad = \quad
	\begin{tikzpicture}
	\matrix[row sep=1em,column sep=0.5em]{
		& \node (1) [category] {}; \\
		& \node (A) [component] {}; \\
		\node (2) [category] {$f$}; & & \node (3) [category] {$f$}; \\
	};
	\graph[use existing nodes]{
		1 -> A -> {2,3}; 
	};
	\end{tikzpicture}
	\]
	\item $\eval{}{} = \left(\eval A B \colon A \times (A \implies B) \to B\right)_{A,B \in \C}$ (example \ref{ex:eval}):
	\[\scriptstyle{
		\begin{tikzcd}
		A \times (A' \implies B) \ar[r,"f\times (1 \implies 1)"] \ar[d,"1\times(f\implies 1)"'] & A' \times (A' \implies B) \ar[d,"\eval {A'} B"] \\
		A \times (A \implies B) \ar[r,"\eval A B"] & B
		\end{tikzcd}}
	\qquad
	\begin{tikzpicture}
	\matrix[row sep=1em, column sep=.5em]{
		\node (1) [category] {$f$}; & & \node (2) [opCategory] {}; & & \node (3) [category] {}; \\
		& \node (A) [component] {}; & & & \node (B) [component] {}; \\
		& & & & \node (4) [category] {}; \\
	};
	\graph[use existing nodes]{
		1 -> A -> 2; 3 -> B -> 4;
	};
	\end{tikzpicture}
	\quad = \quad
	\begin{tikzpicture}
	\matrix[row sep=1em, column sep=.5em]{
		\node (1) [category] {}; & & \node (2) [opCategory] {$f$}; & & \node (3) [category] {}; \\
		& \node (A) [component] {}; & & & \node (B) [component] {}; \\
		& & & & \node (4) [category] {}; \\
	};
	\graph[use existing nodes]{
		1 -> A -> 2; 3 -> B -> 4;
	};
	\end{tikzpicture}
	\]
	\[\scriptstyle{
		\begin{tikzcd}
		A \times (A \implies B) \ar[r,"1\times(1\implies g)"] \ar[d,"\eval A B"'] & A \times (A \implies B') \ar[d,"\eval A {B'}"] \\
		B \ar[r,"g"] & B'
		\end{tikzcd}}
	\qquad
	\begin{tikzpicture}
	\matrix[row sep=1em, column sep=.5em]{
		\node (1) [category] {}; & & \node (2) [opCategory] {}; & & \node (3) [category] {$g$}; \\
		& \node (A) [component] {}; & & & \node (B) [component] {}; \\
		& & & & \node (4) [category] {}; \\
	};
	\graph[use existing nodes]{
		1 -> A -> 2; 3 -> B -> 4;
	};
	\end{tikzpicture}
	\quad = \quad
	\begin{tikzpicture}
	\matrix[row sep=1em, column sep=.5em]{
		\node (1) [category] {}; & & \node (2) [opCategory] {}; & & \node (3) [category] {}; \\
		& \node (A) [component] {}; & & & \node (B) [component] {}; \\
		& & & & \node (4) [category] {$g$}; \\
	};
	\graph[use existing nodes]{
		1 -> A -> 2; 3 -> B -> 4;
	};
	\end{tikzpicture}
	\] 
	\item $n=(n_A \colon \C(A,A) \to \C(A,A))_{A \in \C}$ (example \ref{ex:Church numeral}):
	\[\scriptstyle{
		\begin{tikzcd}[column sep={.5cm}]
		& \C(B,B) \ar[r,"n_B"] & \C(B,B) \ar[dr,"{\C(f,1)}"] \\
		\C(B,A) \ar[ur,"{\C(1,f)}"] \ar[dr,"{\C(f,1)}"'] & & & \C(A,B) \\
		& \C(A,A) \ar[r,"n_A"] & \C(A,A) \ar[ur,"{\C(1,f)}"']
		\end{tikzcd}}	
	\qquad
	\begin{tikzpicture}
	\matrix[row sep=1em, column sep=1em]{
		\node (1) [opCategory] {}; & & \node (2) [category] {$f$};\\
		& \node (A) [component] {};\\
		\node (3) [opCategory] {$f$}; & & \node (4) [category] {};\\
	};
	\graph[use existing nodes]{
		2 -> A -> 1;
		3 -> A -> 4;
	};
	\end{tikzpicture}
	\quad = \quad
	\begin{tikzpicture}
	\matrix[row sep=1em, column sep=1em]{
		\node (1) [opCategory] {$f$}; & & \node (2) [category] {};\\
		& \node (A) [component] {};\\
		\node (3) [opCategory] {}; & & \node (4) [category] {$f$};\\
	};
	\graph[use existing nodes]{
		2 -> A -> 1;
		3 -> A -> 4;
	};
	\end{tikzpicture}
	\]
\end{itemize}

All in all, the dinaturality condition becomes, in graphical terms, as follows: \emph{$\phi$ is dinatural if and only if having in $\graph\phi$ one $f$ in all white boxes at the top and grey boxes at the bottom is the same as having one $f$ in all grey boxes at the top and white boxes at the bottom}. 

Not only does $\graph\phi$ give an intuitive representation of the dinaturality properties of $\phi$, but also of the process of composition of transformations. Given two transformations $\phi \colon F \to G$ and $\psi \colon G \to H$ as in Definition~\ref{def:vertical composition}, the act of computing the pushout~(\ref{eqn:pushout composite type}) corresponds to ``glueing together'' $\graph\phi$ and $\graph\psi$ along the boxes corresponding to the functor $G$ (more precisely, one takes the disjoint union of $\graph\phi$ and $\graph\psi$ and then identifies the $G$-boxes), obtaining a composite graph which we will call ${\graph\psi} \circ {\graph\phi}$. The number of its connected components is, indeed, the result of the pushout. That being done, $\graph{\psi\circ\phi}$ is obtained by collapsing each connected component of $\graph\psi\circ\graph\phi$ into a single black square together with the $F$- and $H$-boxes. The following example shows this process. The graph $\graph\psi\circ\graph\phi$ will play a crucial role into the compositionality problem of $\psi\circ\phi$.

\begin{example}\label{ex:acyclic-example}
	Suppose that $\C$ is cartesian closed, fix an object $R$ in $\C$, consider functors
	\[
	\begin{tikzcd}[row sep=0em,column sep=1em]
	\C \times \Op\C \ar[r,"F"] & \C \\
	(A,B) \ar[r,|->] & A \times (B \Rightarrow R)
	\end{tikzcd}
	\quad
	\begin{tikzcd}[row sep=0em,column sep=1em]
	\C \times \C \times \Op\C \ar[r,"G"] & \C \\
	(A,B,C) \ar[r,|->] & A \times B \times (C \Rightarrow R)
	\end{tikzcd} 
	\quad
	\begin{tikzcd}[row sep=0em,column sep=1.5em]
	\C \ar[r,"H"] & \C \\
	A \ar[r,|->] & A \times R
	\end{tikzcd}
	\]
	and transformations $\phi = \delta \times \id{(-)\Rightarrow R} \colon F \to G$
	and $\psi = \id\C \times \eval {(-)} R \colon G \to H$ of types, respectively,
	\[
	\begin{tikzcd}[row sep=0em]
	2 \ar[r,"\sigma"] & 2 & \ar[l,"\tau"'] 3 \\
	1 \ar[r,|->]      & 1 & \ar[l,|->]     1 \\[-3pt]
	2 \ar[r,|->]      & 2 & \ar[ul,|->,out=180,in=-30]    2 \\[-3pt]
	&	  & \ar[ul,|->,out=180,in=-20]    3
	\end{tikzcd}
	\quad\text{and}\quad
	\begin{tikzcd}[row sep=0em]
	3 \ar[r,"\eta"] & 2 & \ar[l,"\theta"'] 1 \\
	1 \ar[r,|->]    & 1 & \ar[l,|->]       1 \\[-3pt]
	2 \ar[r,|->]    & 2 \\[-3pt]
	3 \ar[ur,|->,out=0,in=210]
	\end{tikzcd}
	\]
	so that
	\[
	\phi_{A,B} = \delta_A \times \id{B\implies R} \colon F(A,B) \to G(A,A,B), \, \psi_{A,B} = \id A \times \eval B R \colon G(A,B,B) \to H(A).
	\]
	Then $\psi \circ \phi$ has type $\begin{tikzcd}[cramped,sep=small]
	2 \ar[r] & 1 & \ar[l] 1
	\end{tikzcd}$ and $\graph{\psi}\circ\graph{\phi}$ is:
	\[
	\begin{tikzpicture}
	\matrix[column sep=2.4mm,row sep=0.4cm]{
		&	\node (A) [category] {}; & & & \node(F) [opCategory] {};\\
		&	\node (B) [component] {}; & & & \node(J) [component] {};\\
		\node (C) [category] {}; & & \node(D) [category] {}; & & \node(E) [opCategory] {};\\
		\node (H) [component] {}; & & & \node(I) [component] {};\\
		\node (G) [category] {}; & & & \\
	};
	\graph[use existing nodes]{
		A -> B -> {C, D};
		C -> H -> G;
		D -> I -> E -> J -> F;
	};
	\end{tikzpicture}
	\]
	The two upper boxes at the top correspond to the arguments of $F$, the three in the middle to the arguments of $G$, and the bottom one to the only argument of $H$. This is a connected graph (indeed, $\psi\circ\phi$ depends only on one variable) and by collapsing it into a single black box we obtain $\graph{\psi\circ\phi}$ as it is according to Definition~\ref{def:standard graph}:
	\[
	\begin{tikzpicture}
	\matrix[column sep=.5em,row sep=1em]{
		\node (1) [category] {}; & & \node (2) [opCategory] {};\\
		& \node (A) [component] {}; \\
		& \node (3) [category] {};\\
	};
	\graph[use existing nodes]{
		1 -> A -> {2,3};
	};
	\end{tikzpicture}
	\]
	We have that $\psi \circ \phi$ is a dinatural transformation. (This is one of the  transformations studied by Girard, Scedrov and Scott in~\cite{girard_normal_1992}.) The following string-diagrammatic argument proves that:
	\[
	\begin{split}
	\begin{tikzpicture}[ampersand replacement=\&]
	\matrix[column sep=2.4mm,row sep=0.4cm]{
		\&	\node (A) [category] {$f$}; \& \& \& \node(F) [opCategory] {};\\
		\&	\node (B) [component] {}; \& \& \& \node(J) [component] {};\\
		\node (C) [category] {}; \& \& \node(D) [category] {}; \& \& \node(E) [opCategory] {};\\
		\node (H) [component] {}; \& \& \& \node(I) [component] {};\\
		\node (G) [category] {}; \& \& \& \\
	};
	\graph[use existing nodes]{
		A -> B -> {C, D};
		C -> H -> G;
		D -> I -> E -> J -> F;
	};
	\end{tikzpicture}
	\quad &= \quad
	\begin{tikzpicture}[ampersand replacement=\&]
	\matrix[column sep=2.4mm,row sep=0.4cm]{
		\&	\node (A) [category] {}; \& \& \& \node(F) [opCategory] {};\\
		\&	\node (B) [component] {}; \& \& \& \node(J) [component] {};\\
		\node (C) [category] {$f$}; \& \& \node(D) [category] {$f$}; \& \& \node(E) [opCategory] {};\\
		\node (H) [component] {}; \& \& \& \node(I) [component] {};\\
		\node (G) [category] {}; \& \& \& \\
	};
	\graph[use existing nodes]{
		A -> B -> {C, D};
		C -> H -> G;
		D -> I -> E -> J -> F;
	};
	\end{tikzpicture}
	\quad = \quad
	\begin{tikzpicture}[ampersand replacement=\&]
	\matrix[column sep=2.4mm,row sep=0.4cm]{
		\&	\node (A) [category] {}; \& \& \& \node(F) [opCategory] {};\\
		\&	\node (B) [component] {}; \& \& \& \node(J) [component] {};\\
		\node (C) [category] {}; \& \& \node(D) [category] {$f$}; \& \& \node(E) [opCategory] {};\\
		\node (H) [component] {}; \& \& \& \node(I) [component] {};\\
		\node (G) [category] {$f$}; \& \& \& \\
	};
	\graph[use existing nodes]{
		A -> B -> {C, D};
		C -> H -> G;
		D -> I -> E -> J -> F;
	};
	\end{tikzpicture}
	\\
	&= \quad \begin{tikzpicture}[ampersand replacement=\&]
	\matrix[column sep=2.4mm,row sep=0.4cm]{
		\&	\node (A) [category] {}; \& \& \& \node(F) [opCategory] {};\\
		\&	\node (B) [component] {}; \& \& \& \node(J) [component] {};\\
		\node (C) [category] {}; \& \& \node(D) [category] {}; \& \& \node(E) [opCategory] {$f$};\\
		\node (H) [component] {}; \& \& \& \node(I) [component] {};\\
		\node (G) [category] {$f$}; \& \& \& \\
	};
	\graph[use existing nodes]{
		A -> B -> {C, D};
		C -> H -> G;
		D -> I -> E -> J -> F;
	};
	\end{tikzpicture} \quad =  \quad
	\begin{tikzpicture}[ampersand replacement=\&]
	\matrix[column sep=2.4mm,row sep=0.4cm]{
		\&	\node (A) [category] {}; \& \& \& \node(F) [opCategory] {$f$};\\
		\&	\node (B) [component] {}; \& \& \& \node(J) [component] {};\\
		\node (C) [category] {}; \& \& \node(D) [category] {}; \& \& \node(E) [opCategory] {};\\
		\node (H) [component] {}; \& \& \& \node(I) [component] {};\\
		\node (G) [category] {$f$}; \& \& \& \\
	};
	\graph[use existing nodes]{
		A -> B -> {C, D};
		C -> H -> G;
		D -> I -> E -> J -> F;
	};
	\end{tikzpicture}
	\end{split}
	\]
	The first equation is due to dinaturality of $\phi$ in its first variable;  the second to dinaturality of $\psi$ in its first variable; the third to dinaturality of $\psi$ in its second variable; the fourth equation holds by dinaturality of $\phi$ in its second variable.
\end{example}

The string-diagrammatic argument above is the essence of our proof of Petrić's theorem: we will interpret $\graph\psi \circ \graph\phi$, for arbitrary transformations $\phi$ and $\psi$ as a \emph{Petri Net} whose set of places is $P$ and of transitions is $T$. The dinaturality of $\psi\circ\phi$ will be expressed as a reachability problem and we will prove that, if $\graph\psi \circ \graph\phi$ is acyclic, then $\psi\circ\phi$ is always dinatural because we can always ``move the $f$'s'' from the upper-white boxes and lower-grey boxes all the way to the upper-grey boxes and lower-white boxes, as we did in Example~\ref{ex:acyclic-example}.

\paragraph{Petri Nets}\label{section: Petri Nets}

Petri Nets were invented by Carl Adam Petri in 1962 in \cite{petri_kommunikation_1962}, and have been used since then to model concurrent systems, resource sensitivity and many dynamic systems. A nice survey of their properties was written by Murata in \cite{murata_petri_1989}, to which we refer the reader for more details and examples. Here we shall limit ourselves only to the definitions and the properties of which we will make use in the paper.

\begin{definition}\label{def:Petri Net}
	A \emph{Petri Net} $N$ is a tuple $(P,T,\inp{(-)},\out{(-)})$ where $P$ and $T$ are distinct, finite sets, and  $\inp{(-)},\out{(-)}\colon T \to \parts{P}$ are functions. Elements of $P$ are called \emph{places}, while elements of $T$ are called \emph{transitions}. For $t$ a transition, $\inp t$ is the set of \emph{inputs} of $t$, and $\out t$ is the set of its \emph{outputs}.  A \emph{marking} for $N$ is a function $M \colon P \to \N$.
\end{definition}

Graphically, the elements of $P$ and $T$ are drawn as light-blue circles and black bars respectively. Notice that the graph of a transformation is, as a matter of fact, a Petri Net. We can represent a marking $M$ by drawing, in each place $p$, $M(p)$ \emph{tokens} (black dots). Note that there is at most one arrow from a node to another.

With little abuse of notation, we extend the input and output notation for places too, where
\[
\inp p = \{ t \in T \mid p \in \out{t}  \},  \qquad
\out p = \{ t \in T \mid p \in  \inp t \}.
\]

A pair of a place $p$ and a transition $t$ where $p$ is both an input and an output of $t$ is called \emph{self-loop}. For the purposes of this article, we shall only consider Petri Nets that contain no self-loops.

\begin{definition}
	Let $N$ be a Petri Net. A place $p$ of $N$ is said to be a \emph{source} if $\inp p = \emptyset$, whereas is said to be a \emph{sink} if $\out p = \emptyset$. A source (or sink) place $p$ is said to be \emph{proper} if $\out p \ne \emptyset$ (or $\inp p \ne \emptyset$, respectively).
\end{definition}

We shall need a notion of (directed) path in a Petri Net, which we introduce now. It coincides with the usual notion of path in a graph.

\begin{definition}
	Let $N$ be a Petri Net. A \emph{path} from a vertex $v$ to a vertex $w$ is a finite sequence of vertices $\pi=(v_0,\dots,v_l)$ where $l \ge 1$, $v_0=v$, $v_l=w$ and for all $i \in \{0,\dots,l-1\}$ $v_{i+1} \in v_i \! \LargerCdot \! \cup \! \LargerCdot \! v_i $. Two vertices are said to be \emph{connected} if there is a path from one to the other. If every vertex in $N$ is connected with every other vertex, then $N$ is said to be \emph{weakly connected}.
	
	A \emph{directed path} from a vertex $v$ to a vertex $w$ is a finite sequence of vertices $\pi=(v_0,\dots,v_l)$ such that $v=v_0$, $w=v_l$ and for all $i \in \{0,\dots,l-1\}\,$ $v_{i+1} \in v_i \! \LargerCdot \!$. In this case we say that the path $\pi$ has length $l$. A directed path from a vertex to itself is called a \emph{cycle}, or \emph{loop}; if $N$ does not have cycles, then it is said to be \emph{acyclic}.
	Two vertices $v$ and $w$ are said to be \emph{directly connected} if there is a directed path either from $v$ to $w$ or from $w$ to $v$. 
\end{definition}

We can give a dynamic flavour to Petri Nets by allowing the tokens to “flow” through the nets, that is allowing markings to change according to the following \emph{transition firing rule}.

\begin{definition}
	Let $N=(P,T,\inp{(-)},\out{(-)})$ be a Petri Net, and $M$ a marking for $N$. A transition $t$ is said to be \emph{enabled} if and only if for all $p \in \inp t$ we have $M(p) \ge 1$. An enabled transition may \emph{fire}; the firing of an enabled transition $t$ removes one token from each $p \in \inp t$ and adds one token to each $p \in \out t$, generating the following new marking $M'$:
	\[
	M'(p) = \begin{cases}
	M(p) -1 & p \in \inp t \\
	M(p)+1  & p \in \out t \\
	M(p)    & \text{otherwise}
	\end{cases}
	\]
\end{definition}

\begin{example}\label{my-example}
	Consider the following net:
	\[
	\begin{tikzpicture}[yscale=0.5,xscale=0.70]
	\foreach \i/\u in {1/1,2/1,3/2}
	{
		\foreach \j/\v in {1/0,2/0,3/0,4/1}
		{
			\node[place,tokens=\u,label=above:$p_\i$](p\i) at (2*\i,0){};
			\node[place,tokens=\v,label=below:$q_\j$](q\j) at (2*\j-2,-4){};
			\node[place,label=below:$q_5$](q5) at (8,-4){};
			\node[place,tokens=1,label=above:$p_4$](p4) at (8,0){};
			\node[place,label=above:$p_5$](p5) at (10,0){};
			
			\node[transition,label=right:{$t$}] at (4,-2) {}
			edge [pre] (p\i) edge [post] (q\j) edge [post] (q5);
			
			\node[transition,label=right:$t'$] at (10,-2) {}
			edge [pre] (q5) edge [pre] (p4) edge [post] (p5);
		}
	}
	\end{tikzpicture}
	\]
	There are two transitions, $t$ and $t'$, but only $t$ is enabled. Firing $t$ will change the state of the net as follows:
	\[
	\begin{tikzpicture}[yscale=0.5,xscale=0.70]
	\foreach \i/\u in {1/0,2/0,3/1}
	{
		\foreach \j/\v in {1/1,2/1,3/1,4/2}
		{
			\node[place,tokens=\u,label=above:$p_\i$](p\i) at (2*\i,0){};
			\node[place,tokens=\v,label=below:$q_\j$](q\j) at (2*\j-2,-4){};
			\node[place,tokens=1,label=below:$q_5$](q5) at (8,-4){};
			\node[place,tokens=1,label=above:$p_4$](p4) at (8,0){};
			\node[place,label=above:$p_5$](p5) at (10,0){};
			
			\node[transition,label=right:{$t$}] at (4,-2) {}
			edge [pre] (p\i) edge [post] (q\j) edge [post] (q5);
			
			\node[transition,label=right:$t'$] at (10,-2) {}
			edge [pre] (q5) edge [pre] (p4) edge [post] (p5);
		}
	}
	\end{tikzpicture}
	\]
	Now $t$ is disabled, but $t'$ is enabled, and by firing it we obtain:
	\[
	\begin{tikzpicture}[yscale=0.5,xscale=0.70]
	\foreach \i/\u in {1/0,2/0,3/1}
	{
		\foreach \j/\v in {1/1,2/1,3/1,4/2}
		{
			\node[place,tokens=\u,label=above:$p_\i$](p\i) at (2*\i,0){};
			\node[place,tokens=\v,label=below:$q_\j$](q\j) at (2*\j-2,-4){};
			\node[place,label=below:$q_5$](q5) at (8,-4){};
			\node[place,label=above:$p_4$](p4) at (8,0){};
			\node[place,tokens=1,label=above:$p_5$](p5) at (10,0){};
			
			\node[transition,label=right:{$t$}] at (4,-2) {}
			edge [pre] (p\i) edge [post] (q\j) edge [post] (q5);
			
			\node[transition,label=right:$t'$] at (10,-2) {}
			edge [pre] (q5) edge [pre] (p4) edge [post] (p5);
		}
	}
	\end{tikzpicture}
	\]
\end{example}

\paragraph{The reachability problem and dinaturality} Suppose we have a Petri Net $N$ and an initial marking $M_0$. The firing of an enabled transition in $N$ will change the distribution of tokens from $M_0$ to $M_1$, according to the firing transition rule, therefore a sequence of firings of enabled transitions yields a sequence of markings. A \emph{firing sequence} is denoted by $\sigma = (t_0,\dots,t_n)$ where the $t_i$'s are transitions which fire.

\begin{definition}
	A marking $M$ for a Petri Net $N$ is said to be \emph{reachable} from a marking $M_0$ if there exists a firing sequence $(t_1,\dots,t_n)$ and markings $M_1,\dots,M_n$ where $M_i$ is obtained from $M_{i-1}$ by firing transition $t_i$, for $i \in \{1,\dots,n\}$, and $M_{n}=M$. 
\end{definition}

The reachability problem for Petri Nets consists in checking whether a marking $M$ is or is not reachable from $M_0$. It has been shown that the reachability problem is decidable \cite{kosaraju_decidability_1982,mayr_algorithm_1981}.

\begin{remark}\label{rem:preliminary-discussion}
	The crucial observation that will be at the core of our proof of Petri\'c's theorem is that the firing of an enabled transition in the graph of a dinatural transformation $\phi$ corresponds, under certain circumstances, to the dinaturality condition of $\phi$ in one of its variables. Take, for instance, the $n$-th numeral transformation (see example~\ref{ex:Church numeral}). Call  the only transition $t$, and consider the following marking $M_0$:
	\[
	\begin{tikzpicture}[scale=0.7]
	\node[opCategory] (1) at (-1,1) {};
	\node[category,tokens=1] (2) at (1,1) {};
	\node[opCategory,tokens=1] (3) at (-1,-1) {};
	\node[category] (4) at (1,-1) {};
	
	\node[component,label=left:$t$] {} edge[pre] (2) edge[pre] (3) edge[post] (1) edge[post] (4);
	\end{tikzpicture}
	\]
	Transition $t$ is enabled, and once it fires we obtain the following marking $M_1$:
	\[
	\begin{tikzpicture}[scale=0.7]
	\node[opCategory] (1) at (-1,1) {};
	\node[category,tokens=1] (2) at (1,1) {};
	\node[opCategory,tokens=1] (3) at (-1,-1) {};
	\node[category] (4) at (1,-1) {};
	
	\node[component,label=left:$t$] {} edge[pre] (2) edge[pre] (3) edge[post] (1) edge[post] (4);
	
	\draw[->,snake=snake,segment amplitude=.4mm,segment length=2mm,line after snake=1mm] (1.5,0) -- node[above]{$t$} node[below]{fires} (3.5,0);
	\begin{scope}[xshift=5cm]
	\node[opCategory,tokens=1] (1) at (-1,1) {};
	\node[category] (2) at (1,1) {};
	\node[opCategory] (3) at (-1,-1) {};
	\node[category,tokens=1] (4) at (1,-1) {};
	
	\node[component,label=left:$t$] {} edge[pre] (2) edge[pre] (3) edge[post] (1) edge[post] (4);
	\end{scope}
	\end{tikzpicture}
	\]
	The striking resemblance with the graphical version of the dinaturality condition for $n$ is evident:
	\[
	\begin{tikzpicture}
	\matrix[row sep=1em, column sep=1em]{
		\node (1) [opCategory] {}; & & \node (2) [category] {$f$};\\
		& \node (A) [component] {};\\
		\node (3) [opCategory] {$f$}; & & \node (4) [category] {};\\
	};
	\graph[use existing nodes]{
		2 -> A -> 1;
		3 -> A -> 4;
	};
	\end{tikzpicture}
	\quad = \quad
	\begin{tikzpicture}
	\matrix[row sep=1em, column sep=1em]{
		\node (1) [opCategory] {$f$}; & & \node (2) [category] {};\\
		& \node (A) [component] {};\\
		\node (3) [opCategory] {}; & & \node (4) [category] {$f$};\\
	};
	\graph[use existing nodes]{
		2 -> A -> 1;
		3 -> A -> 4;
	};
	\end{tikzpicture}
	\]
	By treating the ``morphism $f$ in a box'' as a ``token in a place'' of $\graph n$, we have seen that the firing of $t$ generates an equation in $\Set$, namely the one that expresses the dinaturality of $n$. 
\end{remark}

Suppose now we have two composable transformations $\phi$ and $\psi$ dinatural in all their variables, in a category $\C$, together with a graph. We shall make precise how certain markings of $\graph\psi\circ\graph\phi$ correspond to morphisms in $\C$, and how the firing of an enabled transition corresponds to applying the dinaturality of $\phi$ or $\psi$ in one of their variables, thus creating an equation of morphisms in $\C$. Therefore, if the firing of a single transition generates an equality in the category, a sequence of firings of enabled transitions yields a chain of equalities. By individuating two markings $M_0$ and $M_d$, each corresponding to a leg of the dinaturality hexagon for $\psi\circ\phi$ we want to prove is commutative, and by showing that $M_d$ is reachable from $M_0$, we shall have proved that $\psi\circ\phi$ is dinatural.

We are now ready to present and prove the first main result of this article. For the rest of this section, fix transformations $\phi \colon F_1 \to F_2$ and $\psi \colon F_2 \to F_3$  where 
\begin{itemize}
	\item $F_i \colon \B^{\alpha^i} \to \C$ is a functor for all $i \in \{1,2,3\}$, 
	\item $\phi$ and $\psi$ have type, respectively,
	\[
	\begin{tikzcd}
	\length{\alpha^1} \ar[r,"\sigma_1"] & k_1 & \length{\alpha^{2}} \ar[l,"\tau_1"']
	\end{tikzcd}
	\qquad \text{and} \qquad
	\begin{tikzcd}
	\length{\alpha^2} \ar[r,"\sigma_2"] & k_2 & \length{\alpha^{3}}. \ar[l,"\tau_2"']
	\end{tikzcd}
	\]
\end{itemize}
We shall establish a sufficient condition for the dinaturality of $\psi \circ \phi$ in some of its variables. However, since we are interested in analysing the dinaturality of the composition in each of its variables \emph{separately}, we start by assuming that $\psi\circ\phi$ depends on only one variable, i.e. has type 
$
\begin{tikzcd}[cramped,sep=small]
\length{\alpha^1} \ar[r] & 1 & \length{\alpha^{3}} \ar[l],
\end{tikzcd}
$
and that $\phi$ and $\psi$ are dinatural in all their variables. In this case, we have to show that the following hexagon commutes for all $f \colon A \to B$, recalling that  $\funminplusconst {F_1} B A$ is the result of applying functor $F_1$ in $B$ in all its contravariant arguments and in $A$ in all its covariant ones:
\begin{equation}\label{eqn:compositionality-hexagon}
	\begin{tikzcd}[column sep=1cm]
		& \funminplusconst {F_1} A A \ar[r,"\phi_{A\dots A}"] &  \funminplusconst {F_2} A A \ar[r,"\psi_{A \dots A}"]  &  \funminplusconst {F_3} A A \ar[dr,"\funminplusconst {F_3} 1 f"] \\
		\funminplusconst {F_1} B A \ar[ur,"\funminplusconst {F_1} f 1"] \ar[dr,"\funminplusconst {F_1} 1 f"'] & &   & & \funminplusconst {F_3} A B \\
		& \funminplusconst {F_1} B B \ar[r,"\phi_{B\dots B}"']  & \funminplusconst {F_2} B B \ar[r,"\psi_{B \dots B}"']  & \funminplusconst {F_3} B B \ar[ur,"\funminplusconst {F_3} f 1"']
	\end{tikzcd}
\end{equation}

The theorem we want to prove is then the following.

\begin{theorem}\label{theorem:acyclic implies dinatural}
	Let $\phi$ and $\psi$ be transformations which are dinatural in all their variables and such that $\psi\circ\phi$ depends on only one variable. If \,$\graph\psi \circ \graph\phi$ is acyclic, then $\psi\circ\phi$ is a dinatural transformation. 
\end{theorem}

The above is a direct generalisation of Eilenberg and Kelly's result on \emph{extranatural transformations} \cite{eilenberg_generalization_1966}, which are dinatural transformations where either the domain or the codomain functor is constant. For example, $\eval{}{}$ is extranatural in its first variable. They worked with the additional assumption that $\graph\phi$ and $\graph\psi$ do not contain any ramifications, that is, the white and grey boxes are always linked in pairs, and they also proved that if the composite graph is acyclic, then the composite transformation is again extranatural. Their condition is also ``essentially necessary''  in the sense that if we do create a cycle upon constructing $\graph\psi \circ \graph\phi$, then that means we are in a situation like this:
\[
\begin{tikzpicture}
\matrix[column sep=1em,row sep=1em]{
	& \node[component] (A) {}; \\
	\node[opCategory] (1) {}; & & \node[category] (2) {};\\
	& \node[component] (B) {};\\
};
\graph[use existing nodes]{
	1 -> A -> 2 -> B -> 1;
};
\end{tikzpicture}
\]
where we have a transformation between constant functors. Such a family of morphisms is (extra)natural precisely when it is constant (that is, if every component is equal to the same morphism) on each connected component of the domain category.

As already said in Remark~\ref{rem:preliminary-discussion}, the key to prove this theorem is to see $\graph\psi \circ \graph\phi$ as a Petri Net, reducing the dinaturality of $\psi\circ\phi$ to the reachability problem for two markings we shall individuate. We begin by unfolding the definition of $\graph\psi \circ \graph\phi$: we have $\graph\psi \circ \graph\phi = (P,T,\inp{(-)},\out{(-)})$ where $P = \length{\alpha^1} + \length{\alpha^2} + \length{\alpha^{3}}$, $T = k_1 + k_2$ and, indicating with $\injP i \colon \length{\alpha^i} \to P$ and $\injT i \colon k_i \to T$ the injections defined similarly to $\injP{\length\alpha}$ and $\injP{\length\beta}$ in Definition~\ref{def:standard graph}, 
\begin{equation}\label{input-output-transitions}
	\begin{aligned}
		\inp{(\injT i (t))} &= \, 
		\{ \injP i (p) \mid \sigma_i(p) = t,\, \alpha^i_p = +  \} \, \cup \, \{ \injP {i+1} (p) \mid \tau_i(p) = t,\, \alpha^{i+1}_p = -  \}, \\
		\out{(\injT i (t))} &= \, 
		\{ \injP i (p) \mid \sigma_i(p) = t,\, \alpha^i_p = -  \} \, \cup
		\, \{ \injP {i+1} (p) \mid \tau_i(p) = t,\, \alpha^{i+1}_p = +  \}.
	\end{aligned}
\end{equation}
For the rest of this section, we shall reserve the names $P$ and $T$ for the sets of places and transitions of $\graph\psi \circ \graph\phi$.

\begin{remark}\label{rem:graph of a transformation is FBCF}
	Since $\sigma_i$ and $\tau_i$ are functions, we have that $\length{\inp p}, \length{\out p} \le 1$ and also that $\length{\inp p \cup \out p }\ge 1$ for all $p\in P$. With a little abuse of notation then, if $\inp p = \{t\}$ then we shall simply write $\inp p = t$, and similarly for $\out p$.
\end{remark}

\paragraph{Labelled markings as morphisms} We now show how to formally translate certain markings of $\graph\psi \circ \graph\phi$ in actual morphisms of $\C$. The idea is to treat every token in the net as a fixed, arbitrary morphism $f \colon A \to B$ of $\C$ and then use the idea discussed on p.~\pageref{discussion:informal-reading-morphisms-in-a-box}.

However, not all possible markings of $\graph\psi \circ \graph\phi$ have a corresponding morphism in $\C$. For example, if $M$ is a marking and $p$ is a place such that $M(p)>1$, it makes no sense to ``compute a functor $F_i$ in $f$ twice'' in the argument of $F_i$ corresponding to $p$. Hence, only markings $M \colon P \to \{0,1\}$ can be considered. Moreover, we have to be careful with \emph{where} the marking puts tokens: if a token corresponds to a morphism $f \colon A \to B$, we have to make sure that there are no two consecutive tokens (more generally, we have to make sure that there is at most one token in every directed path), otherwise a naive attempt to assign a morphism to that marking might end up with type-checking problems. For instance, consider the diagonal transformation in a Cartesian category $\C$ (example \ref{ex:delta}) and the following marking:
\[
\begin{tikzpicture}
\matrix[row sep=1em,column sep=0.5em]{
	& \node (1) [category,tokens=1] {}; \\
	& \node (A) [component] {}; \\
	\node (2) [category,tokens=1] {}; & & \node (3) [category,tokens=1] {}; \\
};
\graph[use existing nodes]{
	1 -> A -> {2,3}; 
};
\end{tikzpicture}
\]
The token on the top white box should be interpreted as $\id\C(f) \colon A \to B$, hence the black middle box should correspond to the $B$-th component of the family $\delta$, that is $\delta_B \colon B \to B \times B$. However, the bottom two white boxes are read as $f \times f \colon A \times A \to B \times B$, which cannot be composed with $\delta_B$.

We therefore introduce the notion of \emph{labelled marking}, which consists of a marking together with a labelling of the transitions, such that a certain coherence condition between the two is satisfied. This constraint will ensure that every labelled marking 
corresponds to 
a morphism of $\C$. We will then use only \emph{some} labelled markings to prove our compositionality theorem.

\begin{definition}\label{def:labelled marking}
	Consider $f \colon A \to B$ a morphism in $\C$. A \emph{labelled marking} for $\graph\psi \circ \graph\phi$ is a triple $(M,L,f)$ where functions $M \colon P \to \{0,1\}$ and $L \colon T \to \{A,B\}$ are such that for all $p \in P$
	\[
	M(p)=1 \implies L(\inp p) = A, \, L(\out p ) = B
	\]
	\[
	M(p)=0 \implies  L(\inp p ) = L(\out p ) 
	\]
	These conditions need to be satisfied only when they make sense; for example if $M(p) = 1$ and $\inp p = \emptyset$, condition $L(\inp p) = A$ is to be ignored.
\end{definition}

We are now ready to assign a morphism in $\C$ to every labelled marking by reading a token in a place as a morphism $f$ in one of the arguments of a functor, while an empty place corresponds to the identity morphism of the label of the transition of which the place is an input or an output.
\begin{definition}\label{def:morphism for labelled marking}
	Let $(M,L,f\colon A \to B)$ be a labelled marking. We define a morphism $\mor M L f$ in $\C$ as follows:
	\[
	\mor M L f = F_1(x^1_1,\dots,x^1_{\length{\alpha^1}});\phi_{X^1_1\dots X^1_{k_1}} ; 
	F_2(x^2_1,\dots,x^2_{\length{\alpha^2}}); \psi_{X^2_1\dots X^2_{k_2}} ; F_{3}(x^{3}_1,\dots,x^{3}_{\length{\alpha^{3}}})
	\]
	where
	\[
	x^i_j = \begin{cases}
	f \quad &M(\injP i (j)) = 1 \\
	\id{L(t)} \quad & M(\injP i (j)) = 0 \land t \in \inp {\injP v (j)} \cup \out {\injP v (j)} 
	\end{cases}
	\qquad
	X^i_j = L(\injT i (j)).
	\]
	for all $i \in \{1,2,3\}$ and $j\in\{1,\dots,\length{\alpha^i}\}$. (Recall that $\injP i \colon \length{\alpha^i} \to P$ and $\injT i \colon k_i \to T$ are the injections defined similarly to $\injP{\length\alpha}$ and $\injP{\length\beta}$ in Definition~\ref{def:standard graph}.)
\end{definition}

It is easy to see that $\mor M L f$ is indeed a morphism in $\C$, by checking that the maps it is made of are actually composable using the definition of labelled marking and of $\graph\psi \circ \graph\phi$.

What are the labelled markings corresponding to the two legs of diagram~(\ref{eqn:compositionality-hexagon})? In the lower leg of the hexagon, $f$ appears in all the covariant arguments of $F_1$ and the contravariant ones of $F_{3}$, which correspond in $\graph\psi \circ \graph\phi$ to those places which have no inputs (in Petri nets terminology, \emph{sources}), and all variables of $\phi$ are equal to $B$; in the upper leg, $f$ appears in those arguments corresponding to places with no outputs (\emph{sinks}), and $\psi$ is computed in $A$ in each variable. Hence, the lower leg is $\mor {M_0} {L_0} f$ while the upper leg is $\mor {M_d} {L_d} f$, where:
\begin{equation}\label{eqn:markings-definitions}
	\begin{aligned}
		M_0(p)&=\begin{cases}
			1 & \inp p = \emptyset \\
			0 & \text{otherwise}
		\end{cases} \quad
		&
		M_d(p)&=\begin{cases}
			1 & \out p = \emptyset \\
			0 & \text{otherwise}
		\end{cases}
		\\[.5em]
		L_0(t) &= B 
		&
		L_d(t) &= A
	\end{aligned}
\end{equation}
for all $p\in P$ and $t \in T$. It is an immediate consequence of the definition that $(M_0,L_0,f)$ and $(M_d,L_d,f)$ so defined are labelled markings.

We aim to show that $M_d$ is reachable from $M_0$ by means of a firing sequence that preserves the morphism $\mor {M_0} {L_0} f$. In order to do so, we now prove that firing a $B$-labelled transition in an arbitrary labelled marking $(M,L,f)$ generates a new labelled marking, whose associated morphism in $\C$ is still equal to $\mor M L f$.

\begin{proposition}\label{prop:fired labelled marking is equal to original one}
	Let $(M,L,f)$ be a labelled marking, $t \in T$ an enabled transition such that $L(t) = B$. Consider 
	\begin{equation}\label{markings after firing definition}
		\begin{tikzcd}[row sep=0em,column sep=1em,ampersand replacement=\&]
			P \ar[r,"M'"] \& \{0,1\} \& \& \& \& \& T \ar[r,"L'"] \& \{A,B\} \\
			p \ar[r,|->] \& \begin{cases}
				0 & p \in \inp t \\
				1 & p \in \out t \\
				M(p) & \text{otherwise}
			\end{cases}
			\& \& \& \& \&
			s \ar[r,|->] \& \begin{cases}
				A & s = t \\
				L(s) & s \ne t
			\end{cases}
		\end{tikzcd}
	\end{equation}
	Then $(M',L',f)$ is a labelled marking and $\mor M L f = \mor {M'} {L'} f$.
\end{proposition} 
\begin{proof}
	By definition of labelled marking, if $\out t \ne \emptyset$ and $L(t) = B$ then $M(p) = 0$ for all $p \in \out t$, because if there were a $p \in \out t$ with $M(p) = 1$, then $L(t) = A$. $M'$ is therefore the marking obtained from $M$ when $t$ fires once. It is easy to see that $(M',L',f)$ is a labelled marking by simply checking the definition. 
	
	We have now to prove that $\mor M L f = \mor {M'} {L'} f$. Since $t \in T$, we have $t = \injT u (i)$ for some $u \in \{1,2\}$ and $i \in \{1,\dots,k_u\}$. The fact that $t$ is enabled in $M$, together with the definition of $\graph\psi \circ \graph\phi$ (\ref{input-output-transitions}) and Definition~\ref{def:morphism for labelled marking}, ensures that, in the notations of Definition~\ref{def:morphism for labelled marking},
	\begin{align*}
	\sigma_u(j) = i \land \alpha^u_j = + &\implies x^u_j = f  \\
	\sigma_u(j) = i \land \alpha^u_j = - &\implies x^u_j = \id B \\
	\tau_u(j) = i \land \alpha^{u+1}_j = + &\implies x^{u+1}_j = \id B \\
	\tau_u(j) = i \land \alpha^{u+1}_j = - &\implies x^{u+1}_j = f 
	\end{align*}
	hence we can apply the dinaturality of $\phi$ or $\psi$ (if, respectively, $u=1$ or $u=2$) in its $i$-th variable. To conclude, one has to show that the morphism obtained in doing so is the same as $\mor {M'} {L'} f$, which is just a matter of identity check. The details can be found in the second author's thesis~\cite{santamaria_towards_2019}.\qed
\end{proof}

It immediately follows that a sequence of firings of $B$-labelled transitions gives rise to a labelled marking whose associated morphism is still equal to the original one, as the following Proposition states.

\begin{corollary}\label{cor:reachability-implies-equality}
	Let $\mor M L f$ be a labelled marking, $M'$ a marking reachable from $M$ by firing only $B$-labelled transitions $t_1,\dots,t_m$, $L' \colon T \to \{A,B\}$ defined as:
	\[
	L'(s) = \begin{cases}
	A & s = t_i \text{ for some $i \in \{1,\dots,m\}$} \\
	L(s) & \text{otherwise}
	\end{cases}
	\]
	Then $(M', L',f)$ is a labelled marking and $\mor M L f = \mor {M'} {L'} f$.
\end{corollary}

Now all we have to show is that $M_d$ is reachable from $M_0$ (see~(\ref{eqn:markings-definitions})) by only firing $B$-labelled transitions: it is enough to make sure that each transition is fired at most once to satisfy this condition. We shall work on a special class of Petri Nets, to which our $\graph\psi \circ \graph\phi$ belongs (Remark~\ref{rem:graph of a transformation is FBCF}), where all places have at most one input and at most one output.

\begin{definition}\label{def:FBCF petri net}
	A Petri Net is said to be \emph{forward-backward conflict free} (FBCF) if for all $p$ place $\length{\inp p} \le 1$ and $\length{\out p} \le 1$.
\end{definition}

\begin{theorem}\label{thm:acyclic-implies-reachable}
	Let $N$ be an acyclic FBCF Petri Net and let $M_0$, $M_d$ be the only-source and only-sink markings as in~(\ref{eqn:markings-definitions}). Then $M_d$ is reachable from $M_0$ by firing each transition exactly once.
\end{theorem}
\begin{proof}
	We proceed by induction on the number of transitions in $N$. If $N$ has no transitions at all, then every place is both a source and a sink, and $M_0$ and $M_d$ coincide, therefore there is nothing to prove. Now, let $n \ge 0$, suppose that the theorem holds for Petri Nets that have $n$ transitions and assume that $N$ has $n+1$ transitions.
	
	Define, given $t$ and $t'$ two transitions, $t \le t'$ if and only if there exists a directed path from $t$ to $t'$. The relation $\le$ so defined is reflexive, transitive and antisymmetric (because $N$ is acyclic), hence it is a partial order on $T$, the set of transitions of $N$. Now, $T$ is finite by definition, hence it has at least one minimal element $t_0$. Since $t_0$ is minimal, every (if any) input of $t_0$ is a source, therefore $t_0$ is enabled in $M_0$. Now, fire $t_0$ and call $M_1$ the resulting marking. Consider the subnet $N'$ obtained from $N$ by removing $t_0$ and all its inputs. Since $N$ is forward-backward conflict free, we have that all the outputs of $t_0$ are sources in $N'$. This means that $N'$ is an acyclic FBCF Petri Net: by inductive hypothesis, we have that $M_d$ (restricted to $N'$) is reachable from $M_1$ in $N'$, and therefore $M_d$ is reachable from $M_0$ in $N$.\qed
\end{proof}

\begin{remark}
	Theorem~\ref{thm:acyclic-implies-reachable} is an instance of Hiraishi and Ichikawa's result on reachability for arbitrary markings in arbitrary acyclic Petri Nets~\cite{hiraishi_class_1988}. Our proof is an adapted version of theirs for the special case of FBCF Petri Nets and the particular markings $M_0$ and $M_d$ that put one token precisely in every source and in every sink respectively.
\end{remark}

We are now ready to give an alternative proof to the first half of Petri\'{c}'s theorem~\cite{petric_g-dinaturality_2003} that solved the compositionality problem of dinatural transformations.

\begin{proofMainTheorem}
	Let $f \colon A \to B$ be a morphism in $\C$, and define labelled markings $(M_0,L_0,f)$ and $(M_d,L_d,f)$ as in~(\ref{eqn:markings-definitions}). Then $\mor {M_0} {L_0} f$ is the lower leg of~(\ref{eqn:compositionality-hexagon}), while $\mor {M_d} {L_d} f$ is the upper leg. By theorem~\ref{thm:acyclic-implies-reachable}, marking $M_d$ is reachable from $M_0$ by firing each transition of $\graph\psi \circ \graph\phi$ exactly once, hence by only firing $B$-labelled transitions. By Proposition~\ref{cor:reachability-implies-equality}, we have that the hexagon~(\ref{eqn:compositionality-hexagon}) commutes. \qed
\end{proofMainTheorem}

Theorem~\ref{theorem:acyclic implies dinatural} can then be straightforwardly generalised to the case in which $\psi\circ\phi$ depends on $n$ variables for an arbitrary $n$. Suppose then that the type of $\psi\circ\phi$ is given by the following pushout: 
\begin{equation}\label{eqn:pushout2}
	\mkern0mu \begin{tikzcd}
		& & \length{\alpha^3} \ar[d, "\tau_2"] \\
		& \length{\alpha^2} \ar[d, "\tau_1"'] \ar[dr, phantom, "\ulcorner" very near start] \ar[r, "\sigma_2"] & k_2 \ar[d, dotted, "\xi"] \\
		\length{\alpha^1} \ar[r, "\sigma_1"] & k_1 \ar[r, dotted, "\zeta"] & n
	\end{tikzcd}
\end{equation}
$\graph\psi \circ \graph\phi$ now has $n$ connected components, and a sufficient condition for the dinaturality of $\psi\circ\phi$ in its $i$-th variable is that $\phi$ and $\psi$ are dinatural in all those variables of theirs which are ``involved'', as it were, in the $i$-th connected component of $\graph\psi \circ \graph\phi$ \emph{and} such connected component is acyclic.

\begin{theorem}\label{theorem:acyclicity implies dinaturality GENERAL}
	In the notations above, let $i\in \{1,\dots,n\}$. If $\phi$ and $\psi$ are dinatural in all the variables in, respectively, $\zeta^{-1}\{i\}$ and $\xi^{-1}\{i\}$ (with $\zeta$ and $\xi$ given by the pushout~(\ref{eqn:pushout2})), and if the $i$-th connected component of $\graph\psi \circ \graph\phi$ is acyclic, then $\psi\circ\phi$ is dinatural in its $i$-th variable.
\end{theorem}

We then have a straightforward corollary.

\begin{corollary}
	Let $\phi\colon F \to G$ and $\psi \colon G \to H$ be transformations which are dinatural in all their variables. If $\graph\psi \circ \graph\phi$ is acyclic, then $\psi\circ\phi$ is dinatural in all its variables.
\end{corollary}

One can generalise even further Theorem~\ref{theorem:acyclicity implies dinaturality GENERAL} by considering $k$ consecutive dinatural transformations $\phi_1,\dots,\phi_k$, instead of just two, and reasoning on the acyclicity of the connected components of the composite graph $\graph{\phi_k} \circ \dots \circ  \graph{\phi_1}$, obtained by ``glueing together'' the standard graphs of the $\phi_i$'s along the common interfaces (formally this would be a composite performed in the category $\gc$ introduced in Definition~\ref{definition:graph category}).

\begin{theorem}\label{theorem:compositionality with complicated graphs}
	Let $\phi_j \colon F_j \to F_{j+1}$ be transformations of type 
	$
	\begin{tikzcd}[cramped,sep=small]
	\length{\alpha^j} \ar[r,"\sigma_j"] & n_j & \ar[l,"\tau_j"'] \length{\alpha^{j+1}}
	\end{tikzcd}
	$
	for $j \in \{1,\dots,k\}$.
	Suppose that the type of $\phi_k \circ \dots \phi_1$ is computed by the following pushout-pasting:
	\[
	\begin{tikzcd}
	&   &   &   & \length{\alpha^{k+1}} \ar[d,"\tau_k"'] \\
	&   &   & \length{\alpha^k} \ar[r,"\sigma_k"] \ar[d] \ar[dr, phantom, "\ulcorner" very near start] & n_k \ar[d] \\
	&   & \length{\alpha^3} \ar[r] \ar[ur,sloped,phantom,"\dots"] \ar[d,"\tau_2"'] \ar[dr, phantom, "\ulcorner" very near start] & \dots \ar[r] \ar[d] \ar[dr, phantom, "\ulcorner" very near start] & \dots \ar[d] \\
	& \length{\alpha^2} \ar[r,"\sigma_2"] \ar[d,"\tau_1"'] \ar[dr, phantom, "\ulcorner" very near start] & n_2 \ar[r] \ar[d] \ar[dr, phantom, "\ulcorner" very near start] & \dots \ar[r] \ar[d] \ar[dr, phantom, "\ulcorner" very near start] & \dots \ar[d] \\
	\length{\alpha^1} \ar[r,"\sigma_1"] & n_1 \ar[r] & \dots \ar[r] & \dots \ar[r] & l
	\end{tikzcd}
	\]
	Let $\xi_j \colon n_j \to l$ be the map given by any path of morphisms from $n_j$ to $l$ in the above diagram. If the $i$-th connected component of $\graph{\phi_k} \circ \dots \circ \graph{\phi_1}$ (composite calculated in $\gc$) is acyclic and if for all $j \in  \{1,\dots,k\}$, for all $x \in \xi_j^{-1} \{i\}$ the transformation $\phi_j$ is dinatural in its $x$-th variable, then $\phi_k \circ \dots \circ \phi_1$ is dinatural in its $i$-th variable.
\end{theorem}
\begin{proof}
	The proof is essentially the same of Theorem~\ref{theorem:acyclicity implies dinaturality GENERAL}, where instead of two transformations we have $k$: one defines labelled markings $(M_0, L_0,f)$ and $(M_d,L_d,f)$ corresponding to the two legs of the dinaturality hexagon of $\phi_k \circ \dots \circ \phi_1$ in its $i$-th variable, and uses Theorem~\ref{thm:acyclic-implies-reachable} to prove that $M_d$ is reachable from $M_0$, thus showing the hexagon commutes. \qed
\end{proof}

\begin{remark}
	In~\cite{girard_normal_1992}, the authors had to prove the dinaturality of families of morphisms obtained by composing several transformations that are dinatural by assumption. They showed that the dinaturality hexagons for such composites commute by filling them with a trellis of commutative diagrams, stating functoriality properties and dinaturality of the building blocks. Theorem~\ref{theorem:compositionality with complicated graphs} provides an alternative  way to do that: one can simply draw the composite graph of the involved transformations, notice that the resulting Petri Net is always acyclic, and thus infer the dinaturality of the composite. 
\end{remark}

\paragraph{An ``essentially necessary'' condition for compositionality} The other half of Petri\'c's theorem can also be shown with the help of the theory of Petri Nets. One can prove that if $N$ is a weakly connected FBCF Petri Net with at least one proper source or one proper sink and $M_0$ and $M_d$ are the only-source and only-sink markings as before, then a necessary condition for the reachability of $M_d$ from $M_0$ is that every transition in $N$ must fire at least once. The intuition behind this is that there must be at least one transition $t$ which fires, because $M_0$ and $M_d$ are not equal (in the hypothesis that $N$ has at least one proper sink or proper source), and if a transition $t$ fires once, then all the transitions that are connected to it must fire as well: in order for $t$ to fire it must be enabled, hence those transitions which are between the source places and $t$ must fire to move the tokens to the input places of $t$; equally, if $t$ fires, then also all those transitions ``on the way'' from $t$ to the sink places must fire, otherwise some tokens would get stuck in the middle of the net, in disagreement with $M_d$. As a consequence of this fact, we have a sort of inverse of Theorem~\ref{thm:acyclic-implies-reachable}.

\begin{theorem}\label{thm:reachability-implies-acyclicity}
	Let $N$ be weakly connected with at least one proper source or one proper sink place. If $M_d$ is reachable from $M_0$, then $N$ is acyclic.
\end{theorem}
\begin{proof}
	Suppose that $N$ contains a directed, circular path $\pi=(v_0,\dots,v_{2l})$ where $v_0 = v_{2l}$ is a place. Then each $v_{2i}$ is not a source, given that it is the output of $v_{2i-1}$, hence $M_0(v_{2i})=0$ for all $i \in \{1,\dots,l\}$. This means that $v_{2i+1}$ is disabled in $M_0$, therefore it will not fire when transforming $M_0$ into $M_1$. Then also $M_1(v_{2i})=0$. Using the same argument we can see that none of the transitions in the loop $\pi$ can fire, thus $M_d$ cannot be reached by $M_0$. \qed
\end{proof}

In other words, if $N$ contains a loop—in the hypothesis that $N$ is weakly connected and has at least one proper source or sink place—then $M_d$ is \emph{not} reachable from $M_0$. In the case of $N=\graph\psi \circ \graph\phi$, given the correspondence between the dinaturality condition of $\phi$ and $\psi$ in each of their variables and the firing of the corresponding transitions, this intuitively means that $\psi\circ\phi$ cannot be proved to be dinatural as a sole consequence of the dinaturality of $\phi$ and $\psi$ when $\graph\psi \circ \graph\phi$ is cyclic. Therefore, acyclicity is not only a \emph{sufficient} condition for the dinaturality of the composite transformation, but also ``essentially necessary'': if the composite happens to be dinatural despite the cyclicity of the graph, then this is due to some ``third'' property, like the fact that certain squares of morphisms are pullbacks or pushouts. The interested reader can find a detailed formalisation of this intuition in the second author's thesis~\cite{santamaria_towards_2019}, where a syntactic category generated by the equations determined by the dinaturality conditions of $\phi$ and $\psi$ was considered, and where it was shown that in there $\psi \circ \phi$ is \emph{not} dinatural in a similar way to Petri\'c's approach in~\cite{petric_g-dinaturality_2003}.

\section{Horizontal compositionality of dinatural transformations}\label{chapter horizontal}

Horizontal composition of natural transformations is co-protagonist, together with vertical composition, in the classical Godement calculus. In this section we define a new operation of horizontal composition for dinatural transformations, generalising the well-known version for natural transformations. We also study its algebraic properties, proving it is associative and unitary. Remarkably, horizontal composition behaves better than vertical composition, as it is \emph{always} defined between dinatural transformations of matching type.

\subsection{From the Natural to the Dinatural}

Horizontal composition of natural transformations \cite{mac_lane_categories_1978} 
is a well-known operation which is rich in interesting properties: it is associative, unitary and compatible with vertical composition. As such, it makes $\mathbb{C}\mathrm{at}$ a strict 2-category. Also, it plays a crucial role in the calculus of substitution of functors and natural transformations developed by Kelly in \cite{kelly_many-variable_1972}; in fact, as we have seen in the introduction, it is at the heart of Kelly's abstract approach to coherence. 
An appropriate generalisation of this notion for dinatural transformations seems to be absent in the literature: in this section we propose a working definition, as we shall see. The best place to start is to take a look at the usual definition for the natural case.

\begin{definition}\label{def:horizontal composition natural transformations}
	Consider (classical) natural transformations
	\[
	\begin{tikzcd}
	\A \ar[r,bend left,"F"{above},""{name=F,below}]{} \ar[r,bend right,"G"{below},""{name=G}] 
	& \B \ar[r,bend left,"H"{above},""{name=H,below}]{} \ar[r,bend right,"K"{below},""{name=K}] 
	& \C	
	\arrow[Rightarrow,from=F,to=G,"\phi"]
	\arrow[Rightarrow,from=H,to=K,"\psi"]
	\end{tikzcd}
	\]
	The horizontal composition $\hc \fst \snd \colon HF \to KG$ is the natural transformation whose $A$-th component, for $A \in \A$, is either leg of the following commutative square:
	\begin{equation}\label{eqn:horCompNatTransfSquare}
		\mkern0mu \begin{tikzcd}
			HF(A) \ar[r,"\psi_{F(A)}"] \ar[d,"H(\phi_A)"'] & KF(A) \ar[d,"K(\phi_A)"] \\
			HG(A) \ar[r,"\psi_{G(A)}"] & KG(A)
		\end{tikzcd}
	\end{equation}
\end{definition}

Now, the commutativity of (\ref{eqn:horCompNatTransfSquare}) is due to the
naturality of $\psi$; the fact that $\hc \phi \psi$ is in turn a natural
transformation is due to the naturality of both $\phi$ and $\psi$. However, in
order to \emph{define} the family of morphisms $\hc \phi \psi$, all we have to
do is to apply the naturality condition of $\psi$ to the components of $\phi$,
one by one. We apply the very same idea to  dinatural
transformations, leading to the following preliminary definition
for classical dinatural transformations.  

\begin{definition}\label{def:horCompDef}
	Let $\fst\colon\fstDom \to \fstCoDom$ and $\snd \colon \sndDom \to \sndCoDom$ be dinatural transformations of type
	$
	\begin{tikzcd}[cramped,sep=small]
	2 \ar[r] & 1 & 2 \ar[l]
	\end{tikzcd}
	$,
	where $\fstDom, \fstCoDom \colon \Op\A \times \A \to \B$ and $\sndDom, \sndCoDom \colon \Op\B \times \B \to \C$. The \emph{horizontal composition} $\hc \fst \snd$ is the family of morphisms 
	\[
	\bigl((\hc \fst \snd)_A\colon 	\sndDom(\fstCoDom(A,A), \fstDom(A,A)) \to  \sndCoDom(\fstDom(A,A),\fstCoDom(A,A))\bigr)_{A \in \A}
	\]
	where the general component $(\hc \fst \snd)_A$ is given, for any object $A \in \A$, by either leg of the following commutative hexagon:
	\[
	\begin{tikzcd}[column sep=.8cm,font=\normalsize]
	& \sndDom(\fstDom(A,A),\fstDom(A,A)) \ar[r,"{\snd_{\fstDom(A,A)}}"]  & \sndCoDom(\fstDom(A,A),\fstDom(A,A)) \ar[dr,"{\sndCoDom(1,\fst_A)}"] \\
	\sndDom(\fstCoDom(A,A),\fstDom(A,A)) \ar[ur,"{\sndDom(\fst_A,1)}"]  \ar[dr,"{\sndDom(1,\fst_A)}"']  & & & \sndCoDom(\fstDom(A,A),\fstCoDom(A,A)) \\
	& \sndDom(\fstCoDom(A,A),\fstCoDom(A,A)) \ar[r,"{\snd_{\fstCoDom(A,A)}}"] & \sndCoDom(\fstCoDom(A,A),\fstCoDom(A,A)) \ar[ur,"{\sndCoDom(\fst_A,1)}"']
	\end{tikzcd}
	\]
\end{definition}

\begin{remark}\label{rem:our definition of hc generalises natural case}
	If $F$, $G$, $H$ and $K$ all factor through the second projection $\Op\A \times \A \to \A$ or $\Op\B \times \B \to \B$, then $\phi$ and $\psi$ are just ordinary natural transformations and Definition~\ref{def:horCompDef} reduces to the usual notion of horizontal composition, Definition~\ref{def:horizontal composition natural transformations}.
\end{remark}

As in the classical natural case, we can deduce the dinaturality of $\hc \fst \snd$ from the dinaturality of $\fst$ and $\snd$, as the following Theorem states. (Recall that for $F \colon \A \to \B$ a functor, $\Op F \colon \Op\A \to \Op\B$ is the obvious functor which behaves like $F$.)

\begin{sidewaysfigure}[p]
	\centering
	\begin{tikzpicture}[every node/.style={scale=0.5}]
	\matrix[column sep=1cm, row sep=1.5cm]{
		
		\node(1) {$\H(\G(A,A),\F(A,A))$}; & & & \node(2) {$\H(\F(A,A),\F(A,A))$}; & \node(3) {$\K(\F(A,A),\F(A,A))$}; & & &
		\node(4) {$\K(\F(A,A),\G(A,A))$};\\
		
		& \node(5){$\H(\G(A,A),\F(B,A))$}; & \node(6) {$\H(\F(A,A),\F(B,A))$}; & & & \node(7) {$\K(\F(B,A),\F(A,A))$}; &
		\node (8){$\K(\F(B,A),\G(A,A))$};\\
		
		\node(9) {$\H(\G(A,B),\F(B,A))$}; & & &\node(10){$\H(\F(B,A),\F(B,A))$}; & \node(11){$\K(\F(B,A),\F(B,A))$}; & & 
		&\node(12){$\K(\F(B,A),\G(A,B))$};\\
		
		&\node(13){$\H(\G(B,B),\F(B,A))$}; & \node(14){$\H(\F(B,B),\F(B,A))$}; & & & \node(15){$\K(\F(B,A),\F(B,B))$}; & 
		\node(16){$\K(\F(B,A),\G(B,B))$};\\
		
		\node(17){$\H(\G(B,B),\F(B,B))$}; & & &\node(18){$\H(\F(B,B),\F(B,B))$}; & \node(19){$\K(\F(B,B),\F(B,B))$}; & &
		&\node(20){$\K(\F(B,B),\G(B,B))$};\\ 
	};
	
	\graph[use existing nodes,edge quotes={sloped,anchor=south}]{
		9 ->["$\H(\G(1,f),\F(f,1))$"] 1 ->["$\H(\fst_A,1)$"] 2 ->["$\snd_{\F(A,A)}$"] 3 ->["$\K(1,\fst_A)$"] 4 ->["$\K(\F(f,1),\G(1,f))$",] 12; 
		
		9 ->["$\H(\G(1,f))$"] 5 ->["$\H(\fst_A,1)$"] 6 ->["$\H(1,\F(f,1))$"] 2; 
		3 ->["$\K(\F(f,1),1)$"] 7 ->["$\K(1,\fst_A)$"] 8 ->["$\K(1,\G(1,f))$"] 12; 
		
		6 ->["$\H(\F(f,1),1)$"] 10 ->["$\snd_{\F(B,A)}$"] 11 ->["$\K(1,\F(f,1))$"] 7; 
		
		9 ->["$\H(\G(f,1),1)$"] 13 ->["$\H(\fst_B,1)$"] 14 ->["$\H(\F(1,f),1)$"] 10; 
		
		11->["$\K(1,\F(1,f))$"] 15 ->["$\K(1,\fst_B)$"] 16 ->["$\K(1,\G(f,1))$"] 12; 
		
		14->["$\H(1,F(1,f))$"] 18 ->["$\snd_{\F(B,B)}$"] 19 ->["$\K(\F(1,f),1)$"] 15;
		
		18 <-["$\H(\fst_B,1)$"] 17 <-["$\H(\G(f,1),\F(1,f))$"] 9; 
		12 <-["$\K(\F(1,f),\G(f,1))$"] 20 <-["$\K(1,\fst_B)$"] 19;
		
		1 ->[bend left=20,dashed,"$(\hc \fst \snd)_A$"] 4;
		17->[bend right=20,dashed,"$(\hc \fst \snd)_B$"] 20;
	};
	
	\path[] (9) to [out=-80,in=170] node[anchor=mid,red]{Functoriality of $\H$} (18);
	\path[] (9) to [out=80,in=-170] node[anchor=mid,red]{Functoriality of $\H$} (2);
	\path[] (19) to [out=10,in=260] node[anchor=mid,red]{Functoriality of $\K$} (12);
	\path[] (3) to [out=-10,in=-260] node[anchor=mid,red]{Functoriality of $\K$} (12);
	\path (6) to node[anchor=mid,red]{Dinaturality of $\snd$} (7);
	\path (14) to node[anchor=mid,red]{Dinaturality of $\snd$} (15);
	\path (9) to node[anchor=mid,red]{Dinaturality of $\fst$} (10);
	\path (11) to node[anchor=mid,red]{Dinaturality of $\fst$} (12);
	
	\end{tikzpicture}
	\caption{Proof of Theorem \ref{thm:horCompTheorem}: dinaturality of horizontal composition in the classical case. Here $f\colon A \to B$.}
	\label{fig:DinaturalityHorizontalCompositionFigure}
\end{sidewaysfigure}
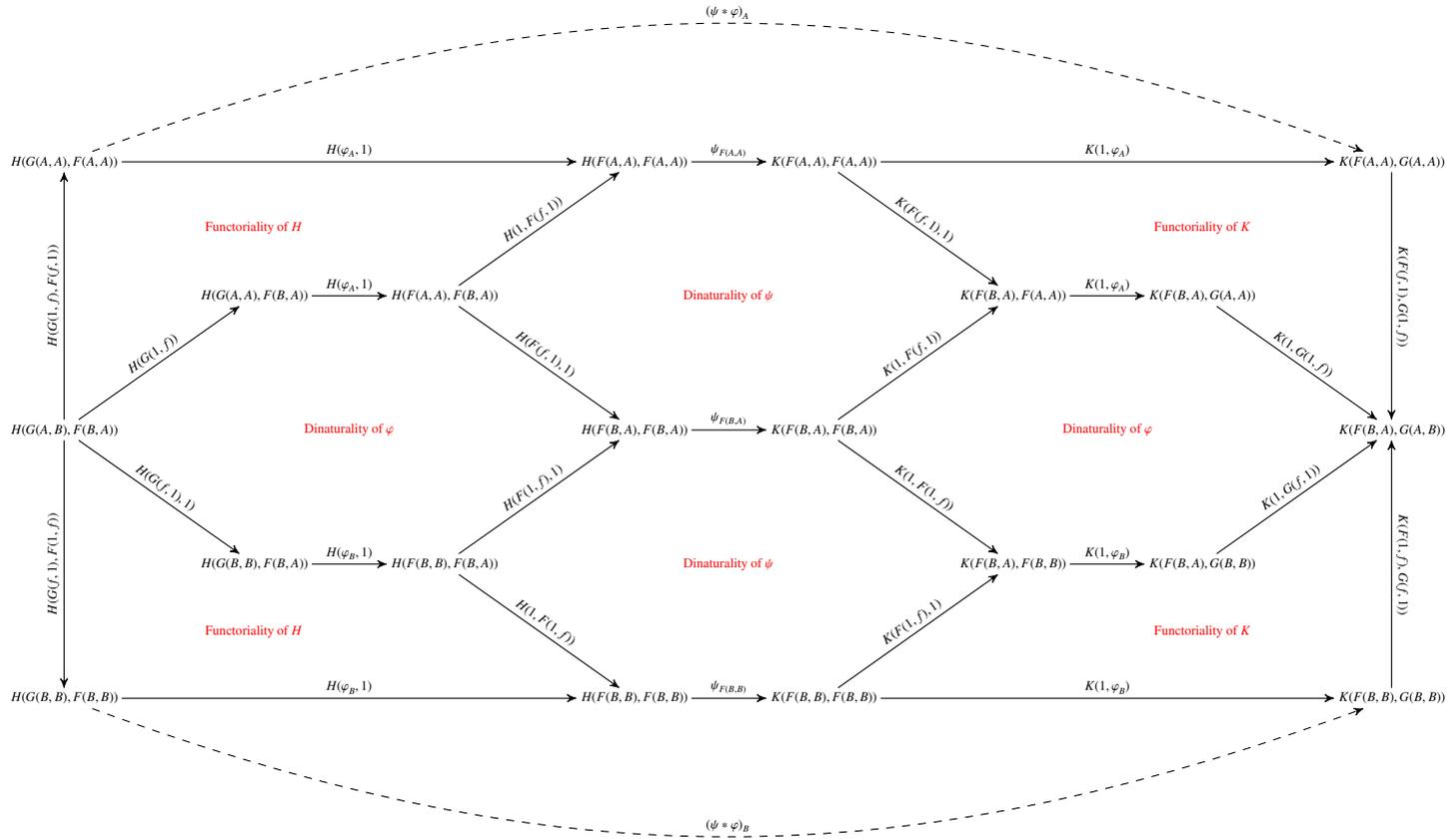

\begin{theorem}\label{thm:horCompTheorem}
	Let $\fst$ and $\snd$ be dinatural transformations as in Definition \ref{def:horCompDef}. Then $\hc \fst \snd$ is a dinatural transformation
	\[
	\hc \fst \snd \colon \sndDom(\Op\fstCoDom , \fstDom) \to \sndCoDom(\Op\fstDom,\fstCoDom)
	\]
	of type 
	$
	\begin{tikzcd}[cramped,sep=small]
	4 \ar[r] & 1 & 4 \ar[l]
	\end{tikzcd}
	$, where $\sndDom(\Op\fstCoDom , \fstDom), \sndCoDom(\Op\fstDom,\fstCoDom) \colon \A^{[+,-,-,+]} \to \C$ are defined on objects as
	\begin{align*}
	\sndDom(\Op\fstCoDom , \fstDom)(A,B,C,D) &= \hcd A B C D \\
	\sndCoDom(\Op\fstDom,\fstCoDom)(A,B,C,D) &= \hcc A B C D
	\end{align*}
	and similarly on morphisms.
\end{theorem}
\begin{proof}
	The proof consists in showing that the diagram that asserts the dinaturality of $\hc \fst \snd$ commutes: this is done in Figure~\ref{fig:DinaturalityHorizontalCompositionFigure}. \qed
\end{proof}

We can now proceed with the general definition, which involves transformations
of arbitrary type. As the idea behind Definition~\ref{def:horCompDef} is to apply
the dinaturality of $\snd$ to the general component of $\fst$ in order to define
$\hc \fst \snd$, if $\snd$ is a transformation with many variables, then we
have many dinaturality conditions we can apply to $\fst$, namely one for each
variable of $\snd$ in which $\snd$ is dinatural. Hence, the general definition
will depend on the variable of $\snd$ we want to use. For the sake of
simplicity, we shall consider only the one-category case, that is when all
functors in the definition involve one category $\C$; the general case follows with no substantial complications
except for a much heavier notation. 

\begin{definition}\label{def:generalHorizontalCompositionDef}
	Let $\F \colon \C^\fstDomVar \to \C$, $\G \colon \C^\fstCoDomVar \to \C$, $\H \colon \C^\sndDomVar \to \C$, $\K \colon \C^\sndCoDomVar \to \C$ be functors, $\fst = (\fst_\bfA)_{\bfA \in \C^n} \colon \fstDom \to \fstCoDom$ be a transformation of type
	$
	\begin{tikzcd}[cramped,sep=small]
	\length\fstDomVar \ar[r,"\fstTypL"] & \fstVarNo & \length\fstCoDomVar \ar[l,"\fstTypR"']
	\end{tikzcd}
	$
	and $\snd = (\snd_{\bf B})_{\bf B \in \C^m}\colon \sndDom \to \sndCoDom$ of type
	$
	\begin{tikzcd}[cramped,sep=small]
	\length\sndDomVar \ar[r,"\sndTypL"] & \sndVarNo & \length\sndCoDomVar \ar[l,"\sndTypR"']
	\end{tikzcd}
	$
	a transformation which is dinatural in its $i$-th variable. 
	Denoting with $\concat$ the concatenation of a family of lists, let
	\[
	\ghcDom  \colon \C^{{\concat_{u=1}^{\length\sndDomVar} \lambda^u}} \to \C, \quad \ghcCoDom \colon \C^{\concat_{v=1}^{\length\sndCoDomVar}\mu^v} \to \C
	\]	
	be functors, defined similarly to $\sndDom(\Op\fstCoDom , \fstDom)$ and  $\sndCoDom(\Op\fstDom,\fstCoDom)$ in Theorem \ref{thm:horCompTheorem}, where for all $j \in \sndVarNo$, $u \in\length\gamma$, $v\in\length\delta$:
	\[	
	\begin{tikzcd}[ampersand replacement=\&,row sep=.5em]
	\F_j= 
	\begin{cases}
	F & j=i \\
	\id\C & j \ne i
	\end{cases}
	\& \G_j= 
	\begin{cases}
	G & j=i \\
	\id\C & j \ne i
	\end{cases}
	\\
	\lambda^u = \begin{cases}
	\alpha & \eta u = i \land \gamma_u=+ \\
	\Not\beta\footnotemark & \eta u = i \land \gamma_u = - \\
	[\gamma_u] & \eta u \ne i
	\end{cases}
	\& \mu^v = \begin{cases}
	\beta & \theta v = i \land \delta_v=+ \\
	\Not\alpha & \theta v = i \land \delta_v = - \\
	[\delta_v] & \theta v \ne i
	\end{cases}
	\end{tikzcd}
	\]
	\footnotetext{Remember that for any $\beta\in\List\{+,-\}$ we denote $\Not\beta$ the list obtained from $\beta$ by swapping the signs.}Define for all $u \in\length\gamma$ and $v\in\length\delta$ the following functions:
	\[
	a_u = \begin{cases}
	\iota_n \sigma & \eta u = i \land \gamma_u=+ \\
	\iota_n \tau & \eta u = i \land \gamma_u = - \\
	\iota_m K_{\eta u} & \eta u \ne i
	\end{cases} \quad
	b_v = \begin{cases}
	\iota_n \tau & \theta v = i \land \delta_v=+ \\
	\iota_n \sigma & \theta v = i \land \delta_v = - \\
	\iota_m K_{\theta v} & \theta v \ne i
	\end{cases}
	\]
	with $K_{\eta u} \colon 1 \to m$ the constant function equal to $\eta u$, while $\iota_n$ and $\iota_m$ are defined as:
	\[
	\begin{tikzcd}[row sep=0em]
	n \ar[r,"\iota_n"] & (i-1)+n+(m-i) \\
	x \ar[|->,r] & i-1+x
	\end{tikzcd}
	\qquad
	\begin{tikzcd}[row sep=0em,ampersand replacement=\&]
	m \ar[r,"\iota_m"] \& (i-1)+n+(m-i) \\
	x \ar[|->,r] \& \begin{cases}
	x & x < i \\
	x +n -1 & x \ge i
	\end{cases}
	\end{tikzcd}
	\]
	The \emph{$i$-th horizontal composition} $\HC {[\fst]} {[\snd]} i$ is the equivalence class of the transformation 
	\[
	\HC \fst \snd i \colon \ghcDom \to \ghcCoDom
	\]	
	of type
	\[
	\begin{tikzcd}[column sep=1.5cm]
	\displaystyle\sum_{u=1}^{\length\gamma} \length{\lambda^u} \ar[r,"{[a_1,\dots, a_{\length\gamma}]}"] & (i-1) + n + (m-i) & \displaystyle\sum_{v=1}^{\length\delta} \length{\mu^v} \ar[l,"{[b_1,\dots, b_{\length\delta}]}"']
	\end{tikzcd}
	\]
	whose general component, $(\HC \fst \snd i)_{\subst B \bfA i}$, is the diagonal of the commutative hexagon obtained by applying the dinaturality of $\snd$ in its $i$-th variable to the general component $\fst_\bfA$ of $\fst$:
	\[
	\begin{tikzcd}[column sep=2em]
	& \H(\subst \bfB {F(\bfA\sigma)} i \eta) \ar[r,"\psi_{\subst \bfB {\F(\bfA\sigma)} i}"] & \K(\subst \bfB {\F(\bfA\sigma)} i \theta) \ar[dr,"\K(\substMV \bfB {\F(\bfA\sigma)} {\phi_\bfA} i \theta)"] \\
	\H(\substMV \bfB {\G(\bfA\tau)} {\F(\bfA\sigma)} i \eta) \ar[ur,"\H(\substMV \bfB {\phi_\bfA} {\F(\bfA\sigma)} i \eta)"] \ar[dr,"\H(\substMV \bfB {\G(\bfA\tau)} {\phi_\bfA} \eta)"'] \ar[rrr,dotted,"(\HC \fst \snd i)_{\subst B \bfA i}"] & & & \K(\substMV \bfB {\F(\bfA\sigma)} {\G(\bfA\tau)} i \theta) \\
	& \H(\subst \bfB {\G(\bfA\tau)} i \eta) \ar[r,"\psi_{\subst \bfB {\G(\bfA\tau)} i}"'] & \K(\subst \bfB {G(\bfA\tau)} i \theta) \ar[ur,"\K(\substMV \bfB {\phi_\bfA} {\G(\bfA\tau)} i \theta)"']
	\end{tikzcd}
	\]
\end{definition}

In other words, the domain of $\HC \fst \snd i$ is obtained by substituting the arguments of $\H$ (the domain of $\snd$) that are in the $i$-th connected component of $\graph\snd$ with $\F$ (the domain of $\fst$) if they are covariant, and with $\Op\G$ (the opposite of the codomain of $\fst$) if they are contravariant; those arguments not in the $i$-th connected component are left untouched. Similarly the codomain. The type of $\HC \fst \snd i$ is obtained by replacing the $i$-th variable of $\snd$ with all the variables of $\fst$ and adjusting the type of $\snd$ with $\fstTypL$ and $\fstTypR$ to reflect this act. In the following example, we see what happens to $\graph\fst$ and $\graph\snd$ upon horizontal composition.

\begin{example}\label{ex:hc example}
	Consider transformations $\delta$ and $\eval{}{}$ (see examples~\ref{ex:delta},\ref{ex:eval}). In the notations of Definition~\ref{def:generalHorizontalCompositionDef}, we have $F=\id\C \colon \C \to \C$, $G = \times \colon \C^{[+,+]} \to \C$, $H \colon \C^{[+,-,+]} \to \C$ defined as $H(X,Y,Z) = X \times (Y \implies Z)$ and $K = \id\C \colon \C \to \C$. The types of $\delta$ and $\eval{}{}$ are respectively
	\[
	\begin{tikzcd}[font=\small]
	1 \ar[r] & 1 & \ar[l] 2
	\end{tikzcd}
	\qquad \text{and} \qquad
	\begin{tikzcd}[row sep=0em,font=\small]
	3 \ar[r] & 2 & 1 \ar[l] \\
	1 \ar[r,|->] & 1 & 1 \ar[dl,|->,out=180,in=30] \\[-3pt]
	2 \ar[ur,|->,out=0,in=210]& 2 & \\[-3pt]
	3 \ar[ur,|->,out=0,in=210]
	\end{tikzcd}
	\]
	The transformation $\eval{}{}$ is extranatural in its first variable and natural in its second: we have two horizontal compositions. $(\HC \delta {\eval{}{}} 1 )_{A,B}$ is given by either leg of the following commutative square:
	\begin{equation}\label{delta inside eval}
		\begin{tikzcd}
			A \times \bigl( (A \times A) \implies B  \bigr) \ar[r,"\delta_A \times (1 \implies 1)"] \ar[d,"1 \times (\delta_A \implies 1)"'] & (A \times A) \times \bigl( (A \times A) \implies B \bigr) \ar[d,"\eval {A \times A} B"] \\
			A \times (A \implies B) \ar[r,"\eval A B"] & B
		\end{tikzcd}
	\end{equation}
	We have $\HC \delta {\eval{}{}} 1 \colon H(\id\C,\times,\id\C) \to \id\C(\id\C)$ where $\id\C(\id\C) = \id\C$ and
	\[
	\begin{tikzcd}[row sep=0em]
	\C^{[+,-,-,+]} \ar[r,"{H(\id\C,\times,\id\C)}"] & \C \\
	(X,Y,Z,W) \ar[|->,r] &  \quad X \times \bigl( (Y \times Z) \implies W \bigr) 
	\end{tikzcd}
	\]
	and it is of type
	\[
	\begin{tikzcd}[row sep=0em]
	4 \ar[r]      			   & 2 & \ar[l] 1 \\
	1 \ar[r,|->]  			   & 1 & 1 \ar[dl,|->,out=180,in=30]          \\[-3pt]
	2 \ar[ur,|->,out=0,in=210] & 2 \\[-3pt]
	3 \ar[uur,|->,out=0,in=230] \\[-3pt]
	4 \ar[uur,|->,out=0,in=230] \\
	\end{tikzcd}
	\]
	Intuitively, $\graph{\HC \delta {\eval{}{}} 1}$ is obtained by substituting $\graph{\delta}=\begin{tikzpicture}
	\matrix[row sep=1em,column sep=0.5em]{
		& \node (1) [category] {}; \\
		& \node (A) [component] {}; \\
		\node (2) [category] {}; & & \node (3) [category] {}; \\
	};
	\graph[use existing nodes]{
		1 -> A -> {2,3}; 
	};
	\end{tikzpicture}$
	into the first connected component of $\graph{\eval{}{}}=\begin{tikzpicture}
	\matrix[row sep=1em, column sep=1em]{
		\node (1) [category] {}; & & \node (2) [opCategory] {}; & & \node (3) [category] {}; \\
		& \node (A) [component] {}; & & & \node (B) [component] {}; \\
		& & & & \node (4) [category] {}; \\
	};
	\graph[use existing nodes]{
		1 -> A -> 2; 3 -> B -> 4;
	};
	\end{tikzpicture}$, by ``bending'', as it were, $\graph\delta$ into the $U$-turn that is the first connected component of $\graph{\eval{}{}}$:
	\[
	\begin{tikzpicture}
	\matrix[column sep=1em,row sep=2em]{
		&\node[category] (1) {}; & & \node[opCategory] (2) {}; & & \node[opCategory] (3) {}; & & \node[category] (4) {};\\
		&\node[component] (A) {};& &                           & &                           & & \node[component] (B) {};\\
		\node[coordinate] (fake1) {}; & & \node[coordinate] (fake2) {}; &\node[coordinate] (fake3) {};&  & \node[coordinate] (fake4) {};                       & & \node[category] (5) {};\\
	};
	\graph[use existing nodes]{
		1 -> A -- fake1 --[out=-90,in=-90] fake4 -> 3;
		A -- fake2 --[out=-90,in=-90] fake3 -> 2;
		4 -> B -> 5;
	};
	\end{tikzpicture}
	\quad \text{or} \quad
	\begin{tikzpicture}
	\matrix[column sep=1em,row sep=2em]{
		\node[category] (1) {};&  & \node[opCategory] (2) {}; & & \node[opCategory] (3) {}; & & \node[category] (4) {};\\
		\node[coordinate] (fake1) {}; & & & \node[component] (A) {}; & & & \node[component] (B) {};\\
		& & & & & & \node[category] (5) {};\\
	};
	\graph[use existing nodes]{
		1 -- fake1 ->[out=-90,in=-90] A -> {2,3};
		4 -> B -> 5;
	};
	\end{tikzpicture}
	\]
	Here the first graph corresponds to the upper leg of (\ref{delta inside eval})  , the second to the lower one. Notice how the component $\eval {A \times A} B$ has now \emph{two} wires, one per each $A$ in the graph on the left. The result is therefore
	\[
	\graph{\HC \delta {\eval{}{}} 1} = 
	\begin{tikzpicture}
	\matrix[column sep=1em,row sep=1.5em]{
		\node[category] (1) {}; & \node[opCategory] (2) {}; & \node[opCategory] (3) {}; & \node[category] (4) {}; \\
		& \node[component] (A) {}; & & \node[component] (B) {};\\
		& & & \node[category] (5) {};\\
	};
	\graph[use existing nodes]{
		1 -> A -> {2,3};
		4 -> B -> 5;
	};
	\end{tikzpicture}
	\]
	Turning now to the other possible horizontal composition, we have that $\HC \delta {\eval{}{}} 2 \colon H(\id\C,\id\C,\id\C) \to \id\C(\times)$ where $ H(\id\C,\id\C,\id\C) = H$ and $\id\C(\times)=\times$ by definition; it is of type
	\[
	\begin{tikzcd}[row sep=0em]
	3 \ar[r]     & 2 & \ar[l] 2 \\
	1 \ar[r,|->] & 1 & 1  \ar[dl,|->,out=180,in=30] \\[-3pt]
	2 \ar[ur,|->,out=0,in=210]& 2 & 2 \ar[l,|->]  \\[-3pt]
	3 \ar[ur,|->,out=0,in=210]
	\end{tikzcd}
	\] 
	and  $(\HC \delta {\eval{}{}} 2)_{A,B}$ is given by either leg of the following commutative square:
	\[
	\begin{tikzcd}[column sep=3em]
	A \times (A \implies B) \ar[r,"1 \times (1 \implies \delta_B)"] \ar[d,"\eval A B"'] & A \times \bigl( A \implies (B \times B) \bigr) \ar[d,"\eval A {B \times B}"] \\
	B \ar[r,"\delta_B"] & B \times B
	\end{tikzcd}
	\]
	Substituting $\graph\delta$ into the second connected component of $\graph{\eval{}{}}$, which is just a ``straight line'', results into the following graph:
	\[
	\graph{\HC \delta {\eval{}{}} 2} = 
	\begin{tikzpicture}
	\matrix[column sep=.5em,row sep=1em]{
		\node[category] (1) {}; & & \node[opCategory] (2) {}; & & \node[category] (3) {}; \\
		& \node[component] (A) {}; & & & \node[component] (B) {};\\
		& & & \node[category] (4) {}; & & \node[category] (5) {};\\
	};
	\graph[use existing nodes]{
		1 -> A -> 2;
		3 -> B -> {4,5};
	};
	\end{tikzpicture}
	\]
\end{example}

\subsection{Dinaturality of horizontal composition}\label{section dinaturality of horizontal composition}

We aim to prove here that our definition of horizontal composition, which we have already noticed generalises the well-known version for classical natural transformations (Remark~\ref{rem:our definition of hc generalises natural case}), is a closed operation on dinatural transformations. For the rest of this section, we shall fix transformations $\fst$ and $\snd$ with the notations used in Definition~\ref{def:generalHorizontalCompositionDef} for their signature; we also fix the ``names'' of the variables of $\fst$ as $\bfA=(A_1,\dots,A_n)$ and of $\snd$ as $\bfB=(B_1,\dots,B_m)$. In this spirit, $i$ is a fixed element of $\{1,\dots,m\}$, we assume $\snd$ to be dinatural in $B_i$ and we shall sometimes refer to $\HC \fst \snd i$ also as $\HC \fst \snd {B_i}$.

As in the classical natural case (Definition~\ref{def:horizontal composition natural transformations}), only the dinaturality of $\snd$ in $B_i$ is needed to \emph{define} the $i$-th horizontal composition of $\fst$ and $\snd$. Here we want to understand in which variables the $i$-th horizontal composition
\[
\HC \fst \snd {B_i} = \bigl( (\HC \fst \snd {B_i})_{\subst \bfB \bfA i}  \bigr)= \bigl( (\HC \fst \snd {B_i})_{B_1,\dots, B_{i-1},A_1,\dots, A_n, B_{i+1}, \dots, B_m} \bigr)
\] 
itself is in turn dinatural. It is straightforward to see that $\HC \fst \snd {B_i}$ is dinatural in all its $B$-variables where $\snd$ is dinatural, since the act of horizontally composing $\fst$ and $\snd$ in $B_i$ has not ``perturbed'' $\sndDom$, $\sndCoDom$ and $\snd$ in any way except in those arguments involved in the $i$-th connected component of $\graph\snd$, see example~\ref{ex:hc example}. Hence we have the following preliminary result.

\begin{proposition}
	If $\snd$ is dinatural in $B_j$, for $j \ne i$, then $\HC \fst \snd {B_i}$ is also dinatural in $B_j$.
\end{proposition}

More interestingly, it turns out that $\HC \fst \snd {B_i}$ is also dinatural in all those $A$-variables where $\fst$ is dinatural in the first place. We aim then to prove the following Theorem.

\begin{theorem}\label{thm:horCompIsDinat}
	If $\fst$ is dinatural in its $k$-th variable and $\snd$ in its $i$-th one, then $\HC \fst \snd i$ is dinatural in its $(i-1+k)$-th variable. In other words, if $\fst$ is dinatural in $A_k$ and $\snd$ in $B_i$, then $\HC \fst \snd {B_i}$ is dinatural in $A_k$.
\end{theorem}

The proof of this theorem relies on the fact that we can reduce ourselves, without loss of generality, to Theorem~\ref{thm:horCompTheorem}. To prove that, we introduce the notion of \emph{focalisation} of a transformation on one of its variables: essentially, the focalisation of a transformation $\varphi$ is a transformation depending on only one variable between functors that have only one covariant and one contravariant argument, obtained by fixing all the parts of the data involving variables different from the one we are focusing on.

\begin{definition}\label{def:focalisation def}
	Let $\varphi = (\varphi_\bfA) = (\varphi_{A_1,\dots,A_p}) \colon T \to S$ be a transformation of type 
	\[
	\begin{tikzcd}
	\length\alpha \ar[r,"\sigma"] & p & \ar[l,"\tau"'] \length\beta
	\end{tikzcd}
	\]
	with $T \colon \C^\alpha \to \C$ and $S \colon \C^\beta \to \C$. Fix $k\in\{1,\dots,p\}$ and objects $A_1,\dots,A_{k-1}$, $A_{k+1},\dots,A_p$ in $\C$. Consider functors $\bar T k$, $\bar S k \colon \Op\C \times \C \to \C$ defined by
	\begin{align*}
	\bar T k (A,B) &= T(\substMV \bfA A B i \sigma) \\
	\bar S k (A,B) &= S(\substMV \bfA A B i \tau)
	\end{align*}
	The \emph{focalisation of $\varphi$ on its $k$-th variable} is the transformation 
	\[
	\bar \varphi k \colon \bar T k \to \bar S k
	\]
	of type 
	$
	\begin{tikzcd}[cramped, sep=small]
	2 \ar[r] & 1 & \ar[l] 2
	\end{tikzcd}
	$
	where 
	\[
	\bar \varphi k_X = \phi_{\subst \bfA X i} = \varphi_{A_1\dots A_{k-1},X,A_{k+1}\dots A_p}.
	\] 
	Sometimes we may write $\bar \varphi {A_k} \colon \bar T {A_k} \to \bar S {A_k}$ too, when we fix as $A_1,\dots,A_p$ the name of the variables of $\varphi$.
\end{definition}

\begin{remark}\label{rem:focalisationIsDinaturalRemark}
	$\varphi$ is dinatural in its $k$-th variable if and only if $\bar \varphi k$ is dinatural in its only variable for all objects $A_1,\dots,A_{k-1},A_{k+1},\dots,A_p$ fixed by the focalisation of $\varphi$.
\end{remark}

The $\bar{} k$ construction depends on the $p-1$ objects we fix, but not to make the notation too heavy, we shall always call those (arbitrary) objects $A_1,\dots,A_{k-1},A_{k+1},\dots,A_n$ for $\bar \fst k$ and $B_1,\dots,B_{i-1}$, $B_{i+1},\dots,B_m$ for $\bar \snd i$.

\begin{lemma}\label{lemma:focalisationLemma}
	It is the case that $\HC \fst \snd i$ is dinatural in its $(i-1+k)$-th variable if and only if $\hc {\bar \fst k} {\bar \snd i}$ is dinatural in its only variable for all objects $B_1,\dots,B_{i-1}$, $A_1,\dots,A_{k-1}$, $A_{k+1},\dots,A_n$, $B_{i+1},\dots,B_m$ in $\C$ fixed by the focalisations of $\fst$ and $\snd$.
\end{lemma}
\begin{proof}
	The proof consists in unwrapping the two definitions and showing that they require the exact same hexagon to commute: see~\cite[Lemma 2.14]{santamaria_towards_2019}. \qed
\end{proof}

We can now prove that horizontal composition preserves dinaturality.

\begin{proofDinTheorem}
	Consider transformations $\bar \fst k$ and $\bar \snd i$. By Remark \ref{rem:focalisationIsDinaturalRemark}, they are both dinatural in their only variable. Hence, by Theorem \ref{thm:horCompTheorem}, $\hc {\bar \fst k} {\bar \snd i}$ is dinatural and by Lemma \ref{lemma:focalisationLemma} we conclude. \qed
\end{proofDinTheorem}

It is straightforward to see that horizontal composition has a left and a right unit, namely the identity (di)natural transformation on the appropriate identity functor. 

\begin{theorem}
	Let $T \colon \B^\alpha \to \C$, $S \colon \B^\beta \to \C$ be functors, and let $\varphi \colon T \to S$ be a transformation of any type. Then
	\[
	\hc \varphi {\id{\id \C}} = \varphi.
	\]
	If $\varphi$ is dinatural in its $i$-th variable, for an appropriate $i$, then also
	\[
	\HC {\id{\id \B}} \varphi i = \varphi.
	\]
\end{theorem}
\begin{proof}
	Direct consequence of the definition of horizontal composition.\qed
\end{proof}

\subsection{Associativity of horizontal composition}\label{section associativity horizontal composition}

Associativity is a crucial property of any respectable algebraic operation. In this section we show that our notion of horizontal composition is at least this respectable. 
We begin by considering classical dinatural transformations $\fst \colon \F \to \G$, $\snd \colon \H \to \K$ and $\trd \colon \U \to \V$, for $\F,\G,\H,\K,\U,\V \colon \Op\C \times \C \to \C$ functors, all of type
$
\begin{tikzcd}[cramped,sep=small]
2 \ar[r] & 1 & \ar[l] 2
\end{tikzcd}
$.

\begin{theorem}\label{thm:associativity simple case}
	$\hc {\left( \hc \fst \snd \right)} \trd = \hc \fst {\left( \hc \snd \trd \right)}$.
\end{theorem}
\begin{proof}
	We first prove that the two transformations have same domain and codomain functors. Since they both depend on one variable, this also immediately implies they have same type.
	
	We have $\hc \fst \snd \colon \H(\Op\G,\F) \to \K(\Op\F,G)$, hence 
	\[
	\hc {\left( \hc \fst \snd \right)} \trd \colon \U\Bigl(\Op{\K(\Op\F,G)},\H(\Op\G,\F)\Bigr) \to \V\Bigl( \Op{\H(\Op\G,\F)}, \K(\Op\F,G) \Bigr).
	\]
	Notice that $\Op{\K(\Op\F,G)} = \Op\K (F, \Op G)$ and $\Op{\H(\Op\G,\F)}  = \Op\H (G, \Op \F)$. Next, we have $\hc \snd \trd \colon \U(\Op\K, \H) \to \V(\Op\H, \K)$. Given that $\U(\Op\K,\H), \V(\Op\H,\K) \colon \C^{[+,-,-,+]} \to \C$, we have
	\[
	\hc \fst {\left( \hc \snd \trd \right)} \colon \underbrace{\U(\Op\K, \H)(\F,\Op\G,\Op\G,\F)}_{\U\bigl(\Op\K(\F,\Op\G),\H(\Op\G,\F)\bigr)} \to \underbrace{\V(\Op\H, \K)(\G,\Op\F,\Op\F,\G)}_{\V\bigl( \Op\H(\G,\Op\F), \K(\Op\F,\G) \bigr)}.
	\] 
	This proves $\hc {\left( \hc \fst \snd \right)} \trd$ and $\hc \fst {\left( \hc \snd \trd \right)}$ have the same signature. 
	
	Only equality of the single components is left to show. Fix then an object $A$ in $\C$. Figure~\ref{fig:Associativity} shows how to pass from $(\trd \ast \snd) \ast \fst$ to $\trd \ast (\snd \ast \fst)$ by pasting three commutative diagrams. In order to save space, we simply wrote ``$\H(\G,\F)$'' instead of the proper ``$\H(\Op\G(A,A),F(A,A))$'' and similarly for all the other instances of functors in the nodes of the diagram in Figure~\ref{fig:Associativity}; we also dropped the subscript for components of $\fst$, $\snd$ and $\trd$ when they appear as arrows, that is we simply wrote $\fst$ instead of $\fst_A$, since there is only one object involved and there is no risk of confusion. \qed
\end{proof}
\begin{sidewaysfigure}
	\footnotesize
	\begin{tikzpicture}
	\matrix[column sep=1cm, row sep=1cm]{
		\node(1) {$\U(\K(\F,\G),\H(\G,\F))$}; & & \node(2) {$\U(\K(\F,\F),\H(\F,\F))$}; & \node(3) {$\U(\H(\F,\F),\H(\F,\F))$};\\
		
		\node(4) {$\U(\K(\F,\F),\H(\G,\F))$}; & & \node(5) {$\U(\H(\F,\F),\H(\F,\F))$}; & \node(6) {$\V(\H(\F,\F),\H(\F,\F))$}; &
		\node(7) {$\V(\H(\F,\F),\K(\F,\F))$}; & & \node(8) {$\V(\H(\G,\F),\K(\F,\G))$};\\
		
		\node(9) {$\U(\H(\F,\F),\H(\G,\F))$};&&&& \node(10){$\V(\H(\G,\F),\H(\F,\F))$};&& \node(11){$\V(\H(\G,\F),\K(\F,\F))$};\\
		
		&  & \node(12){$\U(\H(\G,\F),\H(\G,\F))$}; & \node(13){$\V(\H(\G,\F),\H(\G,\F))$};\\
	};
	\graph[use existing nodes]{
		1 ->["$\U(\K(1,\fst),\H(\fst,1))$"] 2 ->["$\U(\snd,1)$"] 3 ->["$\trd$"] 6;
		1 ->["$\U(\K(1,\fst),1)$"'] 4 ->["$\U(\snd,1)$"'] 9 ->["$\U(1,\H(\fst,1))$",sloped] 5 ->["$\trd$"] 6;
		6 ->["$\V(1,\snd)$"] 7 ->["$\V(\H(\fst,1),\K(1,\fst))$"] 8;
		6 ->["$\V(\H(\fst,1),1)$",sloped] 10 ->["$\V(1,\snd)$"] 11 ->["$\V(1,\K(1,\fst))$"'] 8;
		9 ->["$\U(\H(\fst,1),1)$"',sloped] 12 ->["$\trd$"] 13 ->["$\V(1,\H(\fst,1))$"', sloped] 10;
		
	};
	
	\path (2) -- node[anchor=center,red]{\footnotesize Functoriality of $\U$} (5);
	\path (9) -- node[anchor=center,red]{\footnotesize Dinaturality of $\trd$} (10);
	\path (10) --node[anchor=center,red]{\footnotesize Functoriality of $\V$} (8);
	
	\end{tikzpicture}
	\caption{Associativity of horizontal composition in the classical case. The upper leg is $(\trd \ast \snd) \ast \fst $, whereas the lower one is $\trd \ast (\snd \ast \fst)$.}
	\label{fig:Associativity}
\end{sidewaysfigure}

We can now start discussing the general case for transformations with an arbitrary number of variables; we shall prove associativity by reducing ourselves to Theorem~\ref{thm:associativity simple case} using focalisation (see Definition~\ref{def:focalisation def}). For the rest of this section, fix transformations $\fst$, $\snd$ and $\trd$, dinatural in all their variables, with signatures:
\begin{itemize}
	\item $\fst \colon \fstDom \to \fstCoDom$, for $\fstDom \colon \C^\fstDomVar \to \C$ and $\fstCoDom \colon \C^\fstCoDomVar \to \C$, of type 
	$
	\begin{tikzcd}[cramped,sep=small]
	\length\fstDomVar \ar[r,"\fstTypL"] & \fstVarNo & \ar[l,"\fstTypR"'] \length\fstCoDomVar
	\end{tikzcd}
	$;
	\item $\snd \colon \sndDom \to \sndCoDom$, for $\sndDom \colon \C^\sndDomVar \to \C$ and $\sndCoDom \colon \C^\sndCoDomVar \to \C$, of type
	$
	\begin{tikzcd}[cramped,sep=small]
	\length\sndDomVar \ar[r,"\sndTypL"] & \sndVarNo & \ar[l,"\sndTypR"'] \length\sndCoDomVar
	\end{tikzcd}
	$;
	\item $\trd \colon \trdDom \to \trdCoDom$, for $\trdDom \colon \C^\trdDomVar \to \C$ and $\trdCoDom \colon \C^\trdCoDomVar \to \C$, of type
	$
	\begin{tikzcd}[cramped,sep=small]
	\length\trdDomVar \ar[r,"\trdTypL"] & \trdVarNo & \ar[l,"\trdTypR"'] \length\trdCoDomVar
	\end{tikzcd}
	$
\end{itemize}
For sake of simplicity, let us fix the name of the variables for $\fst$ as $\fstVariables{}{} = (A_1,\dots,A_n)$, for $\snd$ as $\sndVariables{}{} = (B_1,\dots,B_m)$ and for $\trd$ as $\trdVariables{}{} = (C_1,\dots,C_l)$. In this spirit we also fix the variables of the horizontal compositions, so for $i \in \{1,\dots,\sndVarNo\}$, the variables of $\HC \fst \snd i$ are 
\[
\sndVariables i {\fstVariables{}{}} = B_1,\dots,B_{i-1},A_1,\dots,A_n,B_{i+1},\dots,B_m
\]
and, similarly, for $j \in \{1,\dots,\trdVarNo\}$ the variables of $\HC \snd \trd j$ are 
$
\trdVariables j {\sndVariables{}{}}.
$

The theorem asserting associativity of horizontal composition, which we prove in the rest of this section, is the following.

\begin{theorem}\label{thm:associativityTheorem}
	For $i \in \{1,\dots,\sndVarNo\}$ and $j \in \{1,\dots,\trdVarNo\}$,
	\[
	\HC {\left( \HC \fst \snd i \right)} \trd j = \HC \fst {\left(\HC \snd \trd j\right)} {j-1+i}
	\]
	or, in alternative notation,
	\begin{equation}\label{eqn:associativity equation}
		\HC {\left( \HC \fst \snd {B_i} \right)} \trd {C_j} = \HC \fst {\left(\HC \snd \trd {C_j}\right)} {B_i}.
	\end{equation}
\end{theorem}

We shall require the following, rather technical, Lemma, whose proof is a matter of identity checking.

\begin{lemma}\label{lemma:associativity techincal lemma}
	Let $\Phi = (\Phi_{V_1,\dots,V_p})$ and $\Psi = (\Psi_{W_1,\dots,W_q})$ be transformations in $\C$ such that $\Psi$ is dinatural in $W_s$, for $s \in \{1,\dots,q\}$. Let $V_1,\dots,V_{r-1}$, $V_{r+1},\dots,V_p$, $W_1,\dots,W_{s-1}$, $W_{s+1},\dots,W_q$ be objects of $\C$, and let $\bar \Phi {V_r}$ and $\bar \Psi {W_s}$ be the focalisation of $\Phi$ and $\Psi$ in its $r$-th and $s$-th variable respectively using the fixed objects above. Let also $X$ be an object of $\C$. Then 
	\begin{enumerate}[(i)]
		\item $ \left( \hc { \bar \Phi {V_r} } { \bar \Psi {W_s} } \right)_X = \left( \HC \Phi \Psi {W_s} \right)_{W_1,\dots,W_{s-1},V_1,\dots,V_{r-1},X,V_{r+1},\dots,V_p,W_{s+1},\dots,W_q} = \left( \bar {\HC \Phi \Psi {W_s}} {V_r} \right)_X $
		\item $\mathit{(co)dom}\left( \bar {\HC \Phi \Psi {W_s}} {V_r} \right) (x,y) = \mathit{(co)dom}\left( \hc {\bar \Phi {V_r}} {\bar \Psi {W_s}} \right) (x,y,y,x)  $ for any morphisms $x$ and $y$.
	\end{enumerate}
\end{lemma}

\begin{remark}\label{rem:associativity techincal lemma remark}
	Part (i) asserts an equality between \emph{morphisms} and not \emph{transformations}, as $ \hc { \bar \Phi {V_r} } { \bar \Psi {W_s} }$ and $\HC \Phi \Psi {W_s}$ have different types and even different domain and codomain functors.
\end{remark}

\begin{proofAssociativity}
	One can show that $\HC {\left( \HC \fst \snd {B_i} \right)} \trd {C_j}$ and $ \HC \fst {\left(\HC \snd \trd {C_j}\right)} {B_i}$ have the same domain, codomain and type simply by computing them and observing they coincide. In particular, notice that they both depend on the following variables: $\trdVariables j {\sndVariables i {\fstVariables{}{}}}$. Here we show that their components are equal. Let us fix then $
	C_1,\dots, C_{j-1}$, $B_1$, $\dots$, $B_{i-1}$, $A_1$, $\dots$, $A_{k-1}$, $X$, $A_{k+1}$, $\dots$, $A_n$, $B_{i+1}$, $\dots$, $B_m$, $C_{j+1}$, $\dots$, $C_l$ objects in $\C$. Writing just $V$ for this long list of objects, we have, by Lemma~\ref{lemma:associativity techincal lemma}, that
	\[
	\left(\HC \fst {\left(\HC \snd \trd {C_j}\right)} {B_i}\right)_V = \left( \hc {\bar \fst {A_k}} {\bar {\HC \snd \trd {C_j}} {B_i}} \right)_X .
	\]
	Now, we cannot apply again Lemma~\ref{lemma:associativity techincal lemma} to $\bar {\HC \snd \trd {C_j}} {B_i}$ because of the observation in Remark~\ref{rem:associativity techincal lemma remark}, but we can use the definition of horizontal composition to write down explicitly the right-hand side of the equation above: it is the morphism
	\[
	\codom{ \bar {\HC \snd \trd {C_j}} {B_i} } (\id{\bar \F {} (X,X)} , (\bar \fst {A_k})_X) \circ
	\left( \bar{\HC \snd \trd {C_j}}{B_i} \right)_{\bar \F {} (X,X)} \circ
	\dom{ \bar {\HC \snd \trd {C_j}} {B_i} }( {(\bar \fst {A_k})}_X ,\id{\bar \F {} (X,X)})
	\]
	(Remember that $\bar \fst {A_k} \colon \bar \F {A_k} \to \bar \G {A_k}$, here we wrote $\bar F {} (X,X)$ instead of $\bar F {A_k}(X,X)$ to save space.) Now we \emph{can} use Lemma~\ref{lemma:associativity techincal lemma} to ``split the bar'', as it were:
	\begin{multline*}
	\codom{\hc {\bar \snd {B_i}} {\bar \trd {C_j}}} \bigl( {(\bar \fst {A_k})}_X, \id{\bar \F {} (X,X)}, \id{\bar \F {} (X,X)}, {(\bar \fst {A_k})}_X \bigr) \circ \\[.5em]
	\left( \hc {\bar \snd {B_i}} {\bar \trd {C_j}} \right)_{\bar \F {} (X,X)} \circ \\[.5em]
	\dom{\hc {\bar \snd {B_i}} {\bar \trd {C_j}}} \bigl(\id{\bar \F {} (X,X)}, {(\bar \fst {A_k})}_X, {(\bar \fst {A_k})}_X, \id{\bar \F {} (X,X)}\bigr) 
	\end{multline*}
	This morphism is equal, by definition of horizontal composition, to
	\[
	\left( \hc {\bar \fst {A_k}} {\left( \hc {\bar \snd {B_i}} {\bar \trd {C_j}} \right)} \right)_X
	\]
	which, by Theorem~\ref{thm:associativity simple case}, is the same as
	\[
	\left( \hc {\left( \hc {\bar \fst {A_k}} {\bar \snd {B_i}} \right)} {\bar \trd {C_j}} \right)_X.
	\]
	An analogous series of steps shows how this is equal to $\left( \HC {\left(\HC \fst \snd {B_i}\right)} \trd {C_j}  \right)_V$, thus concluding the proof. \qed
\end{proofAssociativity}

\subsection{Whiskering and horizontal composition}

Let $\phi \colon F \to G$, with $F \colon \C^\alpha \to \C$ and $G \colon \C^\beta \to \C$, and $\psi \colon H \to K$, with $H,K \colon \Op\C \times \C\to \C$, be dinatural transformations of type $\begin{tikzcd}[cramped,sep=small]
\length\alpha \ar[r,"\sigma"] & n & \length\beta \ar[l,"\tau"']
\end{tikzcd}$
and
$\begin{tikzcd}[cramped,sep=small]
2 \ar[r] & 1 & \ar[l] 2
\end{tikzcd}$
respectively. Then $\hc \phi \psi \colon H(\Op G , F) \to K(\Op F, G)$ is of type
\[
\begin{tikzcd}
\length\beta + \length\alpha \ar[r,"{[\tau,\sigma]}"] & n & \ar[l,"{[\sigma,\tau]}"'] \length\alpha + \length\beta
\end{tikzcd}
\] 
and its general component $(\hc \phi \psi)_{\bfA}$, with $\bfA=(A_1,\dots,A_n)$, is either leg of
\[
\begin{tikzcd}[column sep=2em]
& H\bigl( F(\bfA\sigma), F(\bfA\sigma) \bigr) \ar[r,"\psi_{F(\bfA\sigma)}"] & K \bigl( F(\bfA\sigma), F(\bfA\sigma) \bigr) \ar[dr,"{K\bigl( F(\bfA\sigma),\phi_\bfA \bigr)}"] \\
H \bigl( G(\bfA\tau), F(\bfA\sigma) \bigr) \ar[ur,"{H \bigl( \phi_\bfA , F(\bfA\sigma) \bigr)}"] \ar[dr,"{ H \bigl( G(\bfA\tau), \phi_\bfA \bigr)}"'] & & & K \bigl( F(\bfA\sigma), G(\bfA\tau) \bigr) \\
& H \bigl( G(\bfA\tau), G(\bfA\tau) \bigr) \ar[r,"\psi_{G(\bfA\tau)}"'] & K \bigl( G(\bfA\tau), G(\bfA\tau) \bigr) \ar[ur,"{K \bigl( \phi_\bfA, G(\bfA\tau) \bigr)}"']
\end{tikzcd}
\]
If we look at the upper leg of the above hexagon, we may wonder: is it, in fact, the general component of the vertical composition
\begin{equation}\label{eqn:whiskering vs horizontal composition}
	\bigl(\hc {(F,\phi)} K \bigr) \circ \bigl( \hc F \psi \bigr) \circ \bigl( \hc {(\phi,F)} H \bigr),
\end{equation}
where by $\hc {(F,\phi)} K$ we mean the simultaneous horizontal composition of $\id K$ with $\id F$ in its first variable and with $\phi$ in its second, namely $\bigl(\HC {\id F} {\id K} 1 \bigr) \circ \bigl( \HC \phi {\id K} 2 \bigr) = \bigl( \HC \phi {\id K} 2 \bigr) \circ \bigl( \HC {\id F} {\id K} 1 \bigr)$? In other words, is horizontal composition a vertical composite of \emph{whiskerings}, analogously to the classical natural case? \emph{No}, but it is intimately related to it, as we show now by computing the composite~\eqref{eqn:whiskering vs horizontal composition}. Let $\bfA=(A_1,\dots,A_n)$,  $\bfB=(B_1,\dots,B_{\length\alpha})$ and $\bfC=(C_1,\dots,C_{\length\alpha})$. Then
\begin{multline*}
\bigl(\hc {(\phi,F)} H \bigr)_{\bfA \concat \bfB} =
\begin{tikzcd}[ampersand replacement=\&,column sep=4em]
H \bigl( G(\bfA\tau), F(\bfB) \bigr) \ar[r,"{H\bigl( \phi_\bfA, F(\bfB)\bigr)}"] \& H \bigl( F(\bfA\sigma), F(\bfB) \bigr)
\end{tikzcd} \\
(\hc F \psi)_\bfC =
\begin{tikzcd}[ampersand replacement=\&,column sep=4em]
H \bigl( F(\bfC), F(\bfC) \bigr) \ar[r,"\psi_{F(\bfC)}"] \& K \bigl( F(\bfC), F(\bfC) \bigr)
\end{tikzcd}
\end{multline*}
Therefore, upon composing $\hc {(\phi,F)} H $ with $\hc F \psi$, we have to impose $\bfA\sigma = \bfC = \bfB$, which means $A_{\sigma_i} = C_i = B_i$ for all $i \in \length\alpha$. If we take also $\bfD=(D_1,\dots,D_{\length\alpha})$ and $\bfE=(E_1,\dots,E_n)$, we have:
\begin{multline*}
\Bigl(\bigl( \hc F \psi \bigr) \circ \bigl( \hc {(\phi,F)} H \bigr)\Bigr)_{\bfA} =
\begin{tikzcd}[ampersand replacement=\&,column sep=6em]
H \bigl( G(\bfA\tau), F(\bfA\sigma) \bigr) \ar[r,"{\psi_{F(\bfA\sigma)} \circ H \bigl( \phi_\bfA, F(\bfA\sigma)\bigr)}"] \& K \bigl( F(\bfA\sigma), F(\bfA\sigma) \bigr)
\end{tikzcd} \\
\bigl(\hc {(F,\phi)} K \bigr)_{\bfD \concat \bfE} =
\begin{tikzcd}[ampersand replacement=\&,column sep=3em]
K \bigl( F(\bfD), F(\bfE\sigma) \bigr) \ar[r,"{K\bigl( F(\bfD), \phi_\bfE \bigr)}"] \& K \bigl( F(\bfD), G(\bfE\tau) \bigr)
\end{tikzcd}
\end{multline*}
So, to compose $\bigl( \hc F \psi \bigr) \circ \bigl( \hc {(\phi,F)} H \bigr)$ with $\bigl(\hc {(F,\phi)} K \bigr)$, we only need $\bfD = \bfA\sigma = \bfE\sigma$. Crucially, for all $k \in n \setminus \mathrm{Im}(\sigma)$, $A_k$ and $E_k$ need not to be equal. This means that, if $l=\length{n \setminus \mathrm{Im}(\sigma)}$, the composite~\eqref{eqn:whiskering vs horizontal composition} is a transformation depending on $l + 2 \cdot (n-l)$ variables (equivalently, $n+(n-l)$ variables), where the two instances of $\phi$, at the beginning and at the end of the general component of~\eqref{eqn:whiskering vs horizontal composition}, are computed in variables $\bfA$ and $\bfE$ respectively where $A_{\sigma i} = E_{\sigma i}$ for all $i \in \length\alpha$, and $A_k \ne E_k$ for $k \in n \setminus \mathrm{Im}(\sigma)$. Instead, the horizontal composition $\hc \phi \psi$ uses \emph{the same} general component of $\phi$ in both instances, as it is obtained by ``applying the dinaturality condition of $\psi$ to $\phi_\bfA$'': it requires a stronger constraint, namely $A_k = E_k$ for all $k\in n$. This means that $\hc \phi \psi$ is a sub-family of~\eqref{eqn:whiskering vs horizontal composition}, and it coincides with it precisely when $\sigma$ is surjective. Analogously, $\hc \phi \psi$ is in general a sub-family of
\[
\bigl(\hc {(\phi,G)} K \bigr) \circ \bigl( \hc G \psi \bigr) \circ \bigl( \hc {(G,\phi)} H \bigr)
\]
and they coincide precisely when $\tau$ is surjective. This issue is part of the wider problem of incompatibility of horizontal and vertical composition, which we discuss in the next section.

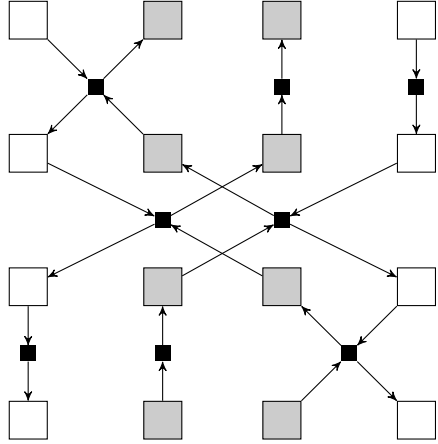
\begin{figure}
	\centering
	\begin{tikzpicture}
	\matrix[row sep=1.5em,column sep=1.5em]{
		\node[category] (1) {}; & & \node[opCategory] (2) {}; & & \node[opCategory] (3) {}; & & \node[category] (4) {}; \\
		& \node[component] (5) {}; & & & \node[component] (6) {}; & & \node[component] (7) {};\\
		\node[category] (8) {}; & & \node[opCategory] (9) {}; & & \node[opCategory] (10) {}; & & \node[category] (11) {}; \\
		& &\node[component] (12) {};  & & \node[component] (13) {};\\
		\node[category] (14) {}; & & \node[opCategory] (15) {}; & & \node[opCategory] (16) {}; & & \node[category] (17) {}; \\
		\node[component] (18) {}; & & \node[component] (19) {}; & & & \node[component] (20) {};\\
		\node[category] (21) {}; & & \node[opCategory] (22) {}; & & \node[opCategory] (23) {}; & & \node[category] (24) {}; \\
	};
	\graph[use existing nodes]{
		1 -> 5 -> 8 -> 12 -> 14 -> 18 -> 21;
		22 -> 19 -> 15 -> 13 -> 9 -> 5 -> 2;
		23 -> 20 -> 16 -> 12 -> 10 -> 6 -> 3;
		4 -> 7 -> 11 -> 13 -> 17 -> 20 -> 24;
	};
	\end{tikzpicture}
	\caption{The composite graph of the three whiskerings in~\eqref{eqn:whiskering vs horizontal composition} for $\phi$ and $\psi$ classical dinatural transformations.}
	\label{fig:composite graph whiskerings}
\end{figure}

\begin{remark}
	If $\phi$ and $\psi$ are classical dinatural transformations as in Definition~\ref{def:horCompDef}, $\hc \phi \psi$ is indeed equal to the vertical composite of whiskerings~\eqref{eqn:whiskering vs horizontal composition}. In this case, one can alternatively infer the dinaturality of $\hc \phi \psi$ (Theorem~\ref{thm:horCompTheorem}) from the easy-to-check dinaturality of $\hc {(\phi,F)} H$, $\hc F \psi$ and $\hc {(F,\phi)} K$ by applying Theorem~\ref{theorem:compositionality with complicated graphs}: one can draw the composite graph of three whiskerings (Figure~\ref{fig:composite graph whiskerings}) and notice that the resulting Petri Net is acyclic. The trellis of commutative diagrams in Figure~\ref{fig:DinaturalityHorizontalCompositionFigure}, which proves Theorem~\ref{thm:horCompTheorem}, is the algebraic counterpart of firing transitions in Figure~\ref{fig:composite graph whiskerings}: the dinaturality of $\phi$ corresponds to firing the top-left and bottom-right transitions, while the dinaturality of $\psi$ to firing the two transitions in the middle layer. 
\end{remark}

\subsection{(In?)Compatibility with vertical composition}\label{section compatibility}

Looking at the classical natural case, there is one last property to analyse: the \emph{interchange law}~\cite{mac_lane_categories_1978}. In the following situation,
\[
\begin{tikzcd}[column sep=1.5cm]
\A \arrow[r, out=60, in=120, ""{name=U, below}]
\arrow[r, ""{name=D, }]
\arrow[r,phantom,""{name=D1,below}]
\arrow[r, bend right=60,""{name=V,above}]
& \B \arrow[r, bend left=60, ""{name=H, below}]
\arrow[r,""{name=E}]
\arrow[r,phantom,""{name=E1,below}]
\arrow[r, bend right=60,""{name=K,above}]
& \C
\arrow[Rightarrow, from=U, to=D,"\phi"]
\arrow[Rightarrow, from=D1, to=V,"\psi"]
\arrow[Rightarrow, from=H, to=E, "\phi'"]
\arrow[Rightarrow, from=E1,to=K, "\psi'"]
\end{tikzcd}
\]
with $\phi,\phi',\psi$ and $\psi'$ natural transformations, we have:
\begin{equation}\label{interchange law}
	\hc {(\psi \circ \phi)} {(\psi' \circ \phi')} = (\hc \psi {\psi'}) \circ (\hc \phi {\phi'}). \tag{$\dagger$}
\end{equation}
The interchange law is the crucial property that makes $\C\mathrm{at}$ a 2-category. It is then certainly highly interesting to wonder whether a similar property holds for the more general notion of horizontal composition for dinatural transformations too.

As we know all too well, dinatural transformations are far from being as well-behaved as natural transformations, given that they do not, in general, vertically compose; on the other hand, their horizontal composition always works just fine. Are these two operations compatible, at least when vertical composition is defined?
The answer, unfortunately, is \emph{No}, at least if by ``compatible'' we mean ``compatible as in the natural case (\ref{interchange law})''. Indeed, consider classical dinatural transformations
\begin{equation}\label{compatibility situation}
	\begin{tikzcd}[column sep=0.75cm]
		\Op\A \times \A \arrow[rr, bend left=60, ""{name=U,below}]
		\arrow[rr, phantom, bend left=60, "F"{above}]
		\arrow[rr, "G"{name=D,anchor=center,fill=white,pos=0.34}]
		\arrow[rr, bend right=60,""{name=V,above}]
		\arrow[rr, bend right=60,"H"{below}]
		&  & \B
		& \Op\B \times \B \arrow[rr, bend left=60, ""{name=H, below,}]
		\arrow[rr,phantom, bend left=60, "J"{above}]
		\arrow[rr,"K"{name=E,anchor=center,fill=white},pos=0.35]
		\arrow[rr, bend right=60,""{name=K,above}]
		\arrow[rr,phantom, bend right=60,"L"{below}]
		&  & \C
		\arrow[Rightarrow, from=U, to=D,"\phi"]
		\arrow[Rightarrow, from=D, to=V,"\psi"]
		\arrow[Rightarrow, from=H, to=E, "\phi'"]
		\arrow[Rightarrow, from=E,to=K, "\psi'"]
	\end{tikzcd}
\end{equation}
such that $\psi\circ\phi$ and $\psi'\circ\phi'$ are dinatural. Then 
\[
\hc \phi {\phi'} \colon J(\Op G,F) \to K(\Op F,G) \qquad
\hc \psi {\psi'} \colon K(\Op H,G) \to L(\Op G,H)
\]
which means that $\hc \phi {\phi'}$ and $\hc \psi {\psi'}$ are not even composable
as families of morphisms, as the codomain of the former is not the domain of the
latter. The problem stems from the fact that the codomain of the horizontal composition $\hc \phi {\phi'}$ depends on the codomain of $\phi'$ and also the domain \emph{and} codomain of $\phi$, which are not the same as the domain and codomain of $\psi$: indeed, in order to be vertically composable, $\phi$ and $\psi$ must share only one functor, and not both. This does not happen in the natural case: the presence of mixed variance, which forces to consider the codomain of $\phi$ in $\hc \phi {\phi'}$ and so on, is the real culprit here. 

The failure of (\ref{interchange law}) is not completely unexpected: after all, our definition of horizontal composition is strictly more general than the classical one for natural transformations, as it extends the audience of functors and transformations it can be applied to quite considerably. Hence it is not surprising that this comes at the cost of losing one of its properties, albeit so desirable. Of course, one can wonder whether a different definition of horizontal composition exists for which (\ref{interchange law}) holds. Although we cannot exclude \emph{a priori} this possibility, the fact that ours not only is a very natural generalisation of the classical definition for natural transformations (as it follows the same idea, see discussion after Definition~\ref{def:horizontal composition natural transformations}), but also enjoys associativity and unitarity, let us think that we \emph{do} have the right definition at hand. (As a side point, behold Figure~\ref{fig:DinaturalityHorizontalCompositionFigure}: its elegance cannot be the fruit of a wrong definition!)

What we suspect, instead, is that a different \emph{interchange law} should be formulated, that can accommodate the hexagonal shape of the dinatural condition. Indeed, what proves (\ref{interchange law}) in the natural case is the naturality of either $\phi'$ or $\psi'$. For instance, the following diagrammatic proof uses the latter, for $\phi \colon F \to G$, $\psi \colon G \to H$, $\phi' \colon J \to K$, $\psi' \colon K \to L$ natural:
\[
\begin{tikzcd}[row sep=2.5em,column sep=2.5em,font=\small]
JF(A) \ar[r,"{\phi'_{F(A)}}"] \ar[rd,dashed,"{(\hc \phi {\phi'})_A}"'] & KF(A) \ar[r,"\psi'_{F(A)}"] \ar[d,"K(\phi_A)"'] & LF(A) \ar[d,"L(\phi_A)"] \\
& KG(A) \ar[r,"\psi'_{G(A)}"] \ar[dr,dashed,"{(\hc \psi {\psi'})}_A"'] & LG(A) \ar[d,"L(\psi_A)"] \\
& & LH(A)
\end{tikzcd}
\]
(The upper leg of the diagram is $	\hc {(\psi \circ \phi)} {(\psi' \circ \phi')}$.) The naturality condition of $\psi'$ is what causes $\phi$ and $\psi'$ to swap places, allowing now $\phi$ and $\phi'$ to interact with each other via horizontal composition; same for $\psi$ and $\psi'$. 

However, for $\phi, \psi, \phi',\psi'$ dinatural as in (\ref{compatibility situation}), this does not happen:
\[
\begin{tikzcd}[column sep=.5cm,font=\small]
& & J(F,F) \ar[r,"\phi'"] & K(F,F) \ar[r,"\psi'"] & L(F,F) \ar[dr,"{L(1,\phi)}"] \\
& J(G,F) \ar[ur,"{J(\phi,1)}"] \ar[rrrr,dashed,"\hc {\phi} {({\psi'}\circ{\phi'})}"] & & & & L(F,G) \ar[dr,"{L(1,\psi)}"]\\
J(H,F) \ar[ur,"{J(\psi,1)}"] & & & & & & L(F,H)
\end{tikzcd}
\]
Here, the upper leg of the diagram is again $\hc {(\psi \circ \phi)} {(\psi' \circ \phi')}$; we have dropped the lower-scripts of the transformations and we have written ``$J(H,F)$'' instead of ``$J(H(A,A),F(A,A))$'' to save space. The dinaturality conditions of $\phi'$ and $\psi'$ do not allow a place-swap for $\phi$ and $\phi'$ or for $\phi$ and $\psi'$; in fact, they cannot be applied at all! The only thing we can notice is that we can isolate $\phi$ from $\phi'$, obtaining the following:
\[
\hc {(\psi \circ \phi)} {(\psi' \circ \phi')} = L(1,\psi) \circ \Bigl(\hc {\phi} {({\psi'}\circ{\phi'})}\Bigr) \circ J(\psi,1).
\]
Notice that the right-hand side is \emph{not} $\hc \psi {\Bigl(\hc {\phi} {({\psi'}\circ{\phi'})}\Bigr)}$, as one might suspect at first glance, simply because the domain of $\hc {\phi} {({\psi'}\circ{\phi'})}$ is not $J$ and its codomain is not $L$. 

It is clear then that the only assumption of $\phi'\circ\phi$ and $\psi'\circ\psi$ being dinatural (for whatever reason) is not enough. One chance of success could come from involving the graph of our transformations; for example, if the composite graphs $\graph\psi \circ \graph\phi$ and $\graph{\psi'} \circ \graph{\phi'}$ are acyclic—hence dinatural, yes, but for a ``good'' reason—then maybe we could be able to deduce a suitably more general, ``hexagonal'' version of (\ref{interchange law}) for dinatural transformations. It also may well be that there is simply no sort of interchange law, of course. This is still an open question, and the matter of further study in the future. In the conclusions we shall make some additional comments in light of the calculus we will build in the rest of the article.

\section{A category of partial dinatural transformations}

Since dinatural transformations do not always compose, they do not form a category. However, the work done in Section~\ref{section vertical compositionality} permits us to define a category whose objects are functors of mixed variance and whose morphisms are transformations that are dinatural only in \emph{some} of their variables, as we shall see. A first attempt would be to construct $\fc \B \C$ by defining:
\begin{itemize}\label{first attempt}
	\item objects:  pairs $(\alpha, F \colon \B^\alpha \to \C)$;
	\item morphisms: a morphism $(\alpha, F) \to (\beta, G)$ would be a tuple  $(\phi, \graph\phi, \Delta_\phi)$ where $\phi \colon F \to G$ is a transformation whose standard graph is $\graph\phi$, and if $n$ is the number of connected components of $\graph\phi$ (hence, the number of variables of $\phi$), then $\Delta_\phi \colon n \to \{0,1\}$ would be the ``discriminant'' function that tells us in which variables $\phi$ is dinatural: if $\Delta_\phi(i)=1$, then $\phi$ is dinatural in its $i$-th variable;
	\item composition: given $(\phi,\graph\phi,\Delta_\phi) \colon (\alpha,F) \to (\beta,G)$ and $(\psi,\graph\psi,\Delta_\psi) \colon (\beta,G) \to (\gamma,H)$ morphisms, their composite would be $(\psi\circ\phi,\graph{\psi\circ\phi},\Delta_{\psi\circ\phi})$, where $\psi\circ\phi$ is simply the vertical composition of transformations $\phi$ and $\psi$, $\graph{\psi\circ\phi}$ is its standard graph, and $\Delta_{\psi\circ\phi}(x)$ is defined to be $1$ if and only if the $x$-th connected component of $\graph\psi \circ \graph\phi$ is acyclic and $\phi$ and $\psi$ are dinatural in all variables involved in the $x$-th connected component of the composite graph $\graph{\psi}\circ\graph\phi$, in the sense of Theorem~\ref{theorem:acyclicity implies dinaturality GENERAL}.
\end{itemize}
However, composition so defined fails to be associative in $\Delta$. Suppose we have three consecutive transformations $\phi$, $\psi$ and $\chi$, dinatural in all their variables, where
\[
\graph\phi = 
\begin{tikzpicture}
\matrix[column sep=2.4mm,row sep=0.4cm]{
	& \node [category] (1) {}; & & & & \node [opCategory] (6)  {}; \\
	& \node [component] (A)  {}; & & & & \node [component] (B) {};\\
	\node [category] (2) {}; & & \node [category] (3){}; & & \node [opCategory] (5)  {}; & & \node [opCategory] (7)  {};\\
};
\graph[use existing nodes]{
	1 -> A -> {2,3};
	5 -> B -> 6;
	7 -> B;
};
\end{tikzpicture}
\quad
\graph\psi =
\begin{tikzpicture}
\matrix[column sep=2.4mm,row sep=0.4cm]{
	\node [category] (2) {}; & & \node [category] (3){}; & & \node [opCategory] (5)  {}; & & \node [opCategory] (7)  {};\\
	\node [component] (C) {}; & & & \node [component] (D) {}; & & &\node [component] (E) {};	\\
	\node [category] (4) {}; & & & & & & \node [opCategory] (8)  {};\\
};
\graph[use existing nodes]{
	2 -> C -> 4;
	3 -> D -> 5;
	8 -> E -> 7;
};
\end{tikzpicture} \quad 
\graph\chi =
\begin{tikzpicture}

\matrix[column sep=2.4mm,row sep=0.4cm]{
	\node[category](2){}; & & & & \node[opCategory](5){}; \\
	& & \node[component](B){};\\
};
\graph[use existing nodes]{
	2 -> B -> 5;
};
\end{tikzpicture}
\]
Of course vertical composition of transformations \emph{is} associative, therefore $(\chi \circ \psi) \circ \phi = \chi \circ (\psi \circ \phi)$ and $\graph{(\chi \circ \psi) \circ \phi} = \graph{\chi \circ (\psi \circ \phi)}$. Yet, $\Delta_{(\chi \circ \psi) \circ \phi} \ne \Delta_{\chi \circ (\psi \circ \phi)}$: indeed, by computing $\graph\chi \circ \graph\psi$ and then collapsing the connected components, we obtain
\[
\graph{\chi\circ\psi} =
\begin{tikzpicture}
\matrix[column sep=2.4mm,row sep=0.4cm]{
	\node [category] (1) {}; & \node[category] (7) {}; & & \node[opCategory] (8) {}; & \node[opCategory] (6) {}; \\
	& & \node[component](D){}; \\
	& & \node[component](B){};\\
};
\graph[use existing nodes]{
	1 -> B -> 6;
	7 -> D -> 8;
};
\end{tikzpicture}
\quad \text{hence }
\graph{\chi \circ \psi} \circ \graph\phi =
\begin{tikzpicture}
\matrix[column sep=2.4mm,row sep=0.4cm]{
	&\node[category](1){}; & & & &\node[opCategory](6){};\\
	&\node[component](A){};& & & &\node[component](C){};\\
	\node[category](2){}; & &\node[category](7){}; & &\node[opCategory](8){}; & &\node[opCategory](5){};\\
	& & &\node[component](D){};\\
	& & &\node[component](B){};\\
};
\graph[use existing nodes]{
	1 -> A -> 2 -> B ;
	B -> 5 -> C -> 6;
	A -> 7 -> D -> 8 -> C;
};
\end{tikzpicture}
\]
Since $\graph{\chi \circ \psi} \circ \graph\phi$ is acyclic, we have that $(\chi\circ\psi)\circ\phi$ is dinatural, thus $\Delta_{(\chi\circ\psi)\circ\phi} \colon 1 \to \{0,1\}$ is the function returning 1. On the other hand, however, we have
\[
\graph\psi \circ \graph\phi =
\begin{tikzpicture}
\matrix[column sep=2.4mm,row sep=0.4cm]{
	& \node [category] (1) {}; & & & & \node [opCategory] (6)  {}; \\
	& \node [component] (A)  {}; & & & & \node [component] (B) {};\\
	\node [category] (2) {}; & & \node [category] (3){}; & & \node [opCategory] (5)  {}; & & \node [opCategory] (7)  {};\\
	\node [component] (C) {}; & & & \node [component] (D) {}; & & &\node [component] (E) {};	\\
	\node [category] (4) {}; & & & & & & \node [opCategory] (8)  {};\\
};

\graph[use existing nodes]{
	1 -> A -> {2,3};
	2 -> C -> 4;
	3 -> D -> 5 -> B -> 6;
	8 -> E -> 7 -> B;
};
\end{tikzpicture}
\quad \text{so } 
\graph{\psi\circ\phi} = 
\begin{tikzpicture}
\matrix[column sep=3.5mm,row sep=0.4cm]{
	\node [category]  (1) {}; & & \node [opCategory] (2) {}; \\
	& \node [component]  (A) {}; \\
	\node [category]  (3) {}; & & \node [opCategory]  (4) {}; \\
};
\graph[use existing nodes]{
	1 -> A -> 3; 4 -> A -> 2;
};
\end{tikzpicture}
\]
which means that, when we glue together $\graph\chi$ and $\graph{\psi\circ\phi}$, we obtain:
\[
\graph{\chi}\circ\graph{\psi\circ\phi}= 
\begin{tikzpicture}
\matrix[column sep=3.5mm,row sep=0.4cm]{
	\node[category](1){}; & & \node[opCategory](6){};\\
	&\node[component](A){};\\
	\node[category](2){}; & & \node[opCategory](5){};\\
	&\node[component](B){};\\
};
\graph[use existing nodes]{
	1->A->2->B->5->A->6;	
};
\end{tikzpicture} 
\]
which is cyclic, so $\Delta_{\chi\circ(\psi\circ\phi)} \colon 1 \to \{0,1\}$ returns 0.

What went wrong? In the graph of $\psi\circ\phi$ there is a path from the bottom-right node to the bottom-left node, which then extends to a cycle once connected to $\graph{\chi}$. That path was created upon collapsing the composite graph $\graph\psi \circ \graph\phi$ into $\graph{\psi \circ \phi}$: 
but in $\graph\psi \circ \graph\phi$ there was no path from the bottom-right node to the bottom-left one. And rightly so: to get a token moved to the bottom-left vertex of $\graph\psi \circ \graph\phi$, we have no need to put one token in the bottom-right vertex.  Therefore, once we have formed $\graph{\psi\circ\phi}$, we have lost  crucial information about which sources and sinks are \emph{directly} connected with which others, because we have collapsed the entire connected component into a single internal transition, with no internal places. As it happens, by computing the composite graph in a different order, instead, no new paths have been created, hence no cycles appear where there should not be. After all, by Theorem~\ref{theorem:acyclicity implies dinaturality GENERAL} we know that $\chi \circ \psi \circ \phi$ is dinatural because it can be written as the composite of two dinatural transformations, namely $\chi \circ \psi$ and $\phi$, whose composite graph is acyclic. 

This tells us that the crucial reason for which associativity fails in our preliminary definition of the category $\fc \B \C$ is that only keeping track of which connected component each of the arguments of the domain and codomain functors belongs to is not enough: we are forgetting too much information, namely the paths that directly connect the white and grey boxes. Hence our transformations will have to be equipped with more complicated Petri Nets than their standard graph that do contain internal places, and upon composition we shall simply link the graphs together along the common interface, without collapsing entire connected components into a single transition. 

Recall from Definition~\ref{def:FBCF petri net} that a FBCF Petri Net is a net where all the places have at most one input and at most one output transition. We now introduce the category of FBCF Petri Nets, using the usual definition of morphism for bipartite graphs.

\begin{definition}
	The category $\PN$ consists of the following data:
	\begin{itemize}
		\item objects are FBCF Petri Nets $N=(P_N,T_N,\inp{},\out{})$ together with a fixed ordering of its connected components. Such an ordering will allow us to speak about the ``$i$-th connected component'' of $N$;
		\item a morphism $f \colon N \to M$ is a pair of functions $(f_P,f_T)$, for $f_P \colon P_N \to P_M$ and $f_T \colon T_N \to T_M$, such that for all $t \in T_N$
		\[
		\inp{f_T(t)} = \{f_P(p) \mid p \in \inp t \} \quad \text{and} \quad \out{f_T(t)} = \{ f_P(p) \mid p \in \out t  \}.
		\]
	\end{itemize}
\end{definition}

Note that if $f \colon N \to M$ is a morphism in $\PN$ then $f$ preserves (undirected) paths, hence for $C$ a connected component of $N$ we have that $f(C)$ is connected. In particular, if $f$ is an isomorphism then $f(C)$ is a connected component of $M$.

\begin{remark}\label{remark:finite sets are in PN}
	We have a canonical inclusion $\finset \to \PN$ by seeing a set as a Petri Net with only places and no transitions.	
\end{remark}

For a function $x \colon A \to B$ of sets we call $\parts x \colon \parts A \to \parts B$ the action of the covariant powerset functor on $x$, that is the function such that $\parts x (S) = \{x(a) \mid a \in S\} $ for $S \subseteq A$. We then have that if $f \colon N \to M$ is a morphism in $\PN$, then
\[
\begin{tikzcd}
T_N \ar[r,"f_T"] \ar[d,"\inp{}"'] & T_M \ar[d,"\inp{}"] \\
\parts{P_N} \ar[r,"\parts{f_P}"] & \parts{P_M}
\end{tikzcd}
\quad \text{and} \quad
\begin{tikzcd}
T_N \ar[r,"f_T"] \ar[d,"\out{}"'] & T_M \ar[d,"\out{}"] \\
\parts{P_N} \ar[r,"\parts{f_P}"] & \parts{P_M}
\end{tikzcd}
\]
commute by definition of the category $\PN$.

It turns out that $\PN$ admits pushouts, hence we can form a category $\cospan\PN$.

\begin{proposition}\label{prop: pushouts in PN}
	Let $N,M,L$ be in $\PN$, and consider the following diagram in $\PN$:
	\begin{equation}\label{diagram: pushout in PN}
		\begin{tikzcd}[column sep=3em]
			(P_N,T_N,\Inp {N} ,\Out {N} ) \ar[r,"{(g_P,g_T)}"] \ar[d,"{(f_P,f_T)}"'] & (P_L,T_L,\Inp L,\Out L ) \ar[d,"{(k_P,k_T)}"] \\
			(P_M,T_M,\Inp M,\Out M) \ar[r,"{(h_P,h_T)}"] & (P_Q,T_Q,\Inp Q,\Out Q)
		\end{tikzcd}
	\end{equation}
	where
	\[
	\begin{tikzcd}
	P_N \ar[r,"g_P"] \ar[d,"f_P"'] & P_L \ar[d,"k_P"] \\
	P_M \ar[r,"h_P"] & P_Q
	\end{tikzcd}
	\quad {and} \quad
	\begin{tikzcd}
	T_N \ar[r,"g_T"] \ar[d,"f_T"'] & T_L \ar[d,"k_T"] \\
	T_M \ar[r,"h_T"] & T_Q
	\end{tikzcd}
	\]
	are pushouts and $\Inp Q \colon T_Q \to \parts{P_Q}$ is the unique map (the dashed one) that makes the following diagram commute:
	\[
	\begin{tikzcd}
	\parts{P_N} \ar[dddr,bend angle=20,bend right,"\parts{f_P}"'] \ar[rrrd,bend angle=20,bend left,"\parts{g_P}"] \\
	& T_N \ar[r,"g_T"] \ar[d,"f_T"'] \ar[ul,"\Inp N"] & T_L \ar[d,"k_T"] \ar[r,"\Inp L"] & \parts{P_L} \ar[dd,"\parts{k_P}"] \\
	& T_M \ar[r,"h_T"] \ar[d,"\Inp M"] & T_Q \ar[dr,dashed] \\
	& \parts{P_M} \ar[rr,"\parts{h_P}"] & & \parts{P_Q}
	\end{tikzcd}
	\]
	$\Out Q \colon T_Q \to \parts{P_Q}$ is defined analogously. Then (\ref{diagram: pushout in PN}) is a pushout.
\end{proposition}
\begin{proof}
	It is easily checked that (\ref{diagram: pushout in PN}) satisfies the definition of pushout.\qed
\end{proof}

Remember from Remark~\ref{remark:finite sets are in PN} that finite sets can be seen as places-only Petri Nets: if $S$ is a set and $N$ is an object in $\PN$, then a morphism $f \colon S \to N$ in $\PN$ is a pair of functions $f=(f_P,f_T)$ where $f_T$ is the empty map $\emptyset \colon \emptyset \to T_N$. Hence, by little abuse of notation, we will refer to $f_P$ simply as $f$. 

For later convenience, we consider the following subcategory of $\cospan\PN$, whose morphisms are essentially Petri Nets $N$ in $\PN$ with ``interfaces'', that is specific places seen as ``inputs'' and ``outputs'' of $N$. Composition will then be computed by ``gluing together'' two consecutive nets along the common interface.

\begin{definition}\label{definition:graph category}
	The category $\gc$ consists of the following data:
	\begin{itemize}
		\item objects are lists in $\List\{+,-\}$;
		\item morphisms $f \colon \alpha \to \beta$ are (equivalence classes of) cospans in $\PN$ of the form
		\[
		\begin{tikzcd}
		\length\alpha \ar[r,"\lambda"] & N & \ar[l,"\rho"'] \length\beta
		\end{tikzcd}
		\]
		where 
		\begin{itemize}[leftmargin=*]
			\item $\lambda \colon \length\alpha \to P_N$ and $\rho \colon \length\beta \to P_N$ are injective functions, hence we can see $\length\alpha$ and $\length\beta$ as subsets of $P_N$; 
			\item $\mathit{sources}(N) = \{ \lambda(i) \mid \alpha_i=+ \} \cup \{ \rho(i) \mid \beta_i = - \}$;
			\item $\mathit{sinks}(N) = \{ \lambda(i) \mid \alpha_i=- \} \cup \{ \rho(i) \mid \beta_i = + \}$.
		\end{itemize}
		Two such cospans are in the same class if and only if they differ by an isomorphism of Petri Nets on $N$ coherent with $\lambda$, $\rho$ and the ordering of the connected components of $N$; 
		\item composition is that of $\cospan\PN$. 
	\end{itemize}
\end{definition}

\begin{proposition}
	Composition in $\gc$ is well defined.
\end{proposition}
\begin{proof}
	Consider 
	$
	\begin{tikzcd}[cramped,sep=small]
	\length\alpha \ar[r,"\lambda"] & M & \ar[l,"\rho"'] \length\beta
	\end{tikzcd}
	$
	and
	$
	\begin{tikzcd}[cramped,sep=small]
	\length\beta \ar[r,"\lambda'"] & L & \ar[l,"\rho'"'] \length\gamma
	\end{tikzcd}
	$
	two morphisms in $\gc$. By Proposition~\ref{prop: pushouts in PN} then, their composite is given by computing the pushouts
	\[
	\begin{tikzcd}
	\length\beta \ar[r,"\lambda'"] \ar[d,"\rho"'] & P_L \ar[d,"k_P"] \\
	P_M \ar[r,"h_P"] & P_Q
	\end{tikzcd}
	\quad \text{and} \quad
	\begin{tikzcd}
	\emptyset \ar[r,"\emptyset"] \ar[d,"\emptyset"'] & T_L \ar[d,"k_T"] \\
	T_M \ar[r,"h_T"] & T_Q
	\end{tikzcd}
	\]
	Now, the injectivity of $\rho$ and $\lambda'$ implies that $k_P$ and $h_P$ are also injective, as the pushout (in $\Set$) of an injective map against another yields injective functions. $P_Q$, in particular, can be seen as the quotient of $P_M + P_L$ where the elements of $P_M$ and $P_L$ with a common pre-image in $\length\beta$ are identified. Next, the pushout of the empty map against itself yields as a result the coproduct, thus $T_Q = T_M + T_L$ where $h_T$ and $k_T$ are the injections. Hence, the input function of the composite is defined as follows:
	\[
	\begin{tikzcd}[row sep=0em,ampersand replacement=\&]
	T_M + T_L \ar[r,"\inp{}"] \& \parts{P_Q} \\
	t \ar[r,|->] \& \begin{cases}
	\inp{_M(t)} & t \in T_M \\
	\inp{_L(t)} & t \in T_L
	\end{cases}
	\end{tikzcd}
	\]
	and similarly for the output function. All in all, therefore, composition in $\gc$ is computed by ``glueing'' together the Petri Nets $M$ and $L$ along the common $\length\beta$-places; the resulting morphism of $\gc$ is
	\[
	\begin{tikzcd}[column sep=3em]
	\length\alpha \ar[r,"h_P \circ \lambda"] & L \circ M & \length\beta \ar[l,"k_P \circ \rho'"']. 
	\end{tikzcd}
	\]
	Now, for all $i \in \length\beta$, if $\beta_i=+$ then $\rho(i)$ is a sink of $M$ and $\lambda'(i)$ a source of $L$; if $\beta_i=-$ instead then $\rho(i)$ is a source of $M$ and $\lambda'(i)$ a sink of $L$: in every case, once we glue together $M$ and $L$ along the $\length\beta$-places to form the composite net $L \circ M$, these become internal places of $L \circ M$, with at most one input and one output transition each (depending whether they are proper sources or sinks in $M$ and $L$). Hence $L \circ M$ is still a FBCF Petri Net, and 
	\begin{align*}
	\mathit{sources}(L \circ M) &= \bigl(\mathit{sources} (L) \setminus \rho(\length\beta) \bigr) \cup
	\bigl(\mathit{sources} (L)\setminus \lambda'(\length\beta) \bigr) \\
	&= \{h_P \circ \lambda(i) \mid \alpha_i = +\} \cup \{k_P \circ \rho'(i) \mid \gamma_i =-\}
	\end{align*}
	and similarly for $\mathit{sinks}(N' \circ N)$. \qed
\end{proof}

\paragraph{Generalised graphs of a transformation}

We can now start working towards the definition of a category $\fc \B \C$ of functors of mixed variance and transformations that are dinatural only on some of their variables; $\fc \B \C$ will be a category over $\gcf$ in the sense that transformations in $\fc \B \C$ will carry along, as part of their data, certain cospans in $\PN$. The category of graphs $\gcf$ will be built from $\fc \B \C$ by forgetting the transformations. As such, $\gcf$ will be defined \emph{after} $\fc \B \C$.

It is clear how to define the objects of $\fc \B \C$: they will be pairs $(\alpha,F \colon \B^\alpha \to \C)$. Morphisms are less obvious to define, as we learnt in our preliminary attempt on p.~\pageref{first attempt}. A morphism $(\alpha,F) \to (\beta,G)$  will consist of a transformation $\phi \colon F \to G$ of type 
$
\begin{tikzcd}[cramped,sep=small]
\length\alpha \ar[r,"\sigma"] & n & \length\beta \ar[l,"\tau"']
\end{tikzcd}
$, 
together though with a morphism 
$
\begin{tikzcd}[cramped,sep=small]
\length\alpha \ar[r,"\overline\sigma "] & N & \length\beta \ar[l,"\overline\tau "']
\end{tikzcd}
$ in $\gc$ coherent with the type of $\phi$, in the sense that the Petri Net $N$, under certain conditions, looks exactly like $\graph\phi$ as in Definition~\ref{def:standard graph} except that it allows for internal places as well. For example, if $\psi_1$ and $\psi_2$ are two arbitrary consecutive transformations, $\graph{\psi_2} \circ \graph{\psi_1}$ will be coherent with the type of $\psi_2\circ\psi_1$. In other words, $N$ will have $n$ connected components, its sources (sinks) are exactly the places corresponding to the positive (negative) entries of $\alpha$ and the negative (positive) entries of $\beta$, and elements in $\length\alpha$ ($\length\beta$) mapped by $\sigma$ ($\tau$) into the same $i \in \{1,\dots,n\}$ will belong to the $i$-th connected component of $N$. A priori $N$ can contain places with no inputs or outputs: this will be useful for the special case of $\phi = \id F$ as we shall see in Theorem~\ref{theorem: {B,C} is a category}; however, if all sources and sinks in $N$ are proper, then $N$ plays the role of a generalised $\graph{\phi}$.

\begin{definition}\label{definition: generalised graph of transformation}
	Let $\phi \colon F \to G$ be a transformation of type
	$
	\begin{tikzcd}[cramped,sep=small]
	\length\alpha \ar[r,"\sigma"] & n & \length\beta \ar[l,"\tau"']
	\end{tikzcd}
	$.
	A cospan  
	$
	\begin{tikzcd}[cramped,sep=small]
	\length\alpha \ar[r,"\overline\sigma"] & N & \length\beta \ar[l,"\overline\tau"']
	\end{tikzcd}
	$ in $\PN$, which is a representative of a morphism in $\gc$  (hence $\overline\sigma$ and $\overline\tau$ are injective),
	is said to be \emph{coherent with the type of $\phi$} if and only if the following conditions are satisfied:
	\begin{itemize}[leftmargin=*]
		\item $N$ has $n$ connected components; 
		\item for all $i \in \length\alpha$ and $j \in \length\beta$, $\overline\sigma(i)$ belongs to the $\sigma(i)$-th connected component of $N$ and $\overline\tau(j)$ belongs to the $\tau(j)$-th connected component of $N$.
	\end{itemize}
	In this case we say that $N$ is a \emph{generalised graph of $\phi$}.
\end{definition}

\begin{example}\label{example: graph and type are a generalised graph}
	For $\phi \colon F \to G$ a transformation of type
	$
	\begin{tikzcd}[cramped,sep=small]
	\length\alpha \ar[r,"\sigma"] & n & \length\beta \ar[l,"\tau"']
	\end{tikzcd}
	$,
	recall that the set of places of $\graph\phi$ is $P = \length\alpha + \length\beta$. If we call $\injP {\length\alpha}$ and $\injP {\length\beta}$ the injections as in Definition~\ref{def:standard graph}, then
	\[
	\begin{tikzcd}
	\length\alpha \ar[r,"\injP{\length\alpha}"] & \Gamma(\phi) & \ar[l,"\injP{\length\beta}"'] \length\beta
	\end{tikzcd}
	\]
	is indeed coherent with the type of $\phi$.
	Also $
	\begin{tikzcd}[cramped,sep=small]
	\length\alpha \ar[r,"\sigma"] & n & \length\beta \ar[l,"\tau"']
	\end{tikzcd}
	$
	itself, seen as a cospan in $\PN$, is coherent with itself. 
\end{example}

\begin{remark}
	If $N$ is a generalised graph of $\phi$ as in the notations of Definition~\ref{definition: generalised graph of transformation} and does not have any place which is a source and a sink at once, then $N$ has exactly $\length\alpha + \length\beta$ sources and sinks and their union coincides with the joint image of $\overline\sigma$ and $\overline\tau$.  Moreover, $\overline\sigma$ and $\overline\tau$ have to make sure that they map elements of their domain into places belonging to the correct connected component: in this way, $N$ reflects the type of $\phi$ in a Petri Net like $\graph\phi$, with the possible addition of internal places.
\end{remark}

We shall now show how composition in $\gc$ preserves generalised graphs, in the following sense.

\begin{proposition}\label{proposition: composition in G preservers generalised graphs}
	Let $\phi \colon F \to G$ and $\psi \colon G \to H$ be transformations of type, respectively, 
	$
	\begin{tikzcd}[cramped,sep=small]
	\length\alpha \ar[r,"\sigma"] & n & \length\beta \ar[l,"\tau"']
	\end{tikzcd}
	$
	and
	$
	\begin{tikzcd}[cramped,sep=small]
	\length\beta \ar[r,"\eta"] & m & \length\gamma \ar[l,"\theta"']
	\end{tikzcd}
	$; let also 
	$u=
	\begin{tikzcd}[cramped,sep=small]
	\length\alpha \ar[r,"\overline\sigma"] & N & \length\beta \ar[l,"\overline\tau"']
	\end{tikzcd}
	$
	and
	$v=
	\begin{tikzcd}[cramped,sep=small]
	\length\beta \ar[r,"\overline\eta"] & N' & \length\gamma \ar[l,"\overline\theta"']
	\end{tikzcd}
	$
	be cospans in $\PN$ coherent with the type of $\phi$ and $\psi$, respectively. Suppose the type of $\psi \circ \phi$ is given by
	\[
	\begin{tikzcd}
	& & \length\gamma \ar[d,"\theta"] \\
	& \length\beta \ar[d,"\tau"'] \ar[r,"\eta"] \ar[dr, phantom, "\ulcorner" very near start] & m \ar[d,"\xi"] \\
	\length\alpha \ar[r,"\sigma"] & n \ar[r,"\zeta"] & l
	\end{tikzcd}
	\]
	and that the composite in $\gc$ of $u$ and $v$ is given by
	\begin{equation}\label{composite generalised graphs}
		\begin{tikzcd}
			& & \length\gamma \ar[d,"\overline\theta"] \\
			& \length\beta \ar[d,"\overline\tau"'] \ar[r,"\overline\eta"] \ar[dr, phantom, "\ulcorner" very near start] & N' \ar[d,"\overline\xi"] \\
			\length\alpha \ar[r,"\overline\sigma"] & N \ar[r,"\overline\zeta"] & N' \circ N
		\end{tikzcd}
	\end{equation}
	Then $v\circ u$ is coherent with the type of $\psi \circ \phi$.
\end{proposition}
\begin{proof}
	As we said in the discussion after Definition~\ref{definition:graph category}, $N' \circ N$ is obtained by gluing together $N$ and $N'$ along the $\length\beta$ places which they have in common. The number of connected components of $N' \circ N$ is indeed $l$ by construction. The morphisms $\overline\zeta$ and $\overline\xi$ in $\PN$ are pairs of injections that map each place and transition of $N$ and $N'$ to itself in the composite $N' \circ N$. This means that $\overline\zeta \overline\sigma(i)$ does belong to the $\zeta\sigma(i)$-th connected component of $N' \circ N$, as the latter contains the $\sigma(i)$-th c.c.\ of $N$; similarly the $\overline\xi \overline\theta (j)$ belongs to the $\xi\theta(j)$-th c.c.\ of $N' \circ N$.  \qed
\end{proof}

The morphisms of our generalised functor category $\fc \B \C$ will be, therefore, transformations $\phi$ equipped with a generalised graph $N$ and a discriminant function that tells us in which variables $\phi$ is dinatural. The Petri Net $N$ will not be arbitrary though: unless $\phi$ is an identity transformation, $N$ can be either $\graph\phi$ or $\graph{\phi_k} \circ \dots \circ \graph{\phi_1}$, for some consecutive transformations $\phi_1,\dots,\phi_k$ such that $\phi = \phi_k \circ \dots \phi_1$. Therefore, only transformations which are \emph{explicitly} recognisable as the composite of two or more families of morphisms are allowed to have an associated Petri Net, containing internal places, that is not their standard graph.

\begin{definition}\label{def: generalised functor category}
	Let $\B$ and $\C$ be categories. The \emph{generalised functor category} $\fc \B \C$ consists of the following data:
	\begin{itemize}[leftmargin=*]
		\item objects are pairs $(\alpha,F)$, for $\alpha \in \List\{+,-\}$ and $F \colon \B^\alpha \to \C$ a functor;
		\item morphisms $(\alpha,F) \to (\beta,G)$ are equivalence classes of tuples
		\[
		\Phi = (\phi,
		\begin{tikzcd}[cramped,sep=small]
		\length\alpha \ar[r,"\sigma"] & n & \length\beta \ar[l,"\tau"']
		\end{tikzcd},
		\begin{tikzcd}[cramped,sep=small]
		\length\alpha \ar[r,"\overline\sigma"] & N & \length\beta \ar[l,"\overline\tau"']
		\end{tikzcd},
		\Delta_\Phi
		)
		\]
		where:
		\begin{itemize}[leftmargin=*]
			\item $\phi \colon F \to G$ is a transformation of type
			$
			\begin{tikzcd}[cramped,sep=small]
			\length\alpha \ar[r,"\sigma"] & n & \length\beta \ar[l,"\tau"']
			\end{tikzcd}
			$,
			\item $
			\begin{tikzcd}[cramped,sep=small]
			\length\alpha \ar[r,"\overline\sigma"] & N & \length\beta \ar[l,"\overline\tau"']
			\end{tikzcd}
			$ is a representative of a morphism in $\gc$ coherent with the type of $\phi$,
			\item $\Delta_\Phi \colon n \to \{0,1\}$ is a function such that $\Delta_\Phi (i) = 1$ implies that the $i$-th connected component of $N$ is acyclic and $\phi$ is dinatural in its $i$-th variable.
			
		\end{itemize}
		Moreover: 
		\begin{itemize}[leftmargin=*]
			\item If $N$ consists of $n$ places and no transitions, then $(\alpha,F) = (\beta,G)$, $\phi=\id F$, $\sigma=\tau=\overline\sigma=\overline\tau=\id{\length\alpha}$ and $\Delta_\Phi = K_1$, the constant function equal to $1$; in this case $\Phi$ is the identity morphism of the object $(\alpha,F)$.
			\item If $N = \graph\phi$, $\overline\sigma = \injP{\length\alpha}$ and $\overline\tau=\injP{\length\beta}$, we say that $\Phi$ is \emph{atomic}.
			\item If $N \ne \graph\phi$ and $\Phi \ne \id{(\alpha,F)}$, then there exist $\Phi_1, \dots, \Phi_k$ atomic such that $\Phi = \Phi_k \circ \dots \circ \Phi_1$ in $\fc \B \C$, according to the composition law to follow in this Definition.
		\end{itemize}
		
		We say that $\Phi \eq \Phi'$, for $\Phi' = (\phi',
		\begin{tikzcd}[cramped,sep=small]
		\length\alpha \ar[r,"\sigma'"] & n & \length\beta \ar[l,"\tau'"']
		\end{tikzcd},
		\begin{tikzcd}[cramped,sep=small]
		\length\alpha \ar[r,"\overline{\sigma'}"] & N' & \length\beta \ar[l,"\overline{\tau'}"']
		\end{tikzcd},
		\Delta_{\Phi'}
		)$,
		if and only if the transformations differ only by a permutation of their variables (in a coherent way with the rest of the data) and $N$ and $N'$ are coherently isomorphic: more precisely, when
		\begin{itemize}[leftmargin=*]
			\item there is a permutation $\pi \colon n \to n$ such that $\sigma'=\pi\sigma$, $\tau'=\pi\tau$, $\phi_{A_1,\dots,A_n}'=\phi_{A_{\pi 1},\dots,A_{\pi n}}$, $\Delta_{\Phi}=\Delta_{\Phi'} \pi$;
			\item there is an isomorphism $f=(f_P,f_T) \colon N \to N'$ in $\PN$ such that the following diagram commutes:
			\[
			\begin{tikzcd}
			\length\alpha \ar[r,"\overline\sigma"] \ar[dr,"\overline{\sigma'}"'] & N \ar[d,"f"] & \length\beta \ar[l,"\overline\tau"'] \ar[dl,"\overline{\tau'}"] \\
			& N' 
			\end{tikzcd}
			\]
			mapping the $i$-th connected component of $N$ to the $\pi(i)$-th connected component of $N'$. 
		\end{itemize}
		\item Composition of $\Phi$ as above and
		\[
		\Psi = 
		(\psi,
		\begin{tikzcd}[cramped,sep=small]
		\length\beta \ar[r,"\eta"] & m & \ar[l,"\theta"'] \length\gamma
		\end{tikzcd},
		\begin{tikzcd}[cramped,sep=small]
		\length\beta \ar[r,"\overline\eta"] & N' & \ar[l,"\overline\theta"'] \length\gamma
		\end{tikzcd},
		\Delta_\Psi
		) \colon (\beta,G) \to (\gamma, H)
		\]
		is component-wise: it is the equivalence class of the tuple 
		\begin{equation}\label{eqn:composition in {B,C}}
			\Psi \circ \Phi = (
			\psi\circ\phi,
			\begin{tikzcd}[cramped,sep=small]
				\length\alpha \ar[r,"\zeta\sigma"] & l & \ar[l,"\xi\theta"'] \length\gamma
			\end{tikzcd},
			\begin{tikzcd}[cramped,sep=small]
				\length\alpha \ar[r,"\overline\zeta \overline\sigma"] & {N'} \circ N & \ar[l,"\overline\xi \overline\theta"'] \length\gamma
			\end{tikzcd},
			\Delta_{\Psi\circ\Phi}
			)
		\end{equation} 
		where $\psi\circ\phi$ is the transformation of type given by the result of the pushout:
		\[
		\begin{tikzcd}
		& & \length\gamma \ar[d,"\theta"] \\
		& \length\beta \ar[d,"\tau"'] \ar[r,"\eta"] \ar[dr, phantom, "\ulcorner" very near start] & m \ar[d,"\xi"] \\
		\length\alpha \ar[r,"\sigma"] & n \ar[r,"\zeta"] & l
		\end{tikzcd}
		\]
		$N' \circ N$ is computed by composing in $\gc$, that is by performing the pushout in $\PN$:
		\[
		\begin{tikzcd}
		& & \length\gamma \ar[d,"\overline\theta"] \\
		& \length\beta \ar[d,"\overline\tau"'] \ar[r,"\overline\eta"] \ar[dr, phantom, "\ulcorner" very near start] & N' \ar[d,"\overline\xi"] \\
		\length\alpha \ar[r,"\overline\sigma"] & N \ar[r,"\overline\zeta"] & N' \circ N
		\end{tikzcd}
		\]
		and the discriminant $\Delta_{\Psi\circ\Phi} \colon l \to \{0,1\}$ is obtained by setting $\Delta_{\Psi\circ\Phi} (x) = 1$ if and only if the $x$-th connected component of $N'\circ N$ is acyclic \emph{and} for all $y \in \zeta^{-1}\{x\}$ and $z \in \xi^{-1}\{x\}$ we have that $\Delta_\Phi(y) = 1 = \Delta_\Psi(z)$. The latter condition is tantamount to asking that $\phi$ and $\psi$ are dinatural in all the variables involved by the $x$-th connected component of the composite graph ${N'}\circ N$ of $\psi\circ\phi$.
	\end{itemize}
\end{definition}

\begin{theorem}\label{theorem: {B,C} is a category}
	$\fc \B \C$ is indeed a category.
\end{theorem}
\begin{proof}
	First of all, if $\Phi$ and $\Psi$ as above are in $\fc \B \C$, it is not difficult to check that the equivalence class of $\Psi \circ \Phi$ as in~(\ref{eqn:composition in {B,C}}) does not depend on the choice of representatives for the classes of $\Phi$ and $\Psi$.
	
	Next, we aim to prove that $\Psi \circ \Phi$ is again a morphism of $\fc \B \C$. By Proposition~\ref{proposition: composition in G preservers generalised graphs} we have that ${N'} \circ N$ is a generalised graph for $\psi\circ\phi$. In order to prove that $\Delta_{\Psi\circ\Phi}$ correctly defines a morphism of $\fc \B \C$, that is that if $\Delta_{\Psi\circ\Phi}(i)=1$ then $\psi\circ\phi$ is indeed dinatural in its $i$-th variable, we first show that composition in $\fc \B \C$ is associative: once we have done that we will use Theorem~\ref{theorem:compositionality with complicated graphs} to conclude.
	
	Consider 
	\begin{align*}
	\Phi_1 &= (\phi_1,
	\begin{tikzcd}[cramped,sep=small,ampersand replacement=\&]
	\length\alpha \ar[r,"\sigma_1"] \& n \& \length\beta \ar[l,"\tau_1"']
	\end{tikzcd},
	\begin{tikzcd}[cramped,sep=small,ampersand replacement=\&]
	\length\alpha \ar[r,"\overline{\sigma_1}"] \& N_1 \& \length\beta \ar[l,"\overline{\tau_1}"']
	\end{tikzcd},
	\Delta_\Phi
	) \colon (\alpha, F) \to (\beta, G), \\
	\Phi_2 &= (
	\phi_2,
	\begin{tikzcd}[cramped,sep=small,ampersand replacement=\&]
	\length\beta \ar[r,"\sigma_2"] \& m \& \length\gamma \ar[l,"\tau_2"']
	\end{tikzcd},
	\begin{tikzcd}[cramped,sep=small,ampersand replacement=\&]
	\length\beta \ar[r,"\overline{\sigma_2}"] \& N_2 \& \length\gamma \ar[l,"\overline{\tau_2}"']
	\end{tikzcd},
	\Delta_{\Phi_2}
	) \colon (\beta, G) \to (\gamma, H), \\
	\Phi_3 &= (
	\phi_3,
	\begin{tikzcd}[cramped,sep=small,ampersand replacement=\&]
	\length\gamma \ar[r,"\sigma_3"] \& p \& \length\delta \ar[l,"\tau_3"']
	\end{tikzcd},
	\begin{tikzcd}[cramped,sep=small,ampersand replacement=\&]
	\length\gamma \ar[r,"\overline{\sigma_3}"] \& N_3 \& \length\delta \ar[l,"\overline{\tau_3}"']
	\end{tikzcd},
	\Delta_{\Phi_3}
	) \colon (\gamma,H) \to (\delta,K).
	\end{align*}
	We know that composition of cospans via pushout is associative, as well as composition of transformations; suppose therefore that $\phi_3 \circ \phi_2 \circ \phi_1$ has type given by:
	\[
	\begin{tikzcd}
	& & & \length\delta \ar[d,"\tau_3"] \\
	& & \length\gamma \ar[r,"\sigma_3"] \ar[d,"\tau_2"'] \ar[dr, phantom, "\ulcorner" very near start] & p \ar[d,"\xi_2"] \\
	& \length\beta \ar[r,"\sigma_2"] \ar[d,"\tau_1"'] \ar[dr, phantom, "\ulcorner" very near start] & m \ar[r,"\zeta_2"] \ar[d,"\xi_1"'] \ar[dr, phantom, "\ulcorner" very near start] & q \ar[d,"\xi_3"] \\
	\length\alpha \ar[r,"\sigma_1"] & n \ar[r,"\zeta_1"] & l \ar[r,"\zeta_3"] & r
	\end{tikzcd}
	\]
	and the generalised graph $N_3 \circ N_2 \circ N_1$ is obtained as the result of the following pushout-pasting:
	\[
	\begin{tikzcd}
	& & & \length\delta \ar[d,"\overline{\tau_3}"] \\
	& & \length\gamma \ar[r,"\overline{\sigma_3}"] \ar[d,"\overline{\tau_2}"'] \ar[dr, phantom, "\ulcorner" very near start] & N_3 \ar[d,"\overline{\xi_2}"] \\
	& \length\beta \ar[r,"\overline{\sigma_2}"] \ar[d,"\overline{\tau_1}"'] \ar[dr, phantom, "\ulcorner" very near start] & N_2 \ar[r,"\overline{\zeta_2}"] \ar[d,"\overline{\xi_1}"'] \ar[dr, phantom, "\ulcorner" very near start] & N_3 \circ N_2 \ar[d,"\overline{\xi_3}"] \\
	\length\alpha \ar[r,"\overline{\sigma_1}"] & N_1 \ar[r,"\overline{\zeta_1}"] & N_2 \circ N_1 \ar[r,"\overline{\zeta_3}"] & N_3 \circ N_2 \circ N_1
	\end{tikzcd}
	\]
	We prove that $\Delta_{\Phi_3 \circ (\Phi_2 \circ \Phi_1)} = \Delta_{(\Phi_3 \circ \Phi_2) \circ \Phi_1}$. We have that $\Delta_{\Phi_3 \circ (\Phi_2 \circ \Phi_1)}(x) = 1$ if and only if, by definition:
	\begin{enumerate}[labelindent=0pt]
		\item[(1)] the $x$-th c.c.\ of $N_3 \circ N_2 \circ N_1$ is acyclic;
		\item[(2)] $\forall y \in \zeta_3^{-1}\{x\} \ldotp \Delta_{\Phi_2 \circ \Phi_1}(y) = 1$;
		\item[(3)] $\forall z \in (\xi_3 \circ \xi_2)^{-1}\{x\} \ldotp \Delta_{\Phi_3}(z) = 1$;
	\end{enumerate}
	which is equivalent to say that:
	\begin{enumerate}[labelindent=0pt]
		\item[(1)] the $x$-th c.c.\ of $N_3 \circ N_2 \circ N_1$ is acyclic;
		\item[(2a)] $\forall y \in l \ldotp \Bigl[ \zeta_3(y) = x \implies \text{$y$-th c.c.\ of $N_2 \circ N_1$ is acyclic} \Bigr] $;
		\item[(2b)] $\forall y \in l \ldotp \Bigl[ \zeta_3(y) = x \implies \forall a \in n \ldotp \Bigl( \zeta_1(a)=y \implies \Delta_{\Phi_1}(a)=1 \Bigr) \Bigr] $;
		\item[(2c)] $\forall y \in l \ldotp \Bigl[ \zeta_3(y) = x \implies \forall b \in m \ldotp \Bigl( \xi_1(b)=y \implies \Delta_{\Phi_2}(b)=1 \Bigr) \Bigr] $;
		\item[(3)] $\forall z \in p \ldotp \Bigl[ \xi_3\bigl(\xi_2(z)\bigr) = x \implies \Delta_{\Phi_3} (z) = 1 \Bigr] $.
	\end{enumerate}
	Call $A$ the conjunction of the conditions above. Next, we have that $\Delta_{(\Phi_3 \circ \Phi_2) \circ \Phi_1}(x)=1$ if and only if:
	\begin{enumerate}[labelindent=0pt]
		\item[(i)] the $x$-th c.c.\ of $N_3 \circ N_2 \circ N_1$ is acyclic;
		\item[(ii)] $\forall a \in n \ldotp \Bigl[ \zeta_3 \bigl( \zeta_1(a) \bigr) =x \implies \Delta_{\Phi_1}(a)=1\Bigr]$;
		\item[(iiia)] $\forall w \in q \ldotp \Bigl[ \xi_3(w)=x \implies \text{ $w$-th c.c.\ of $N_3 \circ N_2$ is acyclic } \Bigr]$;
		\item[(iiib)] $\forall w \in q \ldotp \Bigl[ \xi_3(w)=x \implies \forall b \in m \ldotp \Bigl( \zeta_2(b)=w \implies \Delta_{\Phi_2}(b)=1 \Bigr) \Bigr]$;
		\item[(iiic)] $\forall w \in q \ldotp \Bigl[ \xi_3(w)=x \implies \forall z \in p \ldotp \Bigl( \xi_2(z)=w \implies \Delta_{\Phi_3}(z)=1 \Bigr) \Bigr]$
	\end{enumerate}
	Call $B$ the conjunction of these last five conditions. We prove that $A$ implies $B$; in a similar way one can prove the converse as well.
	\begin{enumerate}[labelindent=0pt]
		\item[(ii)] Let $a \in n$, suppose $\zeta_3\bigl(\zeta_1(a)\bigr) = x$. By (2b), with $y = \zeta_1(a)$, we have $\Delta_{\Phi_1}(a) = 1$.
		\item[(iiia)] Let $w \in q$, suppose $\xi_3(w)=x$. Then the $w$-th c.c.\ of $N_3 \circ N_2$ must be acyclic as it is part of the $x$-th c.c.\ of $N_3 \circ N_2 \circ N_1$, which is acyclic.
		\item[(iiib)] Let $w \in q$, suppose $\xi_3(w)=x$. Let also $b \in m$ and suppose $\zeta_2(b) = w$. Then $x = \xi_3\bigl( \zeta_2(b)\bigr) = \zeta_3 \bigl(\xi_1(b)\bigr)$. By (2c), with $y = \xi_1(b)$, we have $\Delta_{\Phi_2}(b)=1$.
		\item[(iiic)] Let $w \in q$, suppose $\xi_3(w)=x$. Let $z \in p$ be such that $\xi_2(z)=w$. Then $\xi_3\bigl(\xi_2(z)\bigr)=x$: by (3), we have $\Delta_{\Phi_3}(z)=1$.
	\end{enumerate}
	
	Hence composition is associative. Take now $\Phi$ and $\Psi$ consecutive morphisms of $\fc \B \C$ as in the Definition of $\fc \B \C$. Then $\Phi=\Phi_k \circ \dots \circ \Phi_1$ for some $\Phi_j$'s, in particular $\phi=\phi_k \circ \dots \circ \phi_1$ for some $\phi_j$'s, and $\Delta_\Phi (i) =1 $ precisely when the $i$-th connected component of $N$ is acyclic and for all $j \in \{1,\dots,k\}$ the transformation $\phi_j$ is dinatural in all its variables involved in the $i$-th c.c.\ of $N$: one can see this by simply unfolding the definition of $\Delta_{\Phi_k \circ \dots \circ \Phi_1}$, extending the case of $\Delta_{\Phi_3 \circ \Phi_2 \circ \Phi_1}$ above. Similarly for $\Psi=\Psi_{k'} \circ \dots \Psi_1$, with $\psi = \psi_{k'} \circ \dots \psi_1$. We have then that if 
	\[
	N' \circ N = \graph{\psi_{k'}} \circ \dots \circ \graph{\psi_1} \circ \graph{\phi_k} \circ \dots \circ \graph{\phi_1}
	\]
	is acyclic in its $x$-th connected component and for all $y \in \zeta^{-1}\{x\}$ and $z \in \xi^{-1}\{x\}$ we have that $\Delta_\Phi(y) = 1 = \Delta_\Psi(z)$, then all the $\phi_j$'s and $\psi_j$'s are dinatural in all their variables involved in the $x$-th connected component of $N' \circ N$: by Theorem~\ref{theorem:compositionality with complicated graphs}, we have that $\psi\circ\phi$ is dinatural in its $x$-th variable. Hence $\Psi\circ\Phi$ is still a morphism of $\fc \B \C$.
	
	All that is left to prove is that composition is unitary where the identity morphism of $(\alpha,F)$ is given by the equivalence class of
	\[
	(
	\id F, 
	\begin{tikzcd}[cramped,sep=small]
	\length\alpha \ar[r,"\id{}"] & \length\alpha & \ar[l,"\id{}"'] \length\alpha
	\end{tikzcd},
	\begin{tikzcd}[cramped,sep=small]
	\length\alpha \ar[r,"\id{}"] & \length\alpha & \ar[l,"\id{}"'] \length\alpha
	\end{tikzcd},
	K_1
	),
	\]
	which is indeed a morphism of $\fc\B\C$ because, 
	as discussed in Example~\ref{example: graph and type are a generalised graph} we have that $\length\alpha$ is a generalised graph for $\id F$; moreover, the identity transformation is indeed (di)natural in all its variables, therefore the constant function equal to $1$, $K_1$, is a valid discriminant function for $\id{\length\alpha}$.
	
	Let 
	\[
	\Phi = (\phi,
	\begin{tikzcd}[cramped,sep=small]
	\length\alpha \ar[r,"\sigma"] & n & \length\beta \ar[l,"\tau"']
	\end{tikzcd},
	\begin{tikzcd}[cramped,sep=small]
	\length\alpha \ar[r,"\overline\sigma"] & N & \length\beta \ar[l,"\overline\tau"']
	\end{tikzcd},
	\Delta_\Phi
	) \colon (\alpha,F) \to (\beta,G). 
	\]
	We prove that $\Phi \circ \id{(\alpha,F)} = \Phi$ and ${\id{(\beta,G)}} \circ \Phi = \Phi$ (by ``$\Phi$'' here we mean its equivalence class). It is clear that $\Phi \circ {\id{(\alpha,F)}}$ consists of $\phi$ together with its type and generalised graph as specified in $\Phi$. Also, $\Delta_{\Phi \circ \id{(\alpha,F)}}(x) = 1$ precisely when the $x$-th connected component of $N$ is acyclic and $\Delta_{\Phi}(x)=1$, by definition. Given that $\Delta_\Phi(x)=1$ implies that the $x$-th c.c.\ of $N$ is acyclic, we have that  $\Delta_{\Phi \circ \id{(\alpha,F)}} = \Delta_\Phi$. One can prove in a similar way the other identity law. \qed
\end{proof}

\begin{remark}
	The condition ``$\Delta_\Phi(i)=1$ implies that the $i$-th connected component of $N$ is acyclic'' in Definition~\ref{def: generalised functor category} is designed to ignore dinaturality properties that happen to be satisfied ``by accident'', as it were, which could cause problems upon composition. Indeed, suppose that we have a transformation $\phi$ which is the composite of four transformations $\phi_1,\dots,\phi_4$, whose resulting generalised graph, obtained by pasting together $\Gamma(\phi_1),\dots,\Gamma(\phi_4)$, is as follows:
	\[
	N= \quad	\begin{tikzpicture}
	\matrix[column sep=2.4mm,row sep=0.4cm]{
		\node (1) [category] {}; \\
		\node (2) [component] {}; & & & \node (7) [component] {}; \\
		\node (3) [category] {}; & & \node (8) [category] {}; & & \node (6) [opCategory] {}; \\
		& \node (4) [component] {}; & & & \node (5) [component] {}; \\
		&	\node (A) [category] {}; & & & \node(F) [opCategory] {};\\
		&	\node (B) [component] {}; & & & \node(J) [component] {};\\
		\node (C) [category] {}; & & \node(D) [category] {}; & & \node(E) [opCategory] {};\\
		\node (H) [component] {}; & & & \node(I) [component] {};\\
		\node (G) [category] {}; & & & \\
	};
	\graph[use existing nodes]{
		1 -> 2 -> 3 -> 4 -> A -> B -> {C, D};
		C -> H -> G;
		D -> I -> E -> J -> F -> 5 -> 6 -> 7 -> 8 -> 4;
	};
	\node[coordinate](p) at (-2,0) {};
	\node[coordinate](q) at (2,0) {};
	\draw [dashed] (3.west -| p) -- (6.east -| q);
	\draw [dashed] (A.west -| p) -- (F.east -| q);
	\draw [dashed] (C.west -| p) -- (E.east -| q);
	\end{tikzpicture}
	\]
	Call $\Phi$ the tuple in $\fc \B \C$ consisting of $\phi$ with its type 
	$
	\begin{tikzcd}[cramped,sep=small]
	1 \ar[r] & 1 & 1 \ar[l]
	\end{tikzcd}
	$ and $N$ as a generalised graph, as a composite of the atomic morphisms of $\fc \B \C$ given by $\phi_1,\dots,\phi_4$. Suppose that $\phi$ happens to be dinatural in its only variable for some reason (extreme example: the category $\C$ is the terminal category). If in the definition of $\fc \B \C$ the only condition on $\Delta$ were ``$\Delta_\Phi(i) = 1$ implies $\phi$ dinatural in its $i$-th variable'', without requiring that the $i$-th connected component of $N$ be acyclic if $\Delta_\Phi(i)=1$, then equipping $\phi$ in $\Phi$ with a discriminant function $\Delta_\Phi$ defined as
	\[
	\begin{tikzcd}[row sep=0pt]
	1 \ar[r,"\Delta_\Phi"] & 1 \\
	1 \ar[r,|->] & 1
	\end{tikzcd}
	\] would be legitimate. Compose now $\Phi$ with the identity morphism of $\fc \B \C$: by definition we would obtain again $\Phi$ except for the discriminant function, which would be defined as $\Delta_{\Phi \circ \id{}}(1)=0$ because the composite graph, which is $N$, is not acyclic. Composition would not be unitary! The condition ``the $i$-th connected component of $N$ is acyclic whenever $\Delta_\Phi(i)=1$'' in Definition~\ref{def: generalised functor category} is therefore not only sufficient, but also necessary for unitarity of composition in $\fc \B \C$.
	
\end{remark}

\begin{remark}\label{remark:non-atomic morphisms of {B,C}}
	Although it is impossible, in general, to judge whether a transformation is or is not a composite of others by looking at its type, one can distinguish atomic morphisms of $\FC \B \C$ from composite morphisms by looking at the generalised graph $N$ they come  with. Indeed, if 
	\[
	\Phi = (\phi,
	\begin{tikzcd}[cramped,sep=small]
	\length\alpha \ar[r,"\sigma"] & n & \length\beta \ar[l,"\tau"']
	\end{tikzcd},
	\begin{tikzcd}[cramped,sep=small]
	\length\alpha \ar[r,"\overline\sigma"] & N & \length\beta \ar[l,"\overline\tau"']
	\end{tikzcd},
	\Delta_\Phi
	)
	\]
	is a non-identity morphism of $\FC \B \C$, then $\Phi$ is atomic if and only if $N=\graph\phi$. In case $N \ne \graph\phi$, then $N$ contains internal places as a result of composing together ``atomic'' graphs of transformations: that is, we have that $\phi = \phi_k \circ \dots \circ \phi_1$ for some transformations $\phi_i$, and $N=\graph{\phi_k} \circ \dots \circ \graph{\phi_1}$. This decomposition of $\phi$ and $N$ is not necessarily unique.
\end{remark}

\paragraph{The category of graphs}

We can now finally individuate the category $\gcf$ of graphs of transformations. To do so, we will first build a category $\GC$, which will consist of those morphisms in $\gc$ that are the generalised graph of a transformation in $\fc \B \C$, together with a discriminant function. The category of graphs $\gcf$ we seek will be defined as a subcategory of it. 

We begin by defining the notion of \emph{skeleton} of a morphism in $\gc$, as it will be useful later on.

\begin{definition}
	Let 
	$
	f = \begin{tikzcd}[cramped,sep=small]
	\length\alpha \ar[r,"\overline\sigma"] & N & \ar[l,"\overline\tau"'] \length\beta
	\end{tikzcd}
	$
	be a morphism in $\gc$, and let $n$ be the number of connected components of $N$. The \emph{skeleton} of the cospan $f$ is an (equivalence class of) cospan(s) in $\finset$
	\[
	\begin{tikzcd}
	\length\alpha \ar[r,"\sigma"] & n & \length\beta \ar[l,"\tau"']
	\end{tikzcd}
	\]
	where $\sigma(i)$ is the number of the connected component of $N$ to which $\overline\sigma(i)$ belongs to, and similarly is defined $\tau$.
\end{definition}

\begin{remark}
	If $\phi$ is a transformation and $N$ is a generalised graph of $\phi$, then the type of $\phi$ is the skeleton of $N$.
\end{remark}

The category $\GC$ will then consist of only part of the data of $\fc \B \C$, obtained, as it were, by discarding functors and transformations, and only considering the graphs and the discriminant functions.

\begin{definition}\label{definition:graph category definitive}
	The category $\GC$ of graphs consists of the following data.
	\begin{itemize}[leftmargin=*]
		\item Objects are lists in $\List\{+,-\}$.
		\item Morphisms $\alpha \to \beta$ are equivalence classes of pairs
		\[
		\bigl(
		\begin{tikzcd}
		\length\alpha \ar[r,"\overline\sigma"] & N & \ar[l,"\overline\tau"'] \length\beta
		\end{tikzcd},
		\Delta_N  \bigr)
		\]
		where:
		\begin{itemize}
			\item $(\overline\sigma,\overline\tau,N)$ is a morphism in $\gc$,
			\item let $n$ be the number of connected components of $N$: then $\Delta_N \colon n \to \{0,1\}$ is called \emph{discriminant function} and it is such that  $\Delta(i)=1$ implies that the $i$-th connected component of $N$ is acyclic.		
		\end{itemize}
		A pair above is equivalent to another $((\overline\sigma',\overline\tau',N'),\Delta_{N'})$, where $N'$ also has $n$ connected components, if and only if there exists $f \colon N \to N'$ an isomorphism in $\PN$ and $\pi \colon n \to n$ a permutation such that
		\[
		\begin{tikzcd}
		\length\alpha \ar[r,"\overline\sigma"] \ar[dr,"\overline{\sigma'}"'] & N \ar[d,"f"] & \length\beta \ar[l,"\overline\tau"'] \ar[dl,"\overline{\tau'}"] \\
		& N' 
		\end{tikzcd}
		\quad \text{and} \quad
		\begin{tikzcd}
		n \ar[r,"\Delta_N"] \ar[d,"\pi"'] & \{0,1\} \\
		n \ar[ur,"\Delta_N'"'] 
		\end{tikzcd}
		\]
		commute and $f$ maps the $i$-th c.c.\ of $N$ to the $\pi(i)$-th c.c.\ of $N'$.
		\item
		Composition is defined exactly as in $\fc \B \C$. To wit, composition of
		\[
		\bigl(
		\begin{tikzcd}
		\length\alpha \ar[r,"\overline\sigma"] & N & \ar[l,"\overline\tau"'] \length\beta
		\end{tikzcd},
		\Delta_N  \bigr)
		\quad \text{and} \quad
		\bigl(
		\begin{tikzcd}
		\length\beta \ar[r,"\overline\eta"] & N & \ar[l,"\overline\theta"'] \length\gamma
		\end{tikzcd},
		\Delta_{N'}  \bigr)
		\]
		is the equivalence class of the pair
		\[
		(
		\begin{tikzcd}
		\length\alpha \ar[r,"\overline\zeta \overline\sigma"] & {N'} \circ N & \ar[l,"\overline\xi \overline\theta"'] \length\gamma
		\end{tikzcd},
		\Delta_{g \circ f}
		)
		\]
		where $N' \circ N$ is the Petri Net given by the result of the pushout
		\[
		\begin{tikzcd}
		& & \length\gamma \ar[d,"\overline\theta"] \\
		& \length\beta \ar[d,"\overline\tau"'] \ar[r,"\overline\eta"] \ar[dr, phantom, "\ulcorner" very near start] & N' \ar[d,"\overline\xi"] \\
		\length\alpha \ar[r,"\overline\sigma"] & N \ar[r,"\overline\zeta"] & N' \circ N
		\end{tikzcd}
		\] 
		and $\Delta_{N' \circ N}$ is defined as follows. If
		$
		\begin{tikzcd}[cramped,sep=small]
		\length\alpha \ar[r,"\sigma"] & n & \length\beta \ar[l,"\tau"']
		\end{tikzcd}
		$
		and 
		$
		\begin{tikzcd}[cramped,sep=small]
		\length\beta \ar[r,"\eta"] & m & \length\gamma \ar[l,"\theta"']
		\end{tikzcd}
		$
		are the skeletons of $(\overline\sigma,\overline\tau,N)$ and $(\overline\eta,\overline\theta,N')$ respectively, then the skeleton of $(\overline\zeta\overline\sigma,\overline\xi\overline\theta,N'\circ N)$ is given by the pushout 
		\[
		\begin{tikzcd}
		& & \length\gamma \ar[d,"\theta"] \\
		& \length\beta \ar[d,"\tau"'] \ar[r,"\eta"] \ar[dr, phantom, "\ulcorner" very near start] & m \ar[d,"\xi"] \\
		\length\alpha \ar[r,"\sigma"] & n \ar[r,"\zeta"] & l
		\end{tikzcd}
		\]
		(cf.\ Proposition~\ref{proposition: composition in G preservers generalised graphs}). Define therefore $\Delta_{N' \circ N}(x)=1$ if and only if the $x$-th connected component of $N'\circ N$ is acyclic \emph{and} for all $y \in \zeta^{-1}\{x\}$ and $z \in \xi^{-1}\{x\}$ we have that $\Delta_N(y) = 1 = \Delta_{N'}(z)$.
	\end{itemize}
\end{definition}

\begin{definition}
	The category $\gcf$ of graphs is the wide subcategory  of $\GC$ (that is, it contains all the objects of $\GC$) generated by equivalence classes of pairs
	\[
	\bigl(
	\begin{tikzcd}
	\length\alpha \ar[r,"\overline\sigma"] & N & \ar[l,"\overline\tau"'] \length\beta
	\end{tikzcd},
	\Delta_N  \bigr)
	\]
	where $P_N=\length\alpha + \length\beta$, $\overline\sigma=\injP{\length\alpha}$, $\overline\tau = \injP{\length\beta}$ and for all $p$ place, $\length{\inp p} + \length{\out p} = 1$ (equivalently, $N$ has no internal places and every place is either a proper source or a proper sink). Hence, the general morphism of $\gcf$ is either:
	\begin{itemize}
		\item an identity 
		$
		\bigl( 
		\begin{tikzcd}[cramped,sep=small]
		\length\alpha \ar[r,"\id{}"] & \length\alpha & \ar[l,"\id{}"'] \length\alpha
		\end{tikzcd},
		K_1
		\bigr),
		$
		\item a generator satisfying the conditions above; such morphisms are called \emph{atomic},
		\item a finite composite of atomic morphisms. 
	\end{itemize}
\end{definition}

The assignment $(\alpha,F) \mapsto \alpha$ and 
\[
\bigl[(\phi,
\begin{tikzcd}[cramped,sep=small]
\length\alpha \ar[r,"\sigma"] & n & \length\beta \ar[l,"\tau"']
\end{tikzcd},
\begin{tikzcd}[cramped,sep=small]
\length\alpha \ar[r,"\injP{\length\alpha}"] & \ggraph\phi & \length\beta \ar[l,"\injP{\length\beta}"']
\end{tikzcd},
\Delta_\Phi
)\bigr] 
\mapsto \Bigl[\bigl(
\begin{tikzcd}[cramped,sep=small]
\length\alpha \ar[r,"\injP{\length\alpha}"] & \ggraph\phi & \length\beta \ar[l,"\injP{\length\beta}"']
\end{tikzcd}, \Delta_\Phi \bigr)\Bigr]
\]
mapping atomic morphisms of $\FC \B \C$ to atomic morphisms of $\gcf$ uniquely extends to a functor $\gf \colon \FC \B \C \to \gcf$. Moreover, $\gf$ has two special properties, by virtue of the ``modularity'' of our $\FC \B \C$ and $\gcf$ and the fact that all and only atoms in $\FC \B \C$ have atomic images: it reflects compositions and identities. By ``reflects identities'' we mean that if $\Phi \colon (\alpha,F) \to (\alpha,F)$ is such that $\gf(\Phi)=\id{\length\alpha}$, then $\Phi=\id{(\alpha,F)}$. By ``reflects compositions'' we mean that if $\Phi$ is a morphism in $\FC \B \C$ and $\gf(\Phi)$ is not atomic, i.e.\ $\gf(\Phi) = (N_k,\Delta_k) \circ \dots \circ (N_1,\Delta_1)$ with $(N_i,\Delta_i)$ atomic in $\gcf$, then there must exist $\Phi_1,\dots,\Phi_k$ morphisms in $\FC \B \C$ such that:
\begin{itemize}
	\item $\Phi = \Phi_k \circ \dots \circ \Phi_1$,
	\item $\gf(\Phi_i) = (N_i,\Delta_i)$.
\end{itemize}
Hence, say $\Phi = 
(\phi,
\begin{tikzcd}[cramped,sep=small]
\length\alpha \ar[r,"\sigma"] & n & \length\beta \ar[l,"\tau"']
\end{tikzcd},
\begin{tikzcd}[cramped,sep=small]
\length\alpha \ar[r,"\overline\sigma"] & N & \length\beta \ar[l,"\overline\tau"']
\end{tikzcd},
\Delta_\Phi
)
$: then there must exist transformations $\phi_i$ with graph $\graph{\phi_i}$ (hence atomic), dinatural according to $\Delta_i$, such that $\phi = \phi_k \circ \dots \circ \phi_1$, cf.\ Remark~\ref{remark:non-atomic morphisms of {B,C}}. In other words, $\gf$ satisfies the following definition.

\begin{definition}
	Let $\D,\E$ be any categories. A functor $P \colon \D \to \E$ is said to be a \emph{weak Conduché fibration} (WCF) if, given $f \colon A \to B$ in $\D$:
	\begin{itemize}
		\item $P(f)=\id{}$ implies $f=\id{}$;
		\item given a decomposition $P(f)=u \circ v$ in $\E$, we have that there exist $g,h$ in $\D$ such that $f = g \circ h$, $P(g) = u$, $P(h)=v$.
	\end{itemize}	 
	We define $\WCFover \E$ to be the full subcategory of $\catover\E$ whose objects are the categories over $\E$ whose augmentation is a weak Conduché fibration.
\end{definition}

We have then proved the following theorem.

\begin{theorem}
	$\FC \B \C$ is an object of $\,\,\WCFover\gcf$.
\end{theorem}

Conduché fibrations were introduced in~\cite{conduche_au_1972} as a re-discovery after the original work of Giraud~\cite{giraud_methode_1964} on exponentiable functors in slice categories. Our  notion is weaker in not requiring the additional  property of uniqueness of the decomposition $f=g \circ h$ up to equivalence, where we say that two factorisations $g \circ h$ and $g' \circ h'$ are equivalent if there exists a morphism $j \colon \codom h \to \dom {g'}$ such that everything in sight commutes in the following diagram:
\[
\begin{tikzcd}
& \codom{h} \ar[r,"g"] \ar[d,"j"] & B \\
A \ar[ur,"h"] \ar[r,"h'"'] & \dom{g'} \ar[ur,"g'"']
\end{tikzcd}
\]
We will not, in fact, need such uniqueness; moreover, it is not evident whether our $\gf$ is a Conduché fibration or not.

\begin{remark}
	The fact that $\FC \B \C$ is not just an object of $\catover\gcf$, but even of $\WCFover\gcf$, will allow us to build the substitution category $\ring \A \B$ just for categories $\A$ over $\gcf$ whose augmentation is more than a mere functor: it is a weak Conduché fibration. 
	The main advantage of restricting our attention to $\WCFover\gcf$ is that a category $\A$ in it inherits, in a sense, the modular structure of $\gcf$, as we shall see in the next Lemma.
\end{remark}

\begin{definition}
	Let $P \colon \D \to \gcf$ be an object of $\WCFover\gcf$. A morphism $d$ in $\D$ is said to be \emph{atomic} if $P(d)$ is atomic.
\end{definition}

\begin{lemma}\label{lemma:functors determined by atoms in WCF over E}
	Suppose that, in the following diagram, $P$ is a weak Conduché fibration and $Q$ is an ordinary functor.
	\[
	\begin{tikzcd}[column sep={1cm,between origins}]
	\D \ar[rr,"Q"] \ar[dr,"P"'] & & \mathbb F  \\
	& \gcf
	\end{tikzcd}
	\]
	Then $Q$ is completely determined on morphisms by the image of atomic morphisms of $\D$.
\end{lemma}
\begin{proof}
	Let $d \colon D \to D'$ be a morphism in $\D$ with $P(D)=\alpha$, $P(D')=\beta$ and $P(d) = \bigl[ \bigl(
	\begin{tikzcd}[cramped,sep=small]
	\length\alpha \ar[r,"\overline\sigma"] & N & \length\beta \ar[l,"\overline\tau"']
	\end{tikzcd}, \Delta_d \bigr) \bigr]
	$. If $P(d)$ is not atomic, then either $P(d)=\id{}$, in which case $d=\id{}$ (because $P$ is a weak Conduché fibration), or $P(d)=(N_k,\Delta_k) \circ \dots \circ (N_1,\Delta_1)$ for some (not necessarily unique) atomic $(N_i,\Delta_i)$. Hence there must exist $d_1,\dots,d_k$ in $\D$ such that $d=d_k \circ \dots \circ d_1$ and $P(d_i)=(N_i,\Delta_i)$. Then $Q(d)$ will necessarily be defined as $\id{}$ in the first case, or as $Q(d_k) \circ \dots \circ Q(d_1)$ in the second case, otherwise $Q$ would not be a functor. \qed
\end{proof}

\section{The category of formal substitutions}\label{section:category of formal substitutions}

Kelly~\cite{kelly_many-variable_1972}, after defining his generalised functor category $\fc \B \C$ for covariant functors and many-variable natural transformations only, 
proceeds by showing that the functor $\fc \B -$
has a left adjoint, which he denotes with $\ring - \B$. The category $\ring \A \B$ will be essential to capture the central idea of substitution. 

Here we aim to do the same in our more general setting where $\fc \B \C$ comprises mixed-variance functors and many-variable, partial dinatural transformations. First, we give an explicit definition of the functor $\FC \B - \colon \Cat \to \WCFover\gcf$. Given a functor $K \colon \C \to \C'$, we define $\FC \B K \colon \FC \B \C \to \FC \B {\C'}$ to be the functor mapping $(\alpha,F \colon \B^\alpha \to \C)$ to $(\alpha,KF \colon \B^\alpha \to \C')$; and if 
\[
\Phi = (\phi,
\begin{tikzcd}[cramped,sep=small]
\length\alpha \ar[r,"\sigma"] & n & \length\beta \ar[l,"\tau"']
\end{tikzcd},
\begin{tikzcd}[cramped,sep=small]
\length\alpha \ar[r,"\overline\sigma"] & N & \length\beta \ar[l,"\overline\tau"']
\end{tikzcd},
\Delta_\Phi
)
\colon (\alpha,F) \to (\beta,G)
\]
is a morphism in $\FC \B \C$, then $\FC \B K (\Phi)$ is obtained by whiskering $K$ with $\phi$, obtaining therefore a transformation with the same type and generalised graph as before, with the same dinaturality properties:
\[
\FC \B K (\Phi) = (
K\phi,
\begin{tikzcd}[cramped,sep=small]
\length\alpha \ar[r,"\sigma"] & n & \length\beta \ar[l,"\tau"']
\end{tikzcd},
\begin{tikzcd}[cramped,sep=small]
\length\alpha \ar[r,"\overline\sigma"] & N & \length\beta \ar[l,"\overline\tau"']
\end{tikzcd},
\Delta_\Phi
).
\]
In particular, $\FC \B K$ is clearly a functor over $\gcf$. It is a classic exercise in Category Theory to prove that $\fc \B -$ is continuous, see~\cite[Theorem 3.52]{santamaria_towards_2019}, which is a necessary condition for the existence of a left adjoint
\[
\ring - \B \colon \WCFover\gcf \to \Cat.
\]
We shall prove that a left adjoint does exist by first constructing the category $\ring \A \B$ explicitly, and then showing the existence of a universal arrow $(\A \circ \B, F_\A \colon \A \to \FC \B {\ring \A \B})$ from $\A$ to $\FC \B -$: this will yield the desired adjunction. 

To see what $\ring \A \B$ looks like, we follow Kelly's strategy: we aim to prove that there is a natural isomorphism
\[
\Cat ( \ring \A \B, \C) \cong \WCFover\gcf (\A, \FC \B \C)
\]
and we use this to deduce how $\ring \A \B$ must be. Write $\Gamma$ for all augmentations (as weak Conduché fibrations) over $\gcf$, and let $\Phi$ be an element of $\,\WCFover\gcf(\A,\FC \B \C)$. We now spell out all we can infer from this fact. To facilitate reading, and to comply with Kelly's notation in~\cite{kelly_many-variable_1972}, we shall now refer to the $\bfA$-th component of a transformation $\phi$, for $\bfA=(A_1,\dots,A_m)$ say, as $\phi(\bfA)$ instead of $\phi_{\bfA}$.

\begin{enumerate}[(a),wide,labelindent=0pt]
	\item \label{PhiA} For all $A \in \A$, $\Gamma(A)=\alpha$ we have $\Phi A \colon \B^\alpha \to \C$ is a functor, hence 
	\begin{enumerate}[label=(a.\roman*),wide,leftmargin=\parindent]
		\item for every $\bfB=(B_1,\dots,B_{\length\alpha})$ object of $\B^\alpha$, $\Phi A (\bfB)$ is an object of $\C$,\label{PhiA(B1...Balpha)}
		\item for all $\bfg=(g_1,\dots,g_{\length\alpha})$, with $g_i \colon B_i \to B_i'$ a morphism in $\B$, we have 
		\[
		\Phi(A)(\bfg) \colon  \funminplus {\Phi A} {B_i'} {B_i} i {\length\alpha} \to \funminplus {\Phi A} {B_i} {B_i'} i {\length\alpha}
		\]
		is a morphism in $\C$.\label{PhiA(g1...galpha)}
	\end{enumerate}
	This data is subject to functoriality of $\Phi A$, that is:
	\begin{enumerate}[(1),wide,leftmargin=\parindent]
		\item For every $\bfB$ object of $\B^\alpha$, $\Phi A (\id\bfB) = \id{\Phi A (\bfB)}$\label{PhiA(1...1)=1PhiA}.
		\item For $\bfh=(h_1,\dots,h_{\length\alpha})$, with $h_i \colon B_i' \to B_i''$ morphism of $\B$,
		\[
		\funminplus {\Phi A} {g_i \circ_{\Op\B} h_i} {h_i \circ_{\B} g_i} i {\length\alpha}
		=
		\funminplus {\Phi A} {g_i} {h_i} i {\length\alpha} \circ \funminplus {\Phi A} {h_i} {g_i} i {\length\alpha}.
		\]\label{PhiA(hg)=PhiA(h)PhiA(g)}
	\end{enumerate}
	\item \label{Phif}For all $f \colon A \to A'$ in $\A$ with
	$
	\Gamma(f) = \Bigl[\bigl(
	\begin{tikzcd}[cramped,sep=small]
	\length\alpha \ar[r,"\overline\sigma"] & N & \length\beta \ar[l,"\overline\tau"']
	\end{tikzcd},\Delta_f
	\bigr) \Bigr]
	$, we have that $\Phi f$ is an equivalence class of transformations whose graphs are representatives of $\Gamma(f)$, such transformations being dinatural in some variables according to $\Delta_f$.
	Hence for all $\xi = \bigl((\overline\sigma, \overline\tau, N),\Delta_\xi\bigr) \in \Gamma(f)$ we have a transformation $\Phi f_\xi \colon \Phi A \to \Phi A'$ whose type
	$
	\begin{tikzcd}[cramped,sep=small]
	\length\alpha \ar[r,"\sigma"] & n & \length\beta \ar[l,"\tau"']
	\end{tikzcd}
	$ 
	is the skeleton of $(\overline\sigma,\overline\tau,N)$ and
	with discriminant function $\Delta_\xi$ that tells us in which variables $\Phi f_\xi$ is dinatural. Therefore to give $\Phi f$ one has to provide, for all $\xi = \bigl((\overline\sigma, \overline\tau, N),\Delta_\xi\bigr) \in \Gamma(f)$, for every $\bfB=(B_1,\dots,B_n)$ object of $\B^n$, a morphism in $\C$
	\[
	\Phi f_\xi (\bfB) \colon \Phi A (\bfB\sigma) \to \Phi A' (\bfB\tau)
	\]
	such that:
	\begin{enumerate}[(1),start=3,wide,leftmargin=\parindent]
		\item for all $\pi \colon n \to n$ permutation, $\Phi f_{\pi\xi}(\bfB) = \Phi f_\xi (\bfB\pi)$,\label{Phif_pixi(Bi)=Phif_xi(Bpii)}
		\item \label{Phif_xi dinatural} for $\bfB'=(B_1',\dots,B_n')$ in $\B^n$ and for $\bfg=(g_1,\dots,g_n) \colon \bfB \to \bfB'$ in $\B^n$, where if $\Delta_\xi(i)=0$ then $B_i=B_i'$ and $g_i = \id {B_i}$, the following hexagon commutes:
		\[
		\begin{tikzcd}[font=\normalsize, column sep={0.5cm}]
		& \Phi A (\bfB\sigma) \ar[rrr,"{\Phi f_\xi (\bfB)}"] 
		& &&  \Phi A' (\bfB\tau) \ar[dr,"\funminplus{\Phi A'}{B_{\tau i}}{g_{\tau i}} i {\length\beta}"] \\
		\funminplus{\Phi A}{B_{\sigma i}'}{B_{\sigma i}} i {\length\alpha} \ar[ur,"{\funminplus{\Phi A}{g_{\sigma i}}{B_{\sigma i}} i {\length\alpha}}"]
		\ar[dr,"\funminplus{\Phi A} {{B_{\sigma i}}} {{g_{\sigma i}}} i {\length\alpha}"'] & & && 
		& \funminplus {\Phi A'} {B_{\tau i}}{B_{\tau i}'} i {\length\beta} \\
		& \Phi A (\bfB'\sigma) 
		\ar[rrr,"{\Phi f_\xi (\bfB')}"']
		& & &\Phi A'(\bfB'\tau) 
		\ar[ur,"\funminplus{\Phi A'}{g_{\tau i}}{B_{\tau i}} i {\length\beta}"']
		\end{tikzcd}
		\]
	\end{enumerate}
	\item The data provided in \ref{PhiA} and \ref{Phif} is subject to the functoriality of $\Phi$ itself, hence:
	\begin{enumerate}[(1),start=5,wide,leftmargin=\parindent]
		\item $\Phi(\id A) = \id{\Phi A}$, \label{Phi(1A)=1_Phi(A)}
		\item for $f \colon A \to A'$ and $f' \colon A' \to A''$, $\Phi(f' \circ_{\A} f) = {\Phi f'} \circ_{\FC \B \C} {\Phi f}$ \label{Phi(f2 f1)=Phi(f2) Phi(f1)}.
	\end{enumerate}
\end{enumerate}

We now mirror all the data and properties of a functor $\Phi \colon \A \to \FC \B \C$ over $\gcf$ to define the category $\ring \A \B$. 

\begin{definition}\label{definition A ring B}
	Let $\A$ be a category over $\gcf$ via a weak Conduché fibration $\Gamma \colon \A \to \gcf$, and let $\B$ be any category. The category $\ring \A \B$ of \emph{formal substitutions} of elements of $\B$ into those of $\A$ is the free category generated by the following data. We use the same enumeration as above to emphasise the correspondence between each piece of information.
	\begin{itemize}[wide=0pt,leftmargin=*]
		\item[\ref{PhiA(B1...Balpha)}] Objects are of the form $A[\bfB]$, for $A$ an object of $\A$ with $\Gamma(A)=\alpha$, and for  $\bfB=(B_1,\dots,B_{\length\alpha})$ in $\B^\alpha$. As it is standard in many-variable calculi, we shall drop a set of brackets and write $A[B_1,\dots,B_{\length\alpha}]$ instead of $A[(B_1,\dots,B_{\length\alpha})]$.
		\item[\ref{PhiA(g1...galpha)},\ref{Phif}] Morphisms are to be generated by
		\[
		A[\bfg] \colon \funminplussq A {B_i'} {B_i} i {\length\alpha} \to \funminplussq A {B_i} {B_i'} i {\length\alpha}
		\]
		for $A$ in $\A$ with $\Gamma(A)=\alpha$, $\bfg=(g_1,\dots,g_{\length\alpha})$ and $g_i \colon B_i \to B_i'$ in $\B$, 
		and by
		\[
		f_{\xi}[\bfB] \colon A[\bfB\sigma] \to A'[\bfB\tau]
		\]
		for $f \colon A \to A'$ in $\A$, 
		$
		\xi = \bigl(
		\begin{tikzcd}[cramped,sep=small]
		\length\alpha \ar[r,"\overline\sigma"] & N & \length\beta \ar[l,"\overline\tau"']
		\end{tikzcd},\Delta_\xi \bigr)
		$
		a representative of $\Gamma(f)$, $(\sigma,\tau,n)$ the skeleton of $(\overline\sigma,\overline\tau,N)$, $\bfB=(B_1,\dots,B_n)$ object of $\B^n$. 
	\end{itemize}
	Such data is subject to the following conditions:
	\begin{itemize}[wide=0pt,leftmargin=*]
		\item[\ref{Phif_pixi(Bi)=Phif_xi(Bpii)}] For every permutation $\pi \colon n \to n$ and for every $\bfB=(B_1,\dots,B_n)$ object of $\B^n$
		\[
		f_{\pi\xi}[\bfB] = f_\xi[\bfB\pi].
		\] 
		\item[\ref{PhiA(1...1)=1PhiA},\ref{Phi(1A)=1_Phi(A)}] For all $A\in\A$ with $\Gamma(A)=\alpha$ and for every $\bfB=(B_1,\dots,B_{\length\alpha})$ object of $\B^\alpha$ 
		\[
		A[\id\bfB] = \id{A[\bfB]} = {\id A}[\bfB].
		\]
		\item[\ref{PhiA(hg)=PhiA(h)PhiA(g)}] For all $A \in \A$ with $\Gamma(A)=\alpha$, for all $g_i \colon B_i \to B_i'$ and $h_i \colon B_i' \to B_i''$ in $\B$, $i \in \{1,\dots,\length\alpha\}$
		\[
		\funminplussq { A} {g_i \circ_{\Op\B} h_i} {h_i \circ_{\B} g_i} i {\length\alpha}
		=
		\funminplussq { A} {g_i} {h_i} i {\length\alpha} \circ \funminplussq { A} {h_i} {g_i} i {\length\alpha}.
		\]
		\item[\ref{Phi(f2 f1)=Phi(f2) Phi(f1)}] For all $f \colon A \to A'$ and $f' \colon A' \to A''$ in $\A$, for all
		\[
		\bigl( 
		\begin{tikzcd}[cramped,sep=small]
		\length\alpha \ar[r,"\overline\sigma"] & N & \length\beta \ar[l,"\overline\tau"']
		\end{tikzcd},\Delta
		\bigr) \in \Gamma(f) \quad \text{and} \quad
		\bigl(
		\begin{tikzcd}[cramped,sep=small]
		\length\beta \ar[r,"\overline\eta"] & M & \ar[l,"\overline\theta"'] \length\gamma
		\end{tikzcd},\Delta'
		\bigr) \in \Gamma(f'),
		\]
		with $(\sigma,\tau,n)$ and $(\eta,\theta,m)$ the skeletons of, respectively, $(\overline\sigma,\overline\tau,N)$ and $(\overline\eta,\overline\theta,M)$,
		and for all choices of a pushout
		\[
		\begin{tikzcd}
		& & \length\gamma \ar[d,"\theta"] \\
		& \length\beta \ar[d,"\tau"'] \ar[r,"\eta"] \ar[dr,phantom,very near start,"\ulcorner"] & m \ar[d,"\xi"] \\
		\length\alpha \ar[r,"\sigma"] & n \ar[r,"\zeta"] & l
		\end{tikzcd}
		\]
		each choice determining the skeleton of (the first projection of) a representative of $\Gamma(f' \circ f)$, and for all $\bfB=(B_1,\dots,B_l)$ object of $\B^l$
		\[
		f'_{(\eta,\theta)}[\bfB\xi] \circ f_{(\sigma,\tau)}[\bfB\zeta] = (f'\circ f)_{(\zeta\sigma,\xi\theta)}[\bfB].
		\] 
		\item[\ref{Phif_xi dinatural}] For all $f \colon A \to A'$, $\xi= \bigl(
		\begin{tikzcd}[cramped,sep=small]
		\length\alpha \ar[r,"\overline\sigma"] & N & \length\beta \ar[l,"\overline\tau"']
		\end{tikzcd}, \Delta_\xi \bigr) \in \Gamma(f)
		$,
		with
		$(\sigma,\tau,n)$ the skeleton of $(\overline\sigma,\overline\tau,N)$, for all $\bfB=(B_1,\dots,B_n)$, $\bfB'=(B_1',\dots,B_n')$ objects of $\B^n$ and for all $\bfg=(g_1,\dots,g_n) \colon \bfB \to \bfB'$, with $B_i=B_i'$ and  $g_i=\id{B_i}$ if $\Delta_\xi(i)=0$, the following hexagon commutes:
		\begin{equation}\label{f[g1...gn]}
			\begin{tikzcd}[column sep={0.5cm}]
				& A[\bfB\sigma] \ar[rrr,"{f_\xi [\bfB]}"] 
				& &&  A' [\bfB\tau] \ar[dr,"\funminplussq{A'}{B_{\tau i}}{g_{\tau i}} i {\length\beta}"] \\
				\funminplussq{A}{B_{\sigma i}'}{B_{\sigma i}} i {\length\alpha} \ar[ur,"{\funminplussq{A}{g_{\sigma i}}{B_{\sigma i}} i {\length\alpha}}"]
				\ar[dr,"\funminplussq{A} {{B_{\sigma i}}} {{g_{\sigma i}}} i {\length\alpha}"'] & & && 
				& \funminplussq { A'} {B_{\tau i}}{B_{\tau i}'} i {\length\beta} \\
				& A [\bfB'\sigma] 
				\ar[rrr,"{f_\xi [\bfB']}"']
				& & & A'[\bfB'\tau] 
				\ar[ur,"\funminplussq{A'}{g_{\tau i}}{B_{\tau i}} i {\length\beta}"']
			\end{tikzcd}
		\end{equation}
		We will denote the diagonal of \ref{f[g1...gn]} as $f[\bfg]$. 
	\end{itemize}
\end{definition}

\begin{remark}
	By \ref{Phi(1A)=1_Phi(A)} and \ref{PhiA(hg)=PhiA(h)PhiA(g)}, we have
	\[
	A[\bfg] = \id A [\bfg]
	\]
	and by \ref{PhiA(1...1)=1PhiA}, we have
	\[
	f[\bfB] = f[\id\bfB]
	\]
	which is coherent with the usual notation of $A$ for $\id A$.
\end{remark}

Since two consecutive morphisms both of type \ref{PhiA(g1...galpha)} or both of type \ref{Phif} can be merged together into a single one by \ref{PhiA(hg)=PhiA(h)PhiA(g)} and \ref{Phi(f2 f1)=Phi(f2) Phi(f1)}, we have no way, in general, to swap the order of a morphism of type $A[\bfg]$ followed by one of the form $f_\xi[\bfB]$, because the only axiom that relates the two generators is (\ref{f[g1...gn]}). Therefore, all we can say about the general morphism of $\ring \A \B$ is that it is a string of compositions of alternate morphisms of type \ref{PhiA(g1...galpha)}  and \ref{Phif}, subject to the equations \ref{PhiA(1...1)=1PhiA}-\ref{Phi(f2 f1)=Phi(f2) Phi(f1)}. 

\begin{remark}
	If $\A$ is such that $\length{\Gamma(A)}=1$ for all objects $A$ in $\A$, then $\ring \A \B$ is highly reminiscent of the category $\A \otimes \B$ as described by Power and Robinson in~\cite{power_premonoidal_1997}. The authors studied the \emph{other} symmetric monoidal closed structure of $\Cat$, where the exponential $[\B,\C]$ is the category of functors from $\B$ to $\C$ and morphisms are simply transformations (not necessarily natural), and $\otimes \B$ is the tensor functor that is the left adjoint of $[\B,-]$. The category $\A \otimes \B$ has pairs $(A,B)$ of objects of $\A$ and $\B$, and a morphism from $(A,B)$ to $(A',B')$ is a finite sequence of non-identity arrows consisting of alternate chains of consecutive morphisms of $\A$ and $\B$. Composition is given by concatenation followed by cancellation accorded by the composition in $\A$ and $\B$, much like our $\ring \A \B$. The only difference with their case is that we have the additional dinaturality equality \ref{Phif_xi dinatural}. For an arbitrary category $\A$ over $\gcf$, our $\ring \A \B$ would be a sort of generalised tensor product, where the number of objects of $\B$ we ``pair up'' with an object $A$ of $\A$ depends on $\Gamma(A)$.
\end{remark}

We are now ready to show that $\fc \B -$ has indeed a left adjoint. This is going to be a crucial step towards a complete substitution calculus for dinatural transformations; we shall discuss some ideas and conjectures about the following steps in the conclusions.

\begin{theorem}\label{theorem:{B,-} has a left adjoint}
	The functor $\FC \B -$ has a left adjoint
	\[
	\begin{tikzcd}[column sep=2cm,bend angle=30]
	{\WCFover\gcf} \ar[r,bend left,"\ring - \B"{name=A},pos=.493] & \Cat \ar[l,bend left,"\FC \B -"{name=B},pos=.507]	
	\ar[from=A,to=B,phantom,"\bot"]
	\end{tikzcd}
	\]
	therefore there is a natural isomorphism
	\begin{equation}\label{natural isomorphism (A circ B, C) -> (A,{B,C})}
		\Cat \bigl( \ring \A \B , \C \bigr) \cong \WCFover\gcf \bigl( \A, \FC \B \C \bigr).
	\end{equation}
	Moreover, $\ring {} {} \colon \WCFover\gcf \times \Cat \to \Cat$ is a functor.
\end{theorem}
\begin{proof}
	Recall that to give an adjunction $ (\ring - \B) \dashv \FC \B -$ is equivalent to give, for all $\A \in \WCFover\gcf$, a universal arrow $(\ring \A \B, F_\A \colon \A \to \FC \B {\ring \A \B})$ from $\A$ to the functor $\FC \B -$; $F_\A$ being a morphism of $\WCFover\gcf$. This means that, for a fixed $\A$, we have to define a functor over $\gcf$ that makes the following triangle commute:
	\[
	\begin{tikzcd}[column sep={1.5cm,between origins}]
	\A \ar[rr,"F_\A"] \ar[dr,"\Gamma"'] & & \FC \B {\ring \A \B} \ar[dl,"\gf"] \\
	& \gcf
	\end{tikzcd}
	\]
	and that is universal among all arrows from $\A$ to $\FC \B -$: for all arrows $(\C, \Phi \colon \A \to \FC \B \C)$ from $\A$ to $\FC \B -$ ($\Phi$ being a functor over $\gcf$), there must exist a unique morphism in $\Cat$, that is a functor, $H \colon \ring \A \B \to \C$ such that
	\[
	\begin{tikzcd}
	\A \ar[r,"F_\A"] \ar[dr,"\Phi"'] & \FC \B {\ring \A \B} \ar[d,"{\FC \B H}"] \\
	& \FC \B \C 
	\end{tikzcd}
	\]
	commutes. In the proof we will refer to properties \ref{PhiA(1...1)=1PhiA}-\ref{Phi(f2 f1)=Phi(f2) Phi(f1)} as given in the definition of $\ring \A \B$.
	
	Let then $\A$ be a category over $\gcf$ with $\Gamma \colon \A \to \gcf$ a weak Conduché fibration. We define the action of $F_\A$ on objects first. If $A$ is an object of $\A$ with $\Gamma(A)=\alpha$, then the assignment
	\[
	\begin{tikzcd}[row sep=0em]
	\B^{\alpha} \ar[r,"F_\A(A)"] & \ring \A \B \\
	\bfB \ar[|->,r] \ar{d}[description,name=A]{\bfg} & A[\bfB] \ar{d}[description,name=B]{{A[\bfg]}}     \\[2em]
	\bfB' \ar[|->,r] & A[\bfB']
	\arrow[from=A,to=B,|->]
	\end{tikzcd}
	\] 
	is a functor by virtue of \ref{PhiA(1...1)=1PhiA} and \ref{PhiA(hg)=PhiA(h)PhiA(g)}. By little abuse of notation, call $F_\A(A)$ also the pair $(\alpha,F_\A(A))$, which is an object of $\FC \B {\ring \A \B}$.
	
	To define $F_\A$ on morphisms, let $f \colon A \to A'$ be a morphism in $\A$, with $\Gamma(A)=\alpha$, $\Gamma(A')=\beta$, let
	\[
	\xi = \bigl(
	\begin{tikzcd}[cramped,sep=small]
	\length\alpha \ar[r,"\overline\sigma"] & N & \length\beta \ar[l,"\overline\tau"']
	\end{tikzcd},
	\Delta_\xi
	\bigr) \in \Gamma(f),
	\]
	and call 
	$
	\begin{tikzcd}[cramped,sep=small]
	\length\alpha \ar[r,"\sigma"] & n & \length\beta \ar[l,"\tau"']
	\end{tikzcd}
	$
	the skeleton of $(\overline\sigma,\overline\tau,N)$. We define $F_\A (f) \colon F_\A(A) \to F_\A(A')$ to be the equivalent class of the tuple
	\[
	\bigl(
	F_\A (f)_\xi,
	\begin{tikzcd}[cramped,sep=small]
	\length\alpha \ar[r,"\sigma"] & n & \length\beta \ar[l,"\tau"']
	\end{tikzcd},
	\begin{tikzcd}[cramped,sep=small]
	\length\alpha \ar[r,"\overline\sigma"] & N & \length\beta \ar[l,"\overline\tau"']
	\end{tikzcd},
	\Delta_\xi
	\bigr)
	\]
	where $F_\A(f)_\xi$ is a transformation whose general component is
	\[
	\begin{tikzcd}[row sep=1em,column sep=4em]
	F_\A(A)(\bfB\sigma) \ar[d,phantom,"\rotatebox{90}="] \ar[r,"{f_\xi[\bfB]}"] & F_\A(A')(\bfB\tau) \ar[d,phantom,"\rotatebox{90} ="] \\
	A[\bfB\sigma] & A'[\bfB\tau]
	\end{tikzcd}
	\]
	Then $F_\A(f)_\xi$ is indeed dinatural in its $i$-th variable whenever $\Delta_\xi(i)=1$ because of~\ref{Phif_xi dinatural}.  Moreover, $F_\A$ is well-defined on morphisms because of \ref{Phif_pixi(Bi)=Phif_xi(Bpii)} and is in fact a functor thanks to \ref{Phi(1A)=1_Phi(A)} and \ref{Phi(f2 f1)=Phi(f2) Phi(f1)}. Finally, $F_\A(f)$ so defined is indeed a morphism of $\FC \B {\ring \A \B}$: if $f$ is such that $\Gamma(f)$ is atomic, then $F_\A(f)$ is an atomic morphism of $\FC \B {\ring \A \B}$; if instead $\Gamma(f)=(N_k,\Delta_k) \circ \dots \circ (N_1,\Delta_1)$ where $(N_i,\Delta_i)$ is atomic, then there exists a factorisation $f=f_k \circ \dots \circ f_1$ in $\A$ with $\Gamma(f_i)=(N_i,\Delta_i)$ because $\Gamma$ is a weak Conduché fibration. By functoriality of $F_\A$, we have that $F_\A(f)=F_\A(f_k) \circ \dots \circ F_\A(f_1)$, hence it is a composite of atomic morphisms of $\FC \B {\ring \A \B}$.
	
	We now prove that $F_\A$ is universal. Let then $\Phi \colon \A \to \FC \B \C$ be a morphism in $\WCFover\gcf$, that is a functor over $\gcf$. We define $H \colon \ring \A \B \to \C$ as follows:
	\begin{itemize}[wide=0pt,leftmargin=*]
		\item[\ref{PhiA(B1...Balpha)}] For $A \in \A$ with $\Gamma(A)=\alpha$ and $\bfB \in \B^\alpha$, 
		\[
		H\bigl(A[\bfB]\bigr) = \Phi(A)(\bfB);
		\]
		\item[\ref{PhiA(g1...galpha)}] For $A \in \A$ with $\Gamma(A)=\alpha$, for $\bfg$ in $\B^\alpha$, 
		\[
		H\bigl(A[\bfg]\bigr) = \Phi(A)(\bfg);
		\]
		\item[\ref{Phif}] For $f \colon A \to A'$ in $\A$, $\xi = (N_\xi,\Delta_\xi) \in \Gamma(f)$ where $N_\xi$ has $n$ connected components, for $\bfB \in \B^n$,
		\[
		H\bigl(f_\xi[\bfB]\bigr) = \Phi(f)_\xi(\bfB),
		\]
		where $\Phi(f)_\xi$ is the representative of $\Phi(f)$ whose type is given by the skeleton of $N_\xi$, cf.~the discussion on the data entailed by a functor $\Phi \colon \A \to \FC \B \C$ over $\gcf$ preceding Definition~\ref{definition A ring B}.
	\end{itemize}
	$H$ so defined on the generators of $\ring \A \B$ extends to a unique functor provided that $H$ preserves the equalities \ref{PhiA(1...1)=1PhiA}-\ref{Phi(f2 f1)=Phi(f2) Phi(f1)} in $\ring \A \B$, which it does as they have been designed \emph{precisely} to reflect all the properties of a functor $\Phi \colon \A \to \FC \B \C$, and $H$ is defined using $\Phi$ accordingly. Finally, by construction 
	\[
	\begin{tikzcd}
	\A \ar[r,"F_\A"] \ar[dr,"\Phi"'] & \FC \B {\ring \A \B} \ar[d,"{\FC \B H}"] \\
	& \FC \B \C 
	\end{tikzcd}
	\]
	commutes. The uniqueness of $H$ follows from the fact that the commutativity of the above triangle implies that $\Phi(A)=H(F_\A(A))$ for all $A \in \A$ and $\Phi(f)=H(F_\A(f))$, hence any such functor $H$ \emph{must} be defined as we did to make the triangle commutative. 
	
	With such a universal arrow $(\ring \A \B, F_\A \colon \A \to \FC \B {\ring \A \B})$ we can define a functor $\ring - \B$ which is the left adjoint of $\FC \B -$. Given $F \colon \A \to \A'$ a functor over $\gcf$, by universality of $F_\A$ there exists a unique functor $\ring F \B \colon \ring \A \B \to \ring {\A'} \B$ that makes the following square commute:
	\[
	\begin{tikzcd}
	\A \ar[r,"F_\A"] \ar[d,"F"'] & \FC \B {\ring \A \B} \ar[d,"\FC \B {\ring F \B}"] \\
	\A' \ar[r,"F_{\A'}"'] & \FC \B {\ring {\A'} \B}
	\end{tikzcd}
	\]
	Such $\ring F \B$ is defined on objects as $\ring F \B \bigl( A[\bfB] \bigr) = (F_{\A'} \circ F)(A)(\bfB) = FA[\bfB]$ and on morphisms as
	\[
	\ring F \B \bigl( A[\bfg] \bigr) = FA[\bfg], \quad 
	\ring F \B \bigl( f[\bfB] \bigr) = Ff[\bfB].
	\]
	
	Finally, $\ring{}{}$ extends to a functor
	\[
	\begin{tikzcd}[row sep=0em]
	\WCFover\gcf \times \Cat \ar[r,"\ring{}{}"] & \Cat \\
	\A\quad\quad\,\B \ar[r,|->] \ar[d,shift right=17mu,"F"'] \ar[d,shift left=17mu,"G"] & \ring \A \B \ar[d,"\ring F G"] \\[3em] 
	{\A'}\quad\quad\B' \ar[r,|->] & \ring {\A'} {\B'}
	\end{tikzcd}
	\]
	where $\ring F G$ is defined as follows on the generators:
	\begin{itemize}
		\item $\ring F G \bigl( A[\bfB] \bigr) = FA[G\bfB]$,
		\item $\ring F G \bigl( A[\bfg] \bigr) = FA[G\bfg]$,
		\item $\ring F G \bigl( f[\bfB] \bigr) = Ff[G\bfB]$
	\end{itemize}
	(where $G\bfB=(GB_1,\dots,GB_{\length\alpha})$ if $\bfB=(B_1,\dots,B_{\length\alpha})$). It is easy to see that $\ring F G$ is well defined (i.e.\ it preserves equalities in $\ring \A \B$), thanks to the functoriality of $F$ and $G$. It is also immediate to verify that   $\ring{}{}$ is indeed a functor.\qed
\end{proof}

\section{Conclusions}\label{section:coda}

The ultimate goal to achieve a complete substitution calculus of dinatural transformations is to obtain an appropriate functor over $\gcf$
\[
M \colon \ring{\FC \B \C} {\FC \A \B} \to \FC \A \C
\]
which, \emph{de facto}, realises a \emph{formal} substitution of functors into functors and transformations into transformations as an \emph{actual} new functor or transformation. As in Kelly's case, {horizontal} composition of dinatural transformations will be at the core, we believe, of the desired functor; the rules of vertical composition are, instead, already embodied into the definition of $\FC \B \C$.

Such $M$ will arise as a consequence of proving that $\WCFover\gcf$ is a monoidal closed category, much like Kelly did, by showing that the natural isomorphism~(\ref{natural isomorphism (A circ B, C) -> (A,{B,C})}) extends to
\[
\WCFover\gcf (\ring \A \B , \C) \cong \WCFover\gcf (\A , \FC \B \C).
\]
Necessarily then, we will first have to show that the substitution category $\ring \A \B$ is itself an object of $\WCFover\gcf$. Following Kelly's steps described in~\cite[\S 2.1]{kelly_many-variable_1972}, this will be done by extending our functor $\ring{}{} \colon \WCFover\gcf \times \Cat \to \Cat$ to a functor
\[
\ring{}{} \colon \WCFover\gcf \times \WCFover\gcf \to \WCFover\gcf,
\]
exhibiting $\WCFover\gcf$ as a monoidal category, with tensor $\ring{}{}$. 
To do so in his case, Kelly defined $\ring \A \B$ just as before, ignoring the augmentation on $\B$, and then augmented $\ring \A \B$ using the augmentations of $\A$ and $\B$. In fact, what he did, using the category $\Per$ of permutations, was to regard $\Per$ as a category over itself in the obvious way and then to define a functor $P \colon \ring \Per \Per \to \Per$ that computes substitution of permutations into permutations. That done, he set $\Gamma \colon \ring \A \B \to \Per$ as a composite
\[
\begin{tikzcd}
\ring \A \B \ar[d,"\ring {\Gamma_\A} {\Gamma_\B}"'] \ar[r] & \Per \\
\ring \Per \Per \ar[ur,"P"']
\end{tikzcd}
\]
This suggests, as usual, to do the same in our case. Hence, the next step will be to come up with a substitution functor
\[
S \colon \ring \gcf  \gcf \to \gcf,
\]
which is tantamount to define an operation of substitution of graphs, and then define $\Gamma \colon \ring \A \B \to \gcf$ as
\begin{equation}\label{augmentation of A ring B via G ring G}
	\begin{tikzcd}
		\ring \A \B \ar[d,"\ring {\Gamma_\A} {\Gamma_\B}"'] \ar[r] & \gcf\\
		\ring \gcf \gcf \ar[ur,"S"']
	\end{tikzcd}
\end{equation}

A possible hint to how to do this is given by how we defined the horizontal composition of dinatural transformations in Chapter~\ref{chapter horizontal}, and what happened to the graphs of the transformations (that is, we consider the special case of $\A = \B = \FC \C \C$). Looking back at Example~\ref{ex:hc example}, when we computed the first horizontal composition of $\delta$ and $(\eval A B)_{A,B}$, in fact we considered the formal substitution $\eval{}{}\bigl[\delta,([+],\id\C)\bigr]$ in $\ring {\FC \C \C} {\FC \C \C}$, which we then realised into the transformation $\HC \delta {\eval{}{}} 1$. This realisation part is what the desired functor $M$ will do, once properly defined. Now, consider, in $\ring \gcf \gcf$, the formal substitution $\graph{\eval{}{}}\bigl[\graph\delta,[+]\bigr]$, which is the image of $\eval{}{}\bigl[\delta,([+],\id\C)\bigr]$ along the functor $\ring \gf \gf \colon \ring {\FC \C \C} {\FC \C \C} \to \ring \gcf \gcf$. Since $M \colon \ring {\FC \C \C} {\FC \C \C}$ ought to be a functor over $\gcf$, we have that $S\bigl(\graph{\eval{}{}}\bigl[\graph\delta,[+]\bigr]\bigr)$ should be the graph that $\HC \delta {\eval{}{}} 1$ has, which is
\[
\begin{tikzpicture}
\matrix[column sep=1em,row sep=1.5em]{
	\node[category] (1) {}; & \node[opCategory] (2) {}; & \node[opCategory] (3) {}; & \node[category] (4) {}; \\
	& \node[component] (A) {}; & & \node[component] (B) {};\\
	& & & \node[category] (5) {};\\
};
\graph[use existing nodes]{
	1 -> A -> {2,3};
	4 -> B -> 5;
};
\end{tikzpicture}
\]
The intuition for it was that we ``bent'' $\graph\delta$ into the U-turn that is the first connected component of $\graph{\eval{}{}}$. A possible approach to a general definition of substitution of graphs into graphs is the following: given two connected graphs $N_1$, $N_2$ in $\gcf$, the graph $S\bigl(N_1[N_2]\bigr)$ is the result of subjecting $N_2$ to all the ramifications and U-turns of $N_1$; in so doing, one would have to substitute a copy of $N_2$ in every \emph{directed path} of $N_1$. 
This idea is not original, as it was suggested by Bruscoli, Guglielmi, Gundersen and Parigot~\cite{guglielmi_substitution} in private communications to implement substitution of \emph{atomic flows}~\cite{GuglGundStra::Breaking:uq}, which are graphs extracted from certain formal proofs in \emph{Deep Inference}~\cite{Gugl:06:A-System:kl} and they look very much like a morphism in $\gcf$. 

How to put such an intuitive idea into a formal, working definition is the subject of current investigations, and this task has already revealed itself as far from being trivial. Once that is done, the rest should follow relatively easily, and we would expect that the correct compatibility law for horizontal and vertical composition sought in \ref{section compatibility} will become apparent, once the substitution functor $M$ above will be found as part of a monoidal closed structure.

\section*{Acknowledgements}
Most of the material in this article derives from Santamaria's PhD thesis~\cite{santamaria_towards_2019}, written under the supervision of McCusker, and it is, in part, a detailed version of~\cite{mccusker_compositionality_2018}. As such, Santamaria acknowledges the support of a research studentship from the
University of Bath, as well as EPSRC grant EP/R006865/1 and the funding support of the Ministero dell’Universit\`a e della Ricerca of Italy under Grant No.~201784YSZ5, PRIN2017.

The authors would like to thank John Power for suggesting the notations to handle the manipulation of tuples, which we believe provided a great improvement to the exposition of our theory with respect to~\cite{mccusker_compositionality_2018,santamaria_towards_2019}. We would also like to thank Alessio Guglielmi for his valuable insights on the simplification of the proof of Theorem~\ref{thm:acyclic-implies-reachable} with respect to~\cite{mccusker_compositionality_2018,santamaria_towards_2019}.

Finally, we thank Zoran Petri\'c for his kind understanding of our lack of acknowledgement of his results in the past: we hope that with this paper we have finally given him the credit he deserves for his work.

\bibliography{JPAAbiblio}

\begin{thebibliography}{36}
\providecommand{\natexlab}[1]{#1}
\providecommand{\url}[1]{\texttt{#1}}
\expandafter\ifx\csname urlstyle\endcsname\relax
  \providecommand{\doi}[1]{doi: #1}\else
  \providecommand{\doi}{doi: \begingroup \urlstyle{rm}\Url}\fi

\bibitem[Bainbridge et~al.(1990)Bainbridge, Freyd, Scedrov, and
  Scott]{bainbridge_functorial_1990}
E.~S. Bainbridge, P.~J. Freyd, A.~Scedrov, and P.~J. Scott.
\newblock Functorial polymorphism.
\newblock \emph{Theoretical Computer Science}, 70\penalty0 (1):\penalty0
  35--64, Jan. 1990.
\newblock ISSN 0304-3975.
\newblock \doi{10.1016/0304-3975(90)90151-7}.

\bibitem[Blute(1993)]{blute_linear_1993}
R.~Blute.
\newblock Linear logic, coherence and dinaturality.
\newblock \emph{Theoretical Computer Science}, 115\penalty0 (1):\penalty0
  3--41, July 1993.
\newblock ISSN 0304-3975.
\newblock \doi{10.1016/0304-3975(93)90053-V}.

\bibitem[Bruscoli et~al.(Personal communication)Bruscoli, Guglielmi, Gundersen,
  and Parigot]{guglielmi_substitution}
P.~Bruscoli, A.~Guglielmi, T.~E. Gundersen, and M.~Parigot.
\newblock Proposal for substitution in {Deep} {Inference} via atomic flows,
  Personal communication.

\bibitem[Conduché(1972)]{conduche_au_1972}
F.~Conduché.
\newblock Au sujet de l’existence d’adjoints à droite aux foncteurs
  ‘image reciproque’ dans la catégorie des catégories.
\newblock \emph{C. R. Acad. Sci. Paris}, \penalty0 (275):\penalty0 891--894,
  1972.

\bibitem[Dubuc and Street(1970)]{dubuc_dinatural_1970}
E.~Dubuc and R.~Street.
\newblock Dinatural transformations.
\newblock In S.~Mac~Lane, H.~Applegate, M.~Barr, B.~Day, E.~Dubuc, A.~P.
  Phreilambud, R.~Street, M.~Tierney, and S.~Swierczkowski, editors,
  \emph{Reports of the {Midwest} {Category} {Seminar} {IV}}, volume 137 of
  \emph{Lecture {Notes} in {Mathematics}}, pages 126--137. Springer, Berlin,
  Heidelberg, 1970.
\newblock ISBN 978-3-540-04926-5.
\newblock \doi{10.1007/BFb0060443}.

\bibitem[Eilenberg and Kelly(1966)]{eilenberg_generalization_1966}
S.~Eilenberg and G.~M. Kelly.
\newblock A {Generalization} of the {Functorial} {Calculus}.
\newblock \emph{Journal of Algebra}, 3\penalty0 (3):\penalty0 366--375, May
  1966.
\newblock ISSN 0021-8693.
\newblock \doi{10.1016/0021-8693(66)90006-8}.

\bibitem[Freyd et~al.(1992)Freyd, Robinson, and
  Rosolini]{freyd_dinaturality_1992}
P.~J. Freyd, E.~P. Robinson, and G.~Rosolini.
\newblock Dinaturality for free.
\newblock In A.~M. Pitts, M.~P. Fourman, and P.~T. Johnstone, editors,
  \emph{Applications of {Categories} in {Computer} {Science}: {Proceedings} of
  the {London} {Mathematical} {Society} {Symposium}, {Durham} 1991}, London
  {Mathematical} {Society} {Lecture} {Note} {Series}, pages 107--118. Cambridge
  University Press, Cambridge, 1992.
\newblock ISBN 978-0-521-42726-5.
\newblock \doi{10.1017/CBO9780511525902.007}.

\bibitem[Girard et~al.(1992)Girard, Scedrov, and Scott]{girard_normal_1992}
J.-Y. Girard, A.~Scedrov, and P.~J. Scott.
\newblock Normal {Forms} and {Cut}-{Free} {Proofs} as {Natural}
  {Transformations}.
\newblock In N.~M. Yiannis, editor, \emph{Logic from {Computer} {Science}},
  volume~21 of \emph{Mathematical {Sciences} {Research} {Institute}
  {Publications}}, pages 217--241. Springer, New York, NY, 1992.
\newblock ISBN 978-1-4612-7685-2 978-1-4612-2822-6.
\newblock \doi{10.1007/978-1-4612-2822-6_8}.

\bibitem[Giraud(1964)]{giraud_methode_1964}
J.~Giraud.
\newblock Méthode de la descente.
\newblock \emph{Mémoires de la Société Mathématique de France}, 2:\penalty0
  156, 1964.
\newblock ISSN 0249-633X, 2275-3230.
\newblock \doi{10.24033/msmf.2}.

\bibitem[Godement(1958)]{godement_topologie_1958}
R.~Godement.
\newblock \emph{Topologie algébrique et théorie des faisceaux}.
\newblock Number~13 in Publications de l'{Institut} de mathématique de
  l'{Université} de {Strasbourg}. Hermann, Paris, 1958.
\newblock OCLC: 1216318.

\bibitem[Guglielmi(2007)]{Gugl:06:A-System:kl}
A.~Guglielmi.
\newblock A system of interaction and structure.
\newblock \emph{ACM Transactions on Computational Logic}, 8\penalty0
  (1):\penalty0 1:1--64, 2007.
\newblock \doi{10.1145/1182613.1182614}.

\bibitem[Guglielmi et~al.(2010)Guglielmi, Gundersen, and
  Stra{\ss}burger]{GuglGundStra::Breaking:uq}
A.~Guglielmi, T.~Gundersen, and L.~Stra{\ss}burger.
\newblock Breaking paths in atomic flows for classical logic.
\newblock In J.-P. Jouannaud, editor, \emph{25th Annual IEEE Symposium on Logic
  in Computer Science (LICS)}, pages 284--293. IEEE, 2010.
\newblock \doi{10.1109/LICS.2010.12}.

\bibitem[Hiraishi and Ichikawa(1988)]{hiraishi_class_1988}
K.~Hiraishi and A.~Ichikawa.
\newblock A {Class} of {Petri} {Nets} {That} a {Necessary} and {Sufficient}
  {Condition} for {Reachability} is {Obtainable}.
\newblock \emph{Transactions of the Society of Instrument and Control
  Engineers}, 24\penalty0 (6):\penalty0 635--640, June 1988.
\newblock ISSN 0453-4654, 1883-8189.
\newblock \doi{10.9746/sicetr1965.24.635}.

\bibitem[Joyal and Street(1991)]{joyal_geometry_1991}
A.~Joyal and R.~Street.
\newblock The geometry of tensor calculus, {I}.
\newblock \emph{Advances in Mathematics}, 88\penalty0 (1):\penalty0 55--112,
  July 1991.
\newblock ISSN 00018708.
\newblock \doi{10.1016/0001-8708(91)90003-P}.

\bibitem[Joyal et~al.(1996)Joyal, Street, and Verity]{joyal_traced_1996}
A.~Joyal, R.~Street, and D.~Verity.
\newblock Traced monoidal categories.
\newblock \emph{Mathematical Proceedings of the Cambridge Philosophical
  Society}, 119\penalty0 (3):\penalty0 447--468, Apr. 1996.
\newblock ISSN 1469-8064, 0305-0041.
\newblock \doi{10.1017/S0305004100074338}.

\bibitem[Kelly(1972{\natexlab{a}})]{kelly_abstract_1972}
G.~M. Kelly.
\newblock An abstract approach to coherence.
\newblock In G.~M. Kelly, M.~Laplaza, G.~Lewis, and S.~Mac~Lane, editors,
  \emph{Coherence in {Categories}}, volume 281 of \emph{Lecture {Notes} in
  {Mathematics}}, pages 106--147. Springer, Berlin, Heidelberg,
  1972{\natexlab{a}}.
\newblock ISBN 978-3-540-05963-9 978-3-540-37958-4.
\newblock \doi{10.1007/BFb0059557}.

\bibitem[Kelly(1972{\natexlab{b}})]{kelly_many-variable_1972}
G.~M. Kelly.
\newblock Many-variable functorial calculus. {I}.
\newblock In G.~M. Kelly, M.~Laplaza, G.~Lewis, and S.~Mac~Lane, editors,
  \emph{Coherence in {Categories}}, volume 281 of \emph{Lecture {Notes} in
  {Mathematics}}, pages 66--105. Springer, Berlin, Heidelberg,
  1972{\natexlab{b}}.
\newblock ISBN 978-3-540-05963-9.
\newblock \doi{10.1007/BFb0059556}.

\bibitem[Kelly and Laplaza(1980)]{kelly_coherence_1980}
G.~M. Kelly and M.~L. Laplaza.
\newblock Coherence for compact closed categories.
\newblock \emph{Journal of Pure and Applied Algebra}, 19:\penalty0 193--213,
  Dec. 1980.
\newblock ISSN 0022-4049.
\newblock \doi{10.1016/0022-4049(80)90101-2}.

\bibitem[Kosaraju(1982)]{kosaraju_decidability_1982}
S.~R. Kosaraju.
\newblock Decidability of {Reachability} in {Vector} {Addition} {Systems}.
\newblock In \emph{Proceedings of the {Fourteenth} {Annual} {ACM} {Symposium}
  on {Theory} of {Computing}}, {STOC} '82, pages 267--281, New York, NY, USA,
  1982. ACM.
\newblock ISBN 978-0-89791-070-5.
\newblock \doi{10.1145/800070.802201}.

\bibitem[Lataillade(2009)]{lataillade_dinatural_2009}
J.~d. Lataillade.
\newblock Dinatural {Terms} in {System} {F}.
\newblock In \emph{2009 24th {Annual} {IEEE} {Symposium} on {Logic} {In}
  {Computer} {Science}}, pages 267--276, Aug. 2009.
\newblock \doi{10.1109/LICS.2009.30}.

\bibitem[Mac~Lane(1963)]{mac_lane_natural_1963}
S.~Mac~Lane.
\newblock Natural {Associativity} and {Commutativity}.
\newblock \emph{Rice Institute Pamphlet - Rice University Studies}, 49\penalty0
  (4), Oct. 1963.

\bibitem[Mac~Lane(1978)]{mac_lane_categories_1978}
S.~Mac~Lane.
\newblock \emph{Categories for the {Working} {Mathematician}}, volume~5 of
  \emph{Graduate {Texts} in {Mathematics}}.
\newblock Springer-Verlag, New York, 2 edition, 1978.
\newblock ISBN 978-0-387-98403-2.

\bibitem[Mayr(1981)]{mayr_algorithm_1981}
E.~W. Mayr.
\newblock An {Algorithm} for the {General} {Petri} {Net} {Reachability}
  {Problem}.
\newblock In \emph{Proceedings of the {Thirteenth} {Annual} {ACM} {Symposium}
  on {Theory} of {Computing}}, {STOC} '81, pages 238--246, New York, NY, USA,
  1981. ACM.
\newblock \doi{10.1145/800076.802477}.

\bibitem[McCusker and Santamaria(2018)]{mccusker_compositionality_2018}
G.~McCusker and A.~Santamaria.
\newblock On {Compositionality} of {Dinatural} {Transformations}.
\newblock In D.~Ghica and A.~Jung, editors, \emph{27th {EACSL} {Annual}
  {Conference} on {Computer} {Science} {Logic} ({CSL} 2018)}, volume 119 of
  \emph{Leibniz {International} {Proceedings} in {Informatics} ({LIPIcs})},
  pages 33:1--33:22, Dagstuhl, Germany, 2018. Schloss
  Dagstuhl–Leibniz-Zentrum fuer Informatik.
\newblock ISBN 978-3-95977-088-0.
\newblock \doi{10.4230/LIPIcs.CSL.2018.33}.

\bibitem[Mulry(1990)]{mulry_categorical_1990}
P.~S. Mulry.
\newblock Categorical fixed point semantics.
\newblock \emph{Theoretical Computer Science}, 70\penalty0 (1):\penalty0
  85--97, Jan. 1990.
\newblock ISSN 0304-3975.
\newblock \doi{10.1016/0304-3975(90)90154-A}.

\bibitem[Murata(1989)]{murata_petri_1989}
T.~Murata.
\newblock Petri nets: {Properties}, analysis and applications.
\newblock \emph{Proceedings of the IEEE}, 77\penalty0 (4):\penalty0 541--580,
  Apr. 1989.
\newblock ISSN 0018-9219.
\newblock \doi{10.1109/5.24143}.

\bibitem[Paré and Román(1998)]{pare_dinatural_1998}
R.~Paré and L.~Román.
\newblock Dinatural numbers.
\newblock \emph{Journal of Pure and Applied Algebra}, 128\penalty0
  (1):\penalty0 33--92, June 1998.
\newblock ISSN 0022-4049.
\newblock \doi{10.1016/S0022-4049(97)00036-4}.

\bibitem[Petri(1962)]{petri_kommunikation_1962}
C.~A. Petri.
\newblock \emph{Kommunikation mit {Automaten}}.
\newblock PhD thesis, Mathematisches Institut der Universität Bonn, Bonn,
  1962.
\newblock OCLC: 258511501.

\bibitem[Petrić(2003)]{petric_g-dinaturality_2003}
Z.~Petrić.
\newblock G-dinaturality.
\newblock \emph{Annals of Pure and Applied Logic}, 122\penalty0 (1):\penalty0
  131--173, Aug. 2003.
\newblock ISSN 0168-0072.
\newblock \doi{10.1016/S0168-0072(03)00003-4}.

\bibitem[Pistone(2017)]{pistone_dinaturality_2017}
P.~Pistone.
\newblock On {Dinaturality}, {Typability} and beta-eta-{Stable} {Models}.
\newblock In D.~Miller, editor, \emph{2nd {International} {Conference} on
  {Formal} {Structures} for {Computation} and {Deduction} ({FSCD} 2017)},
  volume~84 of \emph{Leibniz {International} {Proceedings} in {Informatics}
  ({LIPIcs})}, pages 29:1--29:17, Dagstuhl, Germany, 2017. Schloss
  Dagstuhl–Leibniz-Zentrum fuer Informatik.
\newblock ISBN 978-3-95977-047-7.
\newblock \doi{10.4230/LIPIcs.FSCD.2017.29}.

\bibitem[Plotkin and Abadi(1993)]{plotkin_logic_1993}
G.~Plotkin and M.~Abadi.
\newblock A logic for parametric polymorphism.
\newblock In M.~Bezem and J.~F. Groote, editors, \emph{Typed {Lambda} {Calculi}
  and {Applications}}, Lecture {Notes} in {Computer} {Science}, pages 361--375.
  Springer Berlin Heidelberg, 1993.
\newblock ISBN 978-3-540-47586-6.

\bibitem[Power and Robinson(1997)]{power_premonoidal_1997}
J.~Power and E.~Robinson.
\newblock Premonoidal categories and notions of computation.
\newblock \emph{Mathematical Structures in Computer Science}, 7\penalty0
  (5):\penalty0 453--468, Oct. 1997.
\newblock ISSN 1469-8072, 0960-1295.
\newblock \doi{10.1017/S0960129597002375}.

\bibitem[Santamaria(2019)]{santamaria_towards_2019}
A.~Santamaria.
\newblock \emph{Towards a {Godement} {Calculus} for {Dinatural}
  {Transformations}}.
\newblock PhD thesis, University of Bath, Bath, July 2019.

\bibitem[Selinger(2010)]{selinger_survey_2010}
P.~Selinger.
\newblock A {Survey} of {Graphical} {Languages} for {Monoidal} {Categories}.
\newblock In B.~Coecke, editor, \emph{New {Structures} for {Physics}}, volume
  813 of \emph{Lecture {Notes} in {Physics}}, pages 289--355. Springer, Berlin,
  Heidelberg, 2010.
\newblock ISBN 978-3-642-12820-2 978-3-642-12821-9.
\newblock \doi{10.1007/978-3-642-12821-9_4}.

\bibitem[Simpson(1993)]{simpson_characterisation_1993}
A.~K. Simpson.
\newblock A characterisation of the least-fixed-point operator by dinaturality.
\newblock \emph{Theoretical Computer Science}, 118\penalty0 (2):\penalty0
  301--314, Sept. 1993.
\newblock ISSN 0304-3975.
\newblock \doi{10.1016/0304-3975(93)90112-7}.

\bibitem[Wadler(1989)]{wadler_theorems_1989}
P.~Wadler.
\newblock Theorems for {Free}!
\newblock In \emph{Proceedings of the {Fourth} {International} {Conference} on
  {Functional} {Programming} {Languages} and {Computer} {Architecture}}, {FPCA}
  '89, pages 347--359, New York, NY, USA, 1989. ACM.
\newblock ISBN 978-0-89791-328-7.
\newblock \doi{10.1145/99370.99404}.

\end{thebibliography}

\vfill

{\small\medskip\noindent \copyright 2021. This manuscript version is made available under the CC-BY-NC-ND 4.0 license \url{http://creativecommons.org/licenses/by-nc-nd/4.0/}.}

\medskip\noindent\includegraphics[scale=0.75]{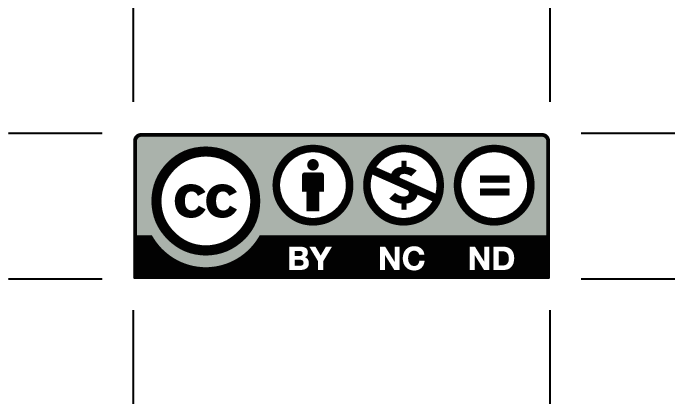}

\end{document}